\documentclass[11pt]{amsart}
\usepackage{graphicx}
\usepackage{amssymb}
\usepackage{epstopdf}
\usepackage{amsmath}
\usepackage{amscd}
\usepackage{amsthm,amsfonts,amsxtra,amssymb,amscd}
\usepackage[all]{xy}
\usepackage{enumerate}
\def\pd#1#2{\dfrac{\partial#1}{\partial#2}}

\usepackage{hyperref}

\theoremstyle{plain}
\newtheorem*{lemma*}{Lemma}
\newtheorem{lemma}[subsection]{Lemma}
\newtheorem*{theorem*}{Theorem}
\newtheorem{theorem}[subsection]{Theorem}
\newtheorem*{proposition*}{Proposition}
\newtheorem{proposition}[subsection]{Proposition}
\newtheorem*{corollary*}{Corollary}
\newtheorem{corollary}[subsection]{Corollary}
\theoremstyle{definition}
\newtheorem*{definition*}{Definition}
\newtheorem{definition}[subsection]{Definition}
\newtheorem*{example*}{Example}
\newtheorem{example}[subsection]{Example}
\theoremstyle{remark}
\newtheorem*{remark*}{Remark}
\newtheorem{remark}[subsection]{Remark}

\title{The algebra of the Box--spline}
\author{C. De Concini}\address{Dip. Mat. Castelnuovo, Univ. di Roma La
Sapienza, Rome, Italy}\email{deconcin@mat.uniroma1.it}
\author{C. Procesi}\address{Dip. Mat. Castelnuovo, Univ. di Roma La
Sapienza, Rome, Italy}\email{procesi@mat.uniroma1.it}
\thanks{The authors are partially supported by the Cofin 40
\%, MIUR}

\begin{document}
\begin{abstract} In this paper we want to revisit results of Dahmen and Micchelli, \cite{DM1},\cite{DM2}, \cite{DM3},  which we reinterpret and make more precise.  We compare these ideas with the work of Brion, Szenes, Vergne and others  \cite{BV}, \cite{SV}, \cite{Ve}, \cite{JKi}, \cite{S1}.\end{abstract}
\maketitle

\tableofcontents

\section{Introduction}

The main purpose of this paper is to revisit a celebrated   Theorem of Dahmen and Micchelli \cite{DM3} which we will state and prove again in a
somewhat stronger form.\smallskip

The theorem is on the following settings:\smallskip

Start from a finite list  $X:=\{a_1,\dots,a_N\}$   of non zero vectors $a_i\in\mathbb R^s$.\smallskip

If $X$ spans $\mathbb R^s$, from $X$ one builds an important function for numerical analysis, the {\it box spline}  implicitly  defined by the formula:
\begin{equation}\label{box} \int_{\mathbb R^s}f(x)B_X(x)dx :=\int_0^1\dots\int_0^1f(\sum_{i=1}^Nt_ia_i)dt_1\dots dt_N,\end{equation}
where $f(x)$ is any continuous function.\medskip

If  0 is not in the  convex hull of the vectors $a_i$  then one has  a simpler    function $T_X(x )$,  the {\it multivariate spline cf
\cite{dhr}} characterized by the formula:
\begin{equation}\label{multiva} \int_{\mathbb R^s}f(x)T_X(x)dx=  \int_{\mathbb R_+^N}  f(\sum_{i=1}^Nt_i a_i)dt,
\end{equation}
where $f(x)$  has compact support.

Both $B_X$ and $T_X$ have a simple geometric interpretation as functions computing the volume of certain variable polytopes.\medskip

If furthermore the vectors $a_i$ happen to be integral  vectors (and 0 is not in their convex hull) one has  a third function, now on $\mathbb
Z^n$,  important  for combinatorics, the {\it partition function} given by:
\begin{equation}\label{partf}\mathcal P_X(v):=\#\{(n_1,\dots,n_N)\,|\, \sum n_ia_i=v,\ n_i\in\mathbb N\}.\end{equation}

One of the goals of the theory is to give computable closed formulas for all these functions and, at the same time, describe some of their
qualitative behavior.\smallskip

The main result  of the theory is that,  these three functions can be described in a combinatorial way as a  finite sum over {\it local pieces}
(see formulas \ref{JK}  and \ref{JK2} ).   In the case of $B_X(x)$ and $T_X(x)$ the local pieces span, together with their derivative,   a
finite dimensional space $D(X)$ of polynomials. In the case of $\mathcal P_X(v)$ they span together with their translates,  a finite dimensional
space $\nabla(X) $ of quasi polynomials (cf. Definition \ref{PE}). \smallskip

The theorem we are referring to, characterizes:
\begin{itemize}
\item  $D(X)$ by differential equations.

\item $\nabla(X) $ by difference equations.
\end{itemize}
In particular Dahmen and Micchelli compute the dimension of both spaces. This dimension  has a simple combinatorial interpretation in terms of
$X$. They also decompose $\nabla(X) $  as a direct sum of natural spaces associated to certain special points $P(X)$ in a suitable torus. $D(X)$
is  then the space associated to the identity.\bigskip

Although this theorem originates from the theory of the box spline, nevertheless it has also an interest in commutative algebra and algebraic
geometry, in particular in the theory of hyperplane arrangements and partition functions.

Here we have been inspired by the  results of Orlik--Solomon on cohomology and the results of Brion, Szenes, Vergne  on partition
functions.\smallskip

In fact a lot of work originated from the seminal paper of Khovanski\u\i,  Pukhlikov \cite{KP}, interpreting the counting formulas as
Riemann--Roch formulas  for toric varieties and of Jeffrey--Kirwan, \cite{JKi}    on moment maps.  These topics are beyond this paper and we
refer  to Vergne's survey article \cite{Ve}.\smallskip

In the theory of hyperplane arrangements, $X$ is viewed as a set of linear equations and the object of study is the arrangement  of hyperplanes
defined by these equations.\medskip

Due to the somewhat  high distance between these two fields, people working in hyperplane arrangements do not seem to be fully aware of the
results on the Box--spline.
\smallskip

On the other hand there are some methods which have been developed in this latter theory which we believe shed some light on the space of
functions used to build the box spline.  Therefore we feel that this presentation may be useful to make a bridge between the two theories. Thus
this paper has been organized partly as a research paper and partly as a survey of some relevant aspects of these two theories.\bigskip

We divide  the discussion in three parts.

In the first on the {\it differentiable case} the main new results are determination of the graded dimension of the space $D(X)$( Theorem
\ref{main}) in terms of the combinatorics of bases extracted from $X$. An algorithmic characterization in terms of differential equations of a
natural basis of the top degree part of $D(X)$  (Proposition \ref{dual}) from which one obtains explicit local expressions for $T_X$ (Theorem
\ref{exxpr}). A duality between $D(X)$ and a subspace of the space of polar parts relative hyperplane arrangement associated to $X$ (Theorem
\ref{main1}). In this section we also give a simple proof of the Dahmen-Micchelli theorem on the dimension on $D(X)$ using elementary
commutative algebra (Theorem \ref{main2}).

The second on the {\it discrete case}, contains various extensions of the results of the first part in the case in which the elements in $X$ lie
in a lattice.  We develop a general approach to linear difference equations and give a method of reduction to the differentiable case (Section
\ref{difequ}). We also give an explicit formula relating partition functions to multivariate splines (Theorem \ref{toddo}). We end with a quick
survey of some of the standard applications as for example can be found in the book {\it Box splines} \cite{dhr}.

In the third part we explain the approach via residues (see Theorem \ref{corrare} and Formula \ref{resper}). We finish
 giving an overview of  relations with wonderful embeddings.
\medskip

\bigskip
{\bf  Warning}:  the reader will see that we  use Fourier transforms as an essentially algebraic tool. \smallskip

The fact is that, for the purpose of this theory, Fourier transform is essentially a duality  between polynomials and differential operators
with constant coefficients.\smallskip

As long as one has to perform only algebraic manipulations  one can avoid difficult issues of convergence and various forms of algebraic or
geometric (residues)  duality  is sufficient.

Thus we usually work with Laplace transforms  which avoids cluttering the notations with unnecessary $i$'s.\medskip

Our conventions for Laplace transforms are the following. We fix a vector space $V$ and set $U:=V^*$. We fix a Eucliden structure on $V$ which
induces Lebesgue measures $dv,du$ on $V,U$ and all their linear subspaces. We set
$$Lf(u):=\int_Ve^{-\langle u\,|\,v\rangle}f(v)dv.$$
$L$ maps functions on $V$ to functions on $U$. We have the basic properties, when  $p\in U, w\in V$, writing $p$, $D_w$   for the linear
function ${ \langle p\,|\,v\rangle}$ and the directional derivative on $V$, (and dually   on $U$):
\begin{equation}\label{L1} L(D_w f)(u)= wLf(u),\  L(p f)(u)=-D_pLf(u),\      \end{equation}
\begin{equation} \label{L2} L(e^pf)(u)= Lf(u-p),\   L(f(v+w))(u)=e^wLf(u).  \end{equation}
\bigskip

\part{The differentiable case.}
\bigskip

\section{Basic definitions}

It is convenient to take a somewhat intrinsic and base free approach to our problems.

Let us fix an $s-$dimensional vector space $U$, let us denote by $V$ its dual and fix a list $X:=\{a_1,\dots,a_N\}$ of non zero elements in $V$
(we allow repetitions in the list as this is important for the applications). We identify  the symmetric algebra $S[V]$ with the ring of
polynomial functions on $U$  and sometimes denote it  by  $A$.

This algebra can also be viewed as the  algebra of polynomial differential operators with constant coefficients on $V$. Similarly $S[U]$ is the
ring of polynomial functions on $V$, or polynomial differential operators with constant coefficients on $U$.\medskip

Given a vector $v\in V$ we denote by $D_v$  the corresponding  directional derivative. This is algebraically characterized,
on $S[U]$, as the derivation which on an element $\phi\in U$ takes the value  $\langle\phi\,|\,v\rangle$.\medskip

One can organize all these facts in the algebraic language of Fourier transform.  Let $W(V),  W(U)$  denote the two algebras of differential
operators with polynomial  coefficients on $V$ and $U$ respectively. Notice  that, from a purely algebraic point of view they are both
generated by $V\oplus U$.

In the first case $V$ is thought of as the space of directional derivatives and then we write $D_v$ instead of $v$, and $U$ as the linear
functions. In $W(U)$   the two roles are exchanged.

The relevant commutation relations are thus:
$$  [D_v,\phi]=\langle \phi\,|\, v\rangle,\qquad        [D_\phi,v]=\langle v |\phi \rangle .         $$

Thus we see that  we have a canonical isomorphism of algebras:
$$\mathcal F:  W(V)\to W(U),\quad D_v\mapsto -v,\quad \phi\mapsto D_\phi.$$

One usually writes $\hat a$ instead of $\mathcal F(a)$.\smallskip

This allows us, given a module $M$ over $W(V),$ to consider its {\it Fourier transform}  $\hat M$  as a module over $W(U)$  by $a.m:=\hat a m$
and conversely.

\subsection{Cocircuits}
We come to the first basic definition  of combinatorial nature.  The importance of this notion will be clear once we start to study the
multivariate spline.  We assume that $X$ spans $V$.
\begin{definition}
We say that a sublist $Y\subset X$ is a {\it cocircuit}, if the elements in $X-Y$ do not span $V$.\end{definition}

The minimal cocircuits can thus be obtained as follows:

Fix a hyperplane $H\subset V$  spanned by elements in $X$  and consider $Y:=\{x\in X\,|\, x\notin  H\}$, it is immediately verified that {\it
this is a cocircuit and every cocircuit contains one of this type}. \smallskip

Sometimes we shall express the fact, that a sublist $Z\subset X$ consists  of all the vectors in $X$ lying in a given  subspace, by saying that
{\it $Z$ is complete}.

Thus a minimal cocircuit is obtained by removing from $X$ a complete set spanning a hyperplane.\smallskip

The set  of all cocircuits will be denoted by $\mathcal E(X)$.
\medskip

\subsection{No broken circuits}  The second basic combinatorial notion has been used extensively  in the theory of hyperplane  arrangements (cf. \cite{Br},\cite{W1},\cite{W2},\cite{Wi}). \smallskip

Let $\underline c:=a_{i_1},\dots, a_{i_k}\in X, \ i_1<i_2\dots <i_k$,   be a sublist of linearly independent
elements.

\begin{definition} We say that $a_i$ {\it breaks} $\underline c$ if there is an index $1\leq e\leq k$ such that:
\begin{itemize}
\item $i\leq i_e$. \item $a_i$ is linearly dependent  on \quad  $a_{i_e},\dots, a_{i_k}.$
\end{itemize}
\end{definition}

In particular, given any basis  $\underline b:=a_{i_1},\dots, a_{i_s}$ extracted from $X$, we set:

$B(\underline b):=\{a\in X\,|\, a\quad\text{breaks}\quad \underline b\}$  and   $n(\underline b)=|B(\underline b)|$ the cardinality of
$B(\underline b)$.

\begin{definition} We say that $\underline b$ is {\it no broken}  if $B(\underline b)=\underline b$ or $n(\underline b)=s$.\end{definition}\smallskip

Let us denote by $\mathcal B(X)$ the set of all bases extracted from  $X$. We shall consider the map $\underline b\mapsto n(\underline b)$ as a
{\it statistic} on $\mathcal B(X)$.

\subsection{The box spline}     Let us recall some  points which are standard using the form presented in \cite{WV} or \cite{dp1}.\smallskip

First let us recall some basic facts on splines (cf. \cite{dhr}).\smallskip

Let  $C(X):=\{\sum_{a\in X}t_aa\,|\, 0\leq t_a,\ \forall a\}$  be the  cone of  linear combinations of vectors in $X$ with positive
coefficients.

We will assume   that 0 is not in the convex hull of the vectors in $X$,  i.e. that  $C(X)$  does not contain lines\smallskip.

We  have already defined, in the introduction, the two basic functions on $V$,  $B_X$  (formula (\ref{box})), and $T_X$,  (formula
(\ref{multiva})).
\smallskip

It is best to think of both   $T_X $  and $B_X $  as {\it tempered distributions} (cf. \cite{Yo}).  Then the definition is valid also if $X$
does not span $V$.

$B_X $ is supported   in the {\it box} $$B(X):=\sum_{i=1}^Nt_ia_i,\ 0\leq t_i\leq 1,\qquad\text{the {\it shadow} of the cube }[0,1]^N,$$
generated by $X$, $T_X$ is supported in $C(X)$.

{\bf Basic example} Let $X=\{a_1,\dots,a_s\}$ be a basis, $d:=|\det(a_1,\dots,a_s)|$:

$B(X)$ is the parallelepiped   with edges the $a_i$,   $C(X)$ is the positive quadrant generated by $X$.

\begin{equation}\label{basic} B_X=d^{-1}\chi_{B(X)},\qquad T_X=d^{-1}\chi_{C(X)}\end{equation}   where, for any given  set $A$, we denote by  $\chi_{A}$   its characteristic function.
\smallskip

If $X=\{a_1,\dots,a_k\}, \ k<s$ is only a linearly independent set, we have to consider   $T_X $  and $B_X $ as meausers on the subspace spanned
by $X$.

$T_X $  and $B_X $ are   functions as soon as $X$  spans  $V$, i.e. when the support pf the distribution has maximal dimension.\medskip

These functions have a nice geometric interpretation.

Let $F:\mathbb R^N\to V$ be defined by $F(t_1,\dots,t_N):=\sum_{i=1}^Nt_ia_i$. Then $B_X(w)$ is the volume of the    polytope $F^{-1}(w)\cap
[0,1]^N$  while $T_X(w)$ is the volume of the    polytope $F^{-1}(w)\cap [0,\infty]^N$ (with a suitable normalization constant).\medskip

It is useful to generalize these notions, introducing  a {\it parametric version} called $E-$splines.
\smallskip

Fix parameters $\underline \mu:=\{\mu_1,\dots,\mu_N\}$ and define the functions (or tempered distributions)  on $V$ by the implicit formulas:

\begin{equation}\label{boxm} \int_Vf(x)B_{X,\underline\mu}(x)dx :=\int_0^1\dots\int_0^1e^{-\sum_{i=1}^Nt_i\mu_i}f(\sum_{i=1}^Nt_ia_i)dt_1\dots dt_N.\end{equation}

\begin{equation}  \int_Vf(x)T_{X,\underline\mu}(x)dx=  \int_{\mathbb R_+^N}  e^{-\sum_{i=1}^Nt_i\mu_i}  f(\sum_{i=1}^Nt_i a_i)dt
\end{equation}

Also these functions have a nice geometric interpretation.

They represent the integral of $e^{-\sum_{i=1}^Nt_i\mu_i}$ on the    polytope $F^{-1}(w)\cap [0,1]^N$   or $F^{-1}(w)\cap [0,\infty]^N$ (with
the same   normalization constant).  Of course for $\underline \mu=0$ we recover the previous definitions.\medskip

An easy computation gives their  Laplace transforms:
\begin{equation}\int_V\!\!\!e^{-\langle x,y\rangle}B_{X,\underline\mu}(x)dy=\int_0^1\dots\int_0^1e^{-\sum_{i=1}^Nt_i(\langle x, a_i\rangle+\mu_i )}dt_1\dots dt_N\end{equation}$$ =\prod_{a\in X}{1- e^{-  a -\mu_a}\over  a  +\mu_a}.$$
and
\begin{equation}\label{laplaT}\int_V\!\!\!e^{-\langle x,y\rangle}T_{X,\underline\mu}(x)dy=\int_0^\infty\dots\int_0^\infty e^{-\sum_{i=1}^Nt_i(\langle x, a_i\rangle+\mu_i )}dt_1\dots dt_N\end{equation}$$ =\prod_{a\in X}{1\over   a  +\mu_a}.$$

We have written shortly $a:= \langle x, a \rangle$, for the linear function on $U$.
\medskip

The use of Laplace rather than  Fourier transforms is justified by the following discussion.

Define the {\sl dual  cone} $ \widehat{C(X)}$ of $C(X)$.
$$\widehat{C(X)}:=\{u\in U\,|\, \langle u\,|\, v\rangle \geq 0,\ \forall v\in C(X)\}. $$ This cone    consists thus of the linear forms that are non negative on $C(X)$. Its interior in not empty since   $C(X)$ contains no lines.
\begin{proposition}
If $T$ is a tempered distribution supported in $C(X)$  its Fourier transform    is an analytic function, of a complex variable $z=x+iy,$ $ x\in
\widehat{C(X)},y \in U$, on the   open set where the real part $x$ lies in the interior $ \widehat{C(X)^0}$ of  $ \widehat{C(X)}$.
\end{proposition}
\proof    This depends on the fact that, if $u\in \widehat{C(X)^0} $ we have that $e^{- \langle u\,|\, v\rangle}$ has exponential decay on
$C(X)$.\qed

In fact for the $T_X$ and all the distributions that we shall encounter  we will have that  their Fourier transform are not only defined in the
region where the real part $x$ lies in   $ \widehat{C(X)^0} $, but in fact they extend to meromorphic functions with poles on  the hyperplanes
$a_i=0$  (or sometimes translates of these hyperplanes).\medskip

In the course of this paper we give an idea of the general algebraic calculus involving these distributions. Under Laplace transform one can
reinterpret the calculus in terms of the structure of certain algebras of rational functions (or exponentials) as $D-$modules.\smallskip

Given $a\in V,\ \mu\in\mathbb C$ let us introduce the notation, which will be discussed more deeply in {\bf Part 2}:
$$\nabla_a^\mu f(x):=f(x)-e^{-\mu}f(x-a).$$

From the expressions of the Laplace transforms one gets that,  the box spline can be obtained from the multivariate spline by a simple
combinatorial formula.

\begin{proposition}
For every subset $S\subset X$ we set $a_S:=\sum_{a\in S}a,$ and $ \mu_S:=\sum_{a\in S}\mu_a$ then:

\begin{equation}\label{boxc} B_{X,\underline \mu}(x)=\prod_{a\in X}\nabla_a^{\mu_a} T_{X,\underline \mu}(x)=\sum_{S\subset X}(-1)^{|S|} e^{-\mu_S}T_{X,\underline \mu}(x-a_S). \end{equation}
\end{proposition} \smallskip

\begin{proof} It follows from the basic rule (\ref{L2})
which gives the commutation relation between the Laplace transform and translations.
\end{proof}

\subsection{{The space $D(X)$.}}
The following definition is of central importance in the work of Dahmen and Micchelli and in this paper
\begin{definition}
We define the space $D(X)$ by the condition:
\begin{equation}\label{ilsist} D(X):=\{p\,|\, D_Yp=0,\ \forall Y\in \mathcal E(X),\quad\text{ the cocircuits} \}.\end{equation}
\end{definition}
In this definition we shall assume that $p$ is a polynomial. In fact, due to the property that  the ideal generated by the elements $D_Y$
contains all large enough products of derivatives (see \ref{span}), one can easily see by induction that any distribution $p$ satisfying  (\ref{ilsist}), is
already a polynomial  (cf. \cite{FKT}). \smallskip

We will see later a generalization with parameters $\underline \mu$  of these equations and also a discrete analogue, using difference operators
rather than derivatives.
\smallskip

Let $m$ be the minimum number of elements in a cocircuit in $X$, assume   $m\geq 2$,  $m$ is characterized by the property that the basic space
$D(X)$  contains {\it all} polynomials of degree $<m$. We shall see in Part 4,  that $m$ controls also the smoothness of the splines that we
have introduced.

\subsection{Two basic modules in correspondence\label{duemod}}

The theory of the Laplace transform tells us how to transform  some basic manipulations on distributions as an algebraic calculus. In our
setting this is best seen introducing the following two $D-$modules  in Fourier duality:\smallskip

The first is the $D-$module $\mathcal D_X:=W(V)T_X$  generated, in the space of tempered distributions, by $T_X$ under the action of the
algebra $W(V)$ of differential operators on $V$ with polynomial coefficients (for basic facts cf.  \cite{Yo}).\smallskip

The second  $D-$module is the  algebra $R_{X} :=S[V][\prod_{a\in X}    a  ^{-1}]$  obtained from the polynomials on $U$ by inverting the element
$d_{X }:= \prod_{a\in X}    a   .$ This is    a module under $W(U)$ and it is the coordinate ring of the open set $\mathcal A_X$  complement of
the union of the hyperplanes of $U$ of equations $a=0,\ a\in X$.\medskip

\begin{theorem}\label{Lapp}   Under Laplace transform,  $\mathcal D_X$    is mapped isomorphically  onto $R_{X}$.
In other words we get a canonical isomorphism of $\hat {\mathcal D}_X$ with  $R_{X}$ as $W(U)$-modules.

$\mathcal D_X$ is the space of tempered distributions which are linear combinations of polynomial functions on the cones  $C(A)$, $A\subset X$ a
linearly independent subset and their distributional derivatives.
\end{theorem}
\begin{proof} The injectivity of the Laplace transform on $\mathcal D_X$ is a standard fact \cite{Yo}.

To see the surjectivity, notice that under the  Laplace transform, by definition and formula (\ref {laplaT}), the image of $\mathcal D_X$ under
the Laplace transform is the smallest $D$-module containing $d_X^{-1}$. Since $R_X$ contains  $d_X^{-1}$ it suffices to see that $d_X^{-1}$
generates $R_X$ as a $D$-module. To see this first notice that if we take  linearly dependent vectors   $a_0=\sum_{i=1}^k\alpha_ia_i,$ in $X$,
we can   write:
$$   {1\over \prod_{i=0}^k a_i  }=   {a_0  \over a_0^2\prod_{i=1}^k a_i  }=\sum_{i=1}^k\alpha_i{ 1 \over a_0^2\prod_{j\neq i. j=1}^k a_i }$$
So, repeating this algorithm, we see that we can write a fraction   ${1/( \prod_ia_i^{h_i})}$  as a linear combination of fractions where the
elements appearing in the denominator are linearly independent. In particular $R_X$ is spanned by the functions of the form $f/(
\prod_ib_i^{k_i+1})$, with $f$ a polynomial, $B=\{b_1,\ldots b_t\}\subset X$  a linearly independent subset and $k_i\geq 0$ for each $i=1,\ldots
,t$.

Now it is  clear that $f/( \prod_ib_i^{k_i+1})$  lies in the $D$-module generated by $1/( \prod_ib_i)$ and
$${1\over \prod_ib_i}={\prod_{a\in X-B}a \over d_X}$$
So, our first claim follows.

The second  follows from our previous discussion once we remark that if $B=\{b_1,\ldots b_t\}\subset X$  is linearly independent, $1/(
\prod_ib_i)$ is the  Laplace  transform of the measure on $V$ which is the push forward of the measure on the subspace spanned by the vectors in
$B$ given by $\chi_{C(B)}\mu_B$ where $\chi_{C(B)}$ is the characteristic function of the cone $C(B)$ and $\mu_B$ is the Lebesque measure
normalized so that the parallelepiped   with edges the $b_i$ has volume 1.  \end{proof}

Let us now introduce a filtration in $R_{X} $   by  $D-$submodules which we will call the {\it filtration by polar order}.

This is defined algebraically as follows. One puts in filtration degree $\leq k$ all the fractions $f\prod_{a\in X}    a  ^{-h_a}, \ h_a\geq 0$
for which the set of vectors $a$,  with $h_a>0$, spans a space of dimension $\leq k$. Denote this part of the filtration by $R_{X,k} $ . Notice
that by our proof of Theorem \ref{Lapp}, $R_{X,s}=R_X. $

From the proof of Theorem \ref{Lapp} one deduces that the corresponding (under inverse Laplace transform) filtration on $\mathcal D_X$ can be
described geometrically as follows.

We cover $C(X)$ by the positive cones $C(A)$ spanned by linearly independent subsets of $A\subset X$.  We define $C(X)_k$ to be the
$k-$dimensional skeleton of the induced stratification, a union of $k-$dimensional cones.

\begin{proposition}\label{supporto} $\mathcal D_{X,k}$ consists of the tempered distributions in $\mathcal D_X$ whose support is contained in $C(X)_k$.\end{proposition}

\subsection{Polar parts\label{polari}}  Proposition \ref{supporto}  implies that, the space  $\mathcal D_{X,s-1}$ is formed by all the distributions in $\mathcal D_X$ which vanish once computed on test functions with support in the {\it set of regular points}  $C(X)^0:=C(X)-C(X)_{s-1}$.

In other words we may identify   $\mathcal D_X/\mathcal D_{X,s-1}$ with a space of distributions on the open set $C(X)^0$.

By duality we then see that the top subquotient, $R_X/R_{X,s-1}$, has a special importance in the theory \cite{BV1}. Let us define it in a much
more general setting which will be useful later.
\smallskip

Let $M$ be an analytic manifold. $Z=Z_1\cup Z_2\cup\cdots \cup Z_m$ be a divisor union of the smooth irreducible components $Z_i$, $i=1,\ldots
,m$. $p$ an isolated point in $Z_1\cap Z_2\cap\cdots \cap Z_m$ with the property that there exist local coordinates around $p$ so that the $Z_i$
are locally given each by a linear equation $a_i=0$.  A subset of the set of divisors $Z_i$ is said to be {\sl non spanning} in $p$, if $p$ is
not an isolated point of their intersection. We then define:
\begin{definition} Define the module of polar parts in $p,$ with respect to $Z,$ and denote it by
$P_{Z,p}$ (or $P_Z$ if $p$ is clear from the context) as the   quotient of the space of germs in $p$ of  meromorphic functions with poles in $Z$
modulo those function whose polar part is supported in a non spanning subset.
\end{definition}
$P_{Z,p}$ is clearly a module over the ring  $\mathcal Z_p$ of germs of differential operators with holomorphic coefficients around $p$. In  the
local coordinates $x_1,\ldots ,x_s$ around $p$, where $Z_i$ has an equation $a_i=0$, $a_i$ i  linear  this space can then be identified with the
module  $P_X:=R_{X}/R_{X,s-1}$ with $X=\{a_1,\ldots ,a_m\}$. Notice the under this identification the operators induced from $\mathcal Z_p$
coincide with the algebra of operators $W(U)$ where $U=T_pX$.
\medskip

\subsection{Basic modules} 
All the modules over Weyl algebras which appear are built out of some basic irreducible modules, in the sense that they have finite composition
series in which  these basic modules appear. It is thus useful to quickly discuss these basic modules.  Let us take $W(V)$, the differential
operators on $V$.  The most basic module on $W(V)$ is the polynomial ring  $S[U]$.  It is the cyclic module   generated by 1 and the annihilator
ideal of  1 is generated by all the derivatives $D_v,\ v\in V$.   Its Fourier transform can be identified with the module   of distributions
supported in $0\in U$, a cyclic module generated by the Dirac distribution $\delta_0$.

Given any point $p\in U$, consider   the 1-dimensional $S[V]$ module $\mathbb C_p$  given by evaluation at $p$ and the induced $W(U)$ module
$N_p:=W(U)\otimes_{S[V]}\mathbb C_p$. $N_p$ is easily seen to  be irreducible and   free of rank 1 as $S[U]-$module.

In fact, in the language of distributions  $N_p$ is identified to the $W(V)$ submodule generated by the Dirac distribution $\delta_p$.

We shall need the following (easy) and standard fact (\cite{Co}).
\begin{lemma}\label{liuno}
1) Given a $W(U)-$module  $M$ and a nonzero element  $u\in M$, if $fv=f(p)v, \forall f\in S[V]$, then $u$ generates a submodule isomorphic to
$N_p$.

2) Given linearly independent vectors  $u_i$ satisfying the previous hypotheses, the submodules that they generate form a direct sum.

3)  The module $N_p$ has as characteristic variety the cotangent space at $p$.

\end{lemma}
The only use we make of part 3) is its consequence that, for distinct points $p,q$, the corresponding modules $N_p,N_q$ are not
isomorphic.\smallskip

Also for  a linear subspace $A\subset V$ or a translate $A+v$ one can define the irreducible module $N_A$ generated by the Dirac distribution $\delta_A$ given by:
$$\langle \delta_A \,|\, f\rangle:=\int_Af(w)dw$$

The annihilator of   $\delta_A$  is generated by the elements $u\in A^\perp\subset U$ vanishing on $A$  (for  $A+v$  we have the elements
$u-u(v),\ u\in A^\perp$) and the elements $D_x,\ x\in A$.   the fact that $N_A$ is irreducible and the previous elements generate the
annihilator of $\delta_A$ can either be verified directly or by remarking that  $N_A$  can be obtained by twisting the polynomial ring  by an
automorphism (a partial Fourier transform) of $W(V)$  defined as follows. One decomposes  $V=A\oplus Z$ for some complement $Z$, then one has
that $W(V)=W(A)\otimes W(Z)$ and apply Fourier transform to the factor  $W(Z)$.\medskip

In explicit coordinates $x_i$ where $A:=\{x_i=0,\ i\leq
k\}$ we have that the partial Fourier transform is the identity on the variables $x_j$ and their derivatives when $j>k$, while
$$x_i\mapsto\pd{}{x_i},\ \pd{}{x_i} \mapsto -x_i,\ \forall i\leq k.$$

As the polynomial ring is a free rank 1 module over the polynomials we see that $N_A$ is thus a free rank 1 module over the partial Fourier
transform of the polynomials, that is the polynomial algebra   $P_A$ in the variables $ \pd{}{x_i} ,\  \forall i\leq k,\quad x_i \ \forall i>k.$

The module   $N_A$ appears naturally in our theory as follows.  Take a basis   $\underline c:=\{a_1,\dots,a_k\}$  of $A$ and consider the  open  cone  $C(\underline c)$ generated by $\underline c$.  Consider next the distribution  $\delta_{C(\underline c)}$ given by integration on this cone under a translation invariant Lebesgue measure for which the parallelepiped  generated by $\{a_1,\dots,a_k\}$  has volume 1.   Then one easily sees that the Laplace transform  of $\delta_{C(\underline c)}$ is $\prod_{i=1}^ka_i^{-1}$ and $N_A$  appears as the module of polar parts for the   polynomial ring with the $a_i$ inverted.

In particular fix some set of linear coordinates $a_1,\dots, a_s$ for $U$.
If for every subset   $S\subset \{1,\dots,s\}$ we set $A_S$ to be the subspace where the variables $a_i=0,\ \forall i\in S$  we see, with the previous notations, that the ring of Laurent polynomials   decomposes as a direct sum  $$\mathbb C[a_1^{\pm 1},\dots, a_s^{\pm 1}]=\oplus_{S\subset \{1,\dots,s\}}  P_{A_S}({1\over \prod_{i\in S}a_i}).$$
Notice   that   $P_{A_S}({1\over \prod_{i\in S}a_i})$ has as basis, the monomials $\prod_ia_i^{k_i}$   with $k_i<0,\ \forall i\in S$ and $k_i\geq 0,\ \forall i\notin S$.

 This decomposition, together with the description of  $A_S$ as a polynomial algebra, gives the explicit partial fraction decomposition for Laurent polynomials, which is in any case elementary.

 Following the same ideas one can develop the general partial fraction decomposition in the case in which we invert any set of linear equations as we shall see in   Proposition \ref{mulp} and Remark \ref{ilrema}.\bigskip

\section{The function $T_X$}

\subsection{Local expansion \label{31}  }  Let us extend the previous ideas to the parametric case.
We take a list $\underline \mu:=\{\mu_a|a\in X\}$ of complex numbers. We introduce the ring $R_{X,\underline \mu}:=S[V][\prod_{a\in
X}(a+\mu_a)^{-1}]$. It is clear that we can introduce on this algebra a filtration completely analogous to that of $R_X=R_{X,0}$.

To a basis  $\underline b:=\{b_1,\dots, b_s\}$,   from $X$, we associate the unique point $p_{\underline b}\in U$ such that $b_i(p_{\underline
b})=-\mu_{b_i}$ for each $i=1,\ldots ,s$. The set $P(X,\underline \mu)$ consisting of the points $p_{\underline b}$ as $\underline b$ varies
among the bases extracted from $X$ is called the set of {\it points of the arrangement}.

For generic  $\underline \mu$  all these points are distinct, while for $\underline \mu=0$  they all coincide with 0. In the other cases we may
have   various ways in  which the points   $P_{\underline b}$ will coincide, depending on the compatibility relations among the parameters
$\underline \mu$.

Given   $p\in P(X,\underline \mu)$ we define the subset  $$X_p:=\{a\in X\,|\, a(p)+\mu_a=0\}\quad (X_p=\cup_{p_{\underline b}=p}\{\underline
b\}.)$$ of the elements in $X$ such that the affine hyperplane $H_a$ of equation $a(p)+\mu_a=0$ contains $p$. It is clear by definition that, if
we restrict to the subset $X_p$, the points of this restricted arrangement reduce to $p$.   Moreover a change of variables     $a':=a+\mu_a$
corresponding to a translation, centers the arrangement in 0.\smallskip

The divisor $Z:=\cup_{a\in X}  H_a$ satisfies the hypotheses of section \ref{polari}, thus we can construct, for each point $p\in P(X,\underline
\mu)$  the corresponding module of polar parts, which we denote by $P_{Z,p}$. In order to understand these modules we need a preliminary
construction.

The following Proposition allows us to reduce the computation of  $T_{X, \underline \mu }$ to that of the various $T_{X_p}$, for $p\in
P(X,\underline \mu)$. Let us denote, for a no broken basis $\underline b$ by $u_{\underline b}$   the class of the element
$\prod_{a\in\underline b}(a+\mu_a)^{-1}$  in the quotient $R_X/R_{X,s-1}$.  Since, given $a_0\in \underline b$ we clearly have that
$(a_0+\mu_{a_0})\prod_{a\in\underline b}(a+\mu_a)^{-1}\in R_{X,s-1}$ we verify easily that:
$$  f u_{\underline b}=f(p)u_{\underline b},\ \forall f\in   S[V]. $$
It follows that, unless $u_{\underline b}=0$ (which is not the case), $u_{\underline b}$ generates a submodule isomorphic to $N_p$. In fact one has:
\begin {proposition}\label{mulp} 1) The mapping $R_X/R_{X,s-1}\to \oplus_{p\in P(X,\underline\mu)}P_{Z,p}$ is an isomorphism of $W(U)$ modules.

2) Each $P_{Z,p}$ is an isotypic component.

3) Each $P_{Z,p}$ is the direct sum of as many copies of $N_p$ as the no broken bases $\underline b$ in $X_p$, each generated by a corresponding
element $u_{\underline b}$.

4) $T_{X, \underline \mu }=\sum_{p\in P(X,\underline \mu)}  c_pe^{ p} T_{X_p},$ with $c_p$ some explicitly computable constants.
\end{proposition}
\begin{proof}
We can use a slightly more precise reduction. Assume we have linearly dependent vectors   $a_0=\sum_{i=1}^k\alpha_ia_i,$ in $X$.   If
$\nu:=\mu_{a_0}-\sum_{i=1}^k\alpha_i\mu_{a_i}\neq 0$ we write:
$$   {1\over \prod_{i=0}^k(a_i+\mu_{a_i})}=\nu^{-1}  {a_0+\mu_{a_0}- \sum_{i=1}^k\alpha_i(a_i+\mu_{a_i})\over \prod_{i=0}^k(a_i+\mu_{a_i})}$$
and then, develop into a sum of $k+1$ terms in each of which one of the elements $a_i+\mu_{a_i}$ has disappeared.  This allows us to separate
the denominators with respect to the points in $P(X,\underline\mu)$.
\smallskip

If $\nu=0$ we can write $$   {1\over \prod_{i=0}^k(a_i+\mu_{a_i})}=   {  \sum_{i=1}^k\alpha_i(a_i+\mu_{a_i})\over
(a_0+\mu_{a_0})\prod_{i=0}^k(a_i+\mu_{a_i})}.$$

We expand and then simplify the numerators with the denominators and obtain that the element $   \prod_{i=0}^k(a_i+\mu_{a_i})^{-1} $   can be
expanded as a linear combination of elements of type  $  \prod_{a_i\in\underline b} (a_i+\mu_{a_i})^{-h_i}, h_i\geq 0 $ with $\underline b$ a
basis.  Only the elements with all the $h_i>0$ give non zero classes in the module of polar parts, moreover when $h_i>0$  from the rules of
derivatives  $\prod_{a_i\in\underline b} (a_i+\mu_{a_i})^{-h_i} $ is obtained by applying a monomial in the derivatives to  $
\prod_{a_i\in\underline b} (a_i+\mu_{a_i})^{-1}$. This shows that the map is onto.  The modules $P_{Z,p}$ are contained  in $R_X/R_{X,s-1}$ and
belong to distinct isotypic components by the previous lemma on the characteristic variety, hence the direct sum of  $R_X/R_{X,s-1}$.

A similar discussion allows us to  reduce  the elements appearing in the denominator to no broken sets. Thus the final point to verify, applying
part 2) of Lemma \ref{liuno},  is that the elements $u_{\underline b}$   are linearly independent in $P_{Z,p}$. This  we shall show as a
consequence of the theory of residues (Theorem \ref{indiplin}),
    completing the proof of 1) and 2).\smallskip

By (\ref {laplaT}) we see that $\prod_{a\in X_p}(a  -a(p))^{-1}$  is the Laplace transform of  $e^{p} T_{X_p}.$

Thus we are reduced to showing  that, for suitable constants $c_p$ we have:
\begin{equation}\label{separa}\prod_{a\in X}{1\over   a  +\mu_a}=\sum_{p\in  P(X,\underline \mu)}  c_p\prod_{a\in X_p}{1\over   a  +\mu_a}=\sum_{p\in  P(X,\underline \mu)}  c_p\prod_{a\in X_p}{1\over   a  -a(p)}.\end{equation}
This follows by induction applying the basic algorithm of separation of denominators.\end{proof}

\begin{remark}\label{sommandi} Notice that the summands $P_{Z,p}$ of $R_X/R_{X,s-1}$ are the spaces of generalized common eigenvalues for the commuting operators induced from $V$ (a vector $v\in V$ has eigenvalue $v(p)$). Thus any $S[V]$ submodule in $R_X/R_{X,s-1}$ decomposes canonically into the direct sum of its intersections with the various $P_{Z,p}$.
\end{remark}

 There is a similar  description of all the pieces $R_{X,k}/R_{X,k-1}$ based on the {\it spaces of the arrangement} generated by the hyperplanes $a=0,
 \ a\in X$ of codimension $k$.
\begin {proposition}\label{ilrema}$R_{X,k}/R_{X,k-1}$  is a direct sum of copies of Fourier transforms of the modules $N_W$, as $W$ runs over the
subspaces of the arrangement, and  for given $W$ the sum is made of terms indexed by the no broken bases in $\underline c\subset X\cap
W^\perp\subset W^\perp$ and generated by the class of $\prod_{a\in  \underline c}a^{-1}$.
\end{proposition}
\proof Consider, for each $W$,  the map
$ R_{X\cap
W^\perp,k}\subset  R_{X,k}\to R_{X,k}/R_{X,k-1}.$  By induction the characteristic variety of all irreducible factors in $R_{X,k-1} $ and  $R_{X\cap W^\perp,k-1}$ are union of conormal spaces to subspaces $A$   of dimension $\leq k-1$.  Since $ R_{X\cap
W^\perp,k}/R_{X\cap W^\perp,k-1}$ is a direct sum of irreducibles with characteristic variety  the conormal space to $W$, it follows that
this map induces an inclusion $ R_{X\cap
W^\perp,k}/R_{X\cap W^\perp,k-1} \subset   R_{X,k}/R_{X,k-1}.$  It is easy to verify that this map is
an isomorphism with the isotypic component of type $N_W$ in $R_{X,k}/R_{X,k-1},$   reducing to the previous case.\qed

\subsection{Local expression for $T_X$\label{loche}}
Formula (\ref{separa}) implies immediately that:
\begin{equation} T_{X, \underline \mu }(x)=\sum_{p\in P(X,\underline \mu)}c_p e^{ p}T_{X_p  } (x)
\end{equation}

This together with formula (\ref{boxc})   gives for
  the box spline:
$$ B_{X, \underline \mu }(x)=\sum_{p\in P(X,\underline \mu)}c_p\sum_{S\subset X-X_p}(-1)^{|S|}e^{-\mu_S+p}B_{X_p  }(x-a_S).$$

Of course this is also a reformulation of the identity of Laplace transforms:
$$\prod_{a\in X}{1- e^{-  a -\mu_a}\over  a  +\mu_a}= \sum_{p\in  P(X,\underline \mu)}  c_p\prod_{a\in X-X_p}(1- e^{-  a -\mu_a})\prod_{a\in X_p}{1- e^{-  a -\mu_a}\over   a  +\mu_a} $$

As a consequence of these formulas the essential problem is the determination of $T_X$ in the non parametric case.
We are  indeed ready to state and prove the main formula one can effectively use for computing the function $T_X$.

We first need some  geometry of the cone $C(X)$. Denote by  $C(\underline b)$ the positive quadrant spanned by the no broken circuit basis  $
\underline b $. The positive cones  $C(\underline b)$  cover the cone $C(X)$  and induce a decomposition of $C(X)$ into polyhedral cones which
happens to be independent of the order chosen and thus of the no broken circuits (this is proved in \cite{dp2}).  The points of $C(X)$, outside
the boundaries of these cones are called {\it regular} and they are a union of open convex cones called the {\it big cells}.\smallskip

The remaining {\it singular points} are the union of  cones $C(A)$  generated by subsets   $A\subset X$  which do not span $V$. We denote by
$\mathcal {NB}$  the set of no broken bases extracted from $X$.

\begin{theorem} \label{exxpr}  Given a point $x$ in the closure of a big cell $\mathfrak c$ we have
\begin{equation}\label{JK}  T_X(x )=\sum_{\underline b\ |\ \mathfrak c\subset C(\underline b)}|\det(\underline b)|^{-1}p_{\underline b,X}  ( -x).\end{equation}
where for each no broken basis $\underline b$, $  p_{\underline b,X}( x)$ is a uniquely defined homogeneous polynomial of degree $|X|-s$ lying
in $D(X)$.  \end{theorem}
\begin{proof} By the continuity of $T_X$ in $C(X)$, it is sufficient to show our   claim in the interior on each big cell. Thus, by formula (\ref{basic}) and our discussion above we need to show that the identity
(\ref{JK}) holds, with $  p_{\underline b,X}( x)\in D(X)$,  in the $D$-module $\mathcal D_{X}/\mathcal D_{X,s-1}$.
\medskip

If we work in the space of polar parts $P_X=R_X/R_{X,s-1}$ we have, by our previous discussion that there exist uniquely defined  polynomials
$p_{\underline b,X}(x)$ homogeneous of degree $|X|-s$ with
\begin{equation}\label{JK2}  {1\over d_X}=\sum_{\underline b \in\mathcal {NB}} p_{\underline b,X}( \partial_x) {1\over d_{  \underline b  }}\end{equation}
In fact this identity already holds in $R_X$.

Applying the inverse of the Laplace transform  (cf.  \ref{basic}, \ref{laplaT}) we obtain identity (\ref {JK}).

To finish notice that  $Y$ is a cocircuit  we have $\prod_{y\in Y}yd_X^{-1}\in R_{X,s-1}$. In other words, the polynomials $  p_{\underline
b,X}$ in formula (\ref{JK2})   lie in $D(X)$.

\end{proof}
From now on, unless there is a possibility of confusion, we shall write $p_{\underline b }$  instead of $p_{\underline b,X}$.

\begin {remark}  The polynomials    $p_{\underline b}( x)$, with $\underline b$   a no broken basis,  are characterized among the polynomials in $D(X)$,  by the further  differential equations
which will be described at a later section (see (\ref{dual})).  In particular we  will show that the polynomials $p_{\underline b}( x)$ are
linearly independent. \end{remark}
\section{The first theorem}

\subsection{The basic equations\label{prima}}

In section 2.4 we have introduced the elements $D_{ Y}:=\prod_{x\in  Y}D_x\in A,\ Y\subset X$, thought of as   differential  operators on $V$.

We have started to see that  the   space \begin{equation}\label{sol} D(X):=\{p\in S[U]\,|\,D_{ Y}p=0,\ \forall\, Y\in\mathcal
E(X)\},\end{equation}  of polynomials $p$ on $V$, which are  solutions of the differential equations, $D_{ Y}p=0$, for all the cocircuits $Y$,
plays a fundamental role in the determination of $T_X$.  In \cite{DM1} the authors prove that:

\begin{theorem}\label{ilteorema} Let  $D(X) \subset S[U]$ be the space  of polynomials on $V$,  solutions of the differential equations  $D_{ Y}f=0$, as $Y$ runs over the cocircuits in $X$.

$D(X)$ is finite dimensional, of dimension the number $d(X)$ of linear bases of $V$ which one can extract from $X$.\end{theorem}

Since the elements $D_{ Y}$ are homogeneous, the space $D(X)$  is also graded and it is interesting to compute its dimension in each degree
$D_k(X)$.

Moreover from general facts we shall see that $D(X)$ is generated, under taking derivatives, by the homogeneous elements of top  degree
$N-s$.\smallskip

As usual one can arrange these dimensions in a generating function, $H_X(q):=\sum_{k}\dim(D_k(X))q^k$. We will show   (Theorem  \ref{main}), by
exhibiting an explicit basis, that this polynomial is given by the  statistic introduced in \S 2.1, i.e.:
\begin{equation}\label{dimen} H_X(q):=\sum_{\underline b\in\mathcal B(X)}q^{N-n(\underline b)}.\end{equation}
In particular the top degree polynomials in $D(X)$ are of degree $N-s$. They appear in the formula of the multivariate spline and thus  have a
very interesting geometric interpretation, since they compute the volume of certain polytopes associated to $X$  (\cite{dp1}).

\medskip

\subsection{A remarkable family\label{remarkable}}

In order to obtain the  proof of \ref{ilteorema} and (\ref{dimen}), we start from the study of   a purely algebraic geometric object.

For notational simplicity let us denote by $A:=S[V]$.

We want to  describe   the  {\it scheme} defined by the ideal $I_X$ of $A $ generated by the elements $M_{Y}:=\prod_{x\in  Y}x $ as $Y$ runs
over all the cocircuits.

\smallskip

Thus we are interested in the  algebra:
\begin{equation}\label{alg}  A_X:=A/I_X.\end{equation}
We shall soon see that $A_X$ is  finite dimensional. Formally $A_X$  is the {\it coordinate ring} of the corresponding scheme.\smallskip

The use of the word {\it scheme} may seem a bit fancy. What we really want to stress by this word is that we have a finite dimensional algebra
(quotient of polynomials) which as functions vanish exactly on some point $p$  but at least infinitesimally they are {\it not constant}. This
appears clearly in the dual picture which produces solutions of differential equations.\medskip

To warm us up in the proof let us verify that this scheme is supported at 0 (which implies that $A_X$ is finite dimensional).

For this, remark that  the variety of zeros of a set of equations, each one of which is itself a product of equations, is the union of the
subvarieties defined by selecting an equation out of each product. Thus what we need is the following:
\begin{lemma}\label{span}  Take one element $x_i$ from each cocircuit $Y_i$, then the resulting  set of elements $x_i$ span $V$  (hence they define \{0\} as subvariety).

\end{lemma}\label{comp}
\proof  If, by contradiction, these elements do not span, their complement is a cocircuit. Since we selected an element from each cocircuit this
is not possible.

\qed

In  fact it is convenient to extend the notion as follows. If  $X=\{a_1,\dots,a_N\}$ and $\underline\mu:=\{\mu_1,\dots,\mu_N\}$ are parameters,
we can define in general the ideal $I_X(\underline\mu)$ given by the equations  $ \prod_{a_j\in   Y}(a_j-\mu_j),\ Y\in\mathcal E(X)$.  This can
either be viewed as an ideal in  $S[V]$ depending on the parameters $\underline\mu$ or as an ideal in the polynomial ring
$S[V][\mu_1,\dots,\mu_N].$

\smallskip

\subsection{The first main theorem}

Theorem (\ref{ilteorema}) follows easily from:

\begin{theorem}\label{main2} For all $\underline\mu$ the ring $A_X( \underline\mu):=A/I_X(\underline\mu)$ has dimension   equal to the number $d(X)$ of bases which can be extracted from $X$.\end{theorem}
\begin{proof}   We will use a standard procedure in commutative algebra.  First we show that  $\dim(A_X( \underline\mu))\geq d(X)$ for {\it generic} $\underline\mu$.  Then, for the special point $\underline\mu=0$ we show  that
$\dim(A_X )\leq d(X)$. Since the dimension is semicontinuous this implies the statement (we will  comment after about the meaning of this
theorem in terms of flatness and the Cohen  Macaulay property).\medskip

As a first step we  claim that, for generic $\mu$,  the ideal $I_X(\underline\mu)$ defines $d(X)$ distinct  points, (where $d(X)$ is the number
of bases extracted from $X$).

This implies of course that  $\dim(A_X( \underline\mu))\geq d(X)$ where equality means that all the points are reduced.
\smallskip

Set theoretically,  the variety described by all the equations given is the union of the varieties described by selecting, for every
$Y\in\mathcal E(X)$  a cocircuit,  an element $a\in  Y$ and setting the equation $a-\mu_a=0$.
\smallskip

From Lemma (\ref{comp}), the resulting list of vectors  $a_{i_1},\dots, a_{i_M}$   (extracted from all the cocircuits) generates the space, thus
the equations $a_{i_j}-\mu_{i_j}=0$  are either incompatible or define a point  of coordinates  $\mu_{i_j}$ in some basis  extracted from $X$.

\smallskip

Conversely, given such a basis every cocircuit must contain at least one of the elements $a_{i_j}$ hence the point of equations
$a_{i_j}=\mu_{i_j}$ must be in the variety.

Now it is clear that, if the $\mu_i$ are generic (i.e. do not satisfy some linear compatibility equations)  these $d(X)$ points are all
distinct.  Proving the first inequality. \bigskip

We now  prove that $\dim(A_X)\leq d(X)$. To do this we   proceed by induction on the cardinality of $X$ and $\dim(V)$. \medskip

Take $y\in X$ and set $Z:=X-\{y\}$.   Set $\overline V:=V/\mathbb Cy$ and denote by $\pi:V\to \overline V$ the quotient homomorphism.

Also, for any sublist  $B\subset Z$ set $\overline B$ equal to the list of non zero vectors in $V/\mathbb Cy$ which are images of vectors in
$B$.

Define finally $\overline A:=A/(y)=S[V/\mathbb Cy]=S[\overline V].$\medskip

\noindent  (1)  {\it  The image of $I_X$ in $\overline A $ is the ideal $I_{\overline Z}$}:
\medskip

The image of $I_X$ is generated by the products $\prod_{x\in Y}\pi(x)$ where $ Y\in\mathcal E(X)$  and $y\notin Y$.  Thus $X-Y$ generates a
proper subspace of $V$ containing $y$. This clearly implies that $\pi(\langle X-Y\rangle)\subsetneq \overline V$ and that $ \overline  Z \cap
\pi(\langle X-Y\rangle)= \overline {X-Y }$.
\smallskip

Therefore  the image $\pi(Y)=\overline  Y$ is a cocircuit for $\overline Z$ and    $$\prod_{x\in Y}\pi(x)= \prod_{\overline x\in \overline
Y}\overline x.$$  This proves that the image of $I_X$ is contained in $I_{\overline Z}$  (the other generators map to 0).\medskip

On the other hand we know that $I_{\overline Z}$  is generated by the products $\prod_{\overline x\in \overline  Y}\overline x$ with $\langle
\overline Z- \overline Y\rangle\subsetneq V/\mathbb Cy$. Setting $Y=X-  \pi^{-1}(\langle \overline Z- \overline Y\rangle)$, we immediately see
that $Y\in \mathcal E(X)$  and $\prod_{\overline x\in \overline  Y}\overline x=\prod_{x\in  Y}\pi(x)$ proving that $I_{\overline Z}$ coincides
with the image of $I_X$.\smallskip

Therefore  by induction $\dim(A/(I_X+(y)))=\dim(\overline A_{\overline  Z})=d(\overline Z)$.\bigskip

\noindent (2)  {\it The  surjective map   $A\stackrel{y}\to yA\stackrel{p_X}\to yA_X$ of $A$ modules (multiplication by $y$ followed by the
projection $p_X:A\to A_X$) factors through the projection $p_Z:A\to A_Z$  followed by   a (surjective) map}:

$$j:A_Z\to yA_X.$$\medskip

To see this it suffices to show that, if $Y\in \mathcal E(Z)$  then $y\cdot \prod_{x\in Y}x\in I_{X}$. \smallskip

We have two cases. If $y\in \langle Z-Y\rangle$, then $Y=X-  \langle Z-Y\rangle $ and $\prod_{x\in Y}x\in I_{X}$.

If $y\notin \langle Z-Y\rangle$, then $X\cap \langle Z-Y\rangle=Z-Y$ and  $X-  \langle Z-Y\rangle =X-(Z-Y)=Y\cup \{y\}\in \mathcal E(X)$ so that
$y\cdot \prod_{x\in Y}x=  \prod_{x\in Y\cup \{y\}}x\in I_{X}$.\smallskip

In particular this  proves that $\dim(A_Z)\geq \dim( yA_X).$
\medskip

Let us now use these facts to estimate $\dim(  A_X) .$\smallskip

Consider   the exact sequence
$$0\to yA_X\to A_X\to A/(I_{X},y)\to 0.$$
By what we have seen we have
$$\dim A_X =  \dim A/(I_{X},y) +  \dim yA_X \leq \dim (\overline  A_{\overline Z})+\dim (A_Z) .$$

By our inductive hypothesis, $A_Z$ has dimension equal to the number of bases of $V$ contained in $Z$ i.e. the number of bases in $X$ not
containing $y$ among its elements.

Also by induction $\dim (\overline  A_{\overline Z})$ equals to the number of bases of $V/(x)$ contained in $\overline Z$.

Given such a basis $\underline b'=\{\pi(a_{i_1}),\ldots ,\pi(a_{i_{s-1}})\},$ the set $\underline b=\{y,a_{i_1},\ldots a_{i_{s-1}}\}$ is a basis
of $V$ contained in $X$ and containing $y$ as an element. Vice versa, given a basis $\underline b=\{y,a_{i_1},\ldots a_{i_{s-1}}\}$ as above,
then $\underline b'=\{\pi(a_{i_1}),\ldots ,\pi(a_{i_{s-1}})\},$ is a basis of $V/(x)$ contained in $\overline Z$.

Thus   $\dim  (\overline  A_{\overline Z})$ equals the number of bases of $V$ contained in $X$ and containing $y$.

Summarizing $ \dim A_Z  +\dim  (\overline  A_{\overline Z})=d(X).$  The theorem is proved.

\end{proof}

Remark that, as a consequence of the proof, which we will use later, we also now have  that $j:A_Z\to yA_X $  is bijective.\medskip

As announced we want to discuss the meaning of this theorem for the variety $V_X$ given by $I_X(\underline \mu)$ thought of as ideal in
$S[V][\mu_1,\dots,\mu_N]$.

This variety is easily seen to be what is called a {\it polygraph}. It lies in $U\times \mathbb C^N$ and can be described as follows.  Given a
basis $\underline b:=\{a_{i_1},\dots,a_{i_s}\}$ of $V$, extracted from $X$, let $a^1,\dots,a^s$ be the associated dual basis in $U$. Define a
linear map  $i_{\underline b}:\mathbb C^N\to U$ by:
$$i_{\underline b}(\mu_1,\dots,\mu_N):=\sum_{j=1}^n\mu_{i_j}a^{i_j}. $$
Let $\Gamma_{\underline b}$ be its graph,  then $V_X=\cup_{\underline b}\Gamma_{\underline b}$.

$V_X$ comes equipped with a projection map $\rho$ to $\mathbb C^N$  whose fibers are the schemes defined by the ideals $I_{\overline
X}(\underline \mu).$

\begin{remark}\label{rem}
Theorem \ref{main2} implies that, $\rho$ is {\it flat}  of degree $d(X)$ and  $I_X(\underline \mu)$ is the full ideal of equations of  $V_X$.
Furthermore $V_X$  is  Cohen Macaulay.

\end{remark}

The fact is remarkable since it is very difficult for a polygraph to satisfy these conditions. When it does,  this has usually deep
combinatorial implications (see \cite{Ha}).\bigskip

One should make some remarks about the algebras $A_X(\underline\mu)$ in general.\smallskip

\noindent First  some basic commutative algebra. Take an ideal $I\subset \mathbb C[x_1,\dots,x_m]$ of a polynomial ring. $\mathbb
C[x_1,\dots,x_m]/I$ is finite dimensional if and only if the variety of zeroes  of $I$ is a finite set of points $p_1,\dots,p_k$. In this case
moreover we have a canonical decomposition $$\mathbb C[x_1,\dots,x_m]/I =\oplus_{i=1}^k \mathbb C[x_1,\dots,x_m]/I_{p_i}$$ where, for each $p\in
\{p_1,\dots,p_k\}$, the ring  $\mathbb C[x_1,\dots,x_m]/I_{p } $ is the local ring   associated to the point $p $.

Let  $p$ have coordinates  $x_i=\mu_i$,   the local ring $\mathbb C[x_1,\dots,x_m]/I_{p } $  is characterized, in terms of linear algebra, by
the property that the elements $x_i$ have generalized eigenvalue $\mu_i$. Thus the previous decomposition is just the Fitting decomposition,
into generalized eigenspaces, for the commuting operators $x_i$.\medskip

In the case of the algebra $A_X(\underline\mu)$  quotinet of $S[V]$ by the ideal $I_X(\underline \mu)$ generated by the elements, $\prod_{y\in
Y}(y+\mu_y)$ we see that, if for a point $p$ of the resulting (finite) variety  we have  that $y(p)+\mu_y\neq 0$, in the local ring of $p$   the
element $y+\mu_y$ is invertible and so it can be dropped from the equations. We easily  deduce that:
\begin{proposition} The local component $A_X(\underline\mu)(p)$  equals the algebra $A_{X_p}(\underline\mu_p)$ defined by the sublist $X_p:=\{x\in X\,|\, x(p)+\mu_x=0\}$.

Furthermore by a change of variables  $X_p(\mu):=\{ x+\mu_x,\ x\in X_p\}$ we can even identify $A_{X_p}(\underline\mu_p)=A/ \tau_p(I_{X_p })$.
\end{proposition}
By $ \tau_p(I_{X_p })$ we mean the ideal $I_{X_p }$ translated at $p$  by the automorphism of $A$ sending $x\mapsto x-x(p), \forall x\in V$
(hence $x\mapsto x+\mu_x, \ \forall x\in X_p$).

\section{Solutions of differential equations}
\subsection{Differential equations with constant  coefficients.}
Let us reinterpret Theorem \ref{main2} in terms of differential equations.
\smallskip

Let us consider the polynomial ring $S[V]$, its graded dual is $S[U]$.

By duality $S[V]$ is identified with differential operators with constant coefficients on $V$ and, using coordinates, the duality pairing can be
explicitly described as follows. Given a polynomial $p(\pd{}{x_1},\dots,\pd{}{x_s})$ in the derivatives and one $q(x_1,\dots,x_s)$ in the
variables the pairing \begin{equation}\label{pairing}\langle  p(\pd{}{x_1},\dots,\pd{}{x_s})\,|\,  q(x_1,\dots,x_s)\rangle\end{equation} is
obtained by applying $p$ to $q$ and then evaluating at 0.
\smallskip

Of course the algebraic dual of  $S[V]$ is a rather enormous space of a rather formal nature. It can be expressed  best  as  a space of {\it
formal power series} by associating to  an element $f\in (S[V])^*$ the formal expression:
$$\langle f\,|\, e^x\rangle=\sum_{k=0}^\infty {\langle f\,|\,  x^k\rangle\over k!},\quad x\in V.  $$
Where now $ {\langle f\,|\,  x^k\rangle}/ k!$  is a genuine homogeneous polynomial of degree $k$ on $V$.

The following facts are easy to see:
\begin{proposition}

\begin{enumerate}\item
Given   a vector $v\in V$,   the   transpose  of multiplication by $v$ is the directional derivative $D_v$ in $(S[V])^*$.

\item  Given $\phi\in U$, denote by $\tau_\phi$ the automorphism of $S[V]$ induced by   translation $x\mapsto x- x(\phi)=x- \langle    \phi\,|\,
x \rangle, \forall x\in V$, and by  $\tau_{\phi}^*$ its  transpose. We have $ \tau_{\phi}^* f= e^{- \langle    \phi\,|\, x \rangle}    f.$
\end{enumerate}
\end{proposition}

\begin{proof} 1) follows from the fact that for $\phi\in U$   $D_v(\phi)=\langle\phi\,|\,v\rangle$. and the chain rule.

As for 2), $$\langle \tau_{\phi}^* f\,|\, e^x\rangle=\langle  f\,|\, \tau_{\phi}(e^x)\rangle=  \langle   f\,|\, e^{x- \langle    \phi\,|\, x
\rangle }\rangle= e^{- \langle    \phi\,|\, x \rangle} \langle   f\,|\, e^{x }\rangle.  $$\end{proof} Observe that, if $J$ is an ideal of $S[V]$
defining a subvariety $Z\subset U$, we have that $\tau_\phi(J)$ defines the subvariety $Z+\phi$.\smallskip

If we have a  quotient $S[V]/I$  by an ideal $I$ we get an injection   of $i:(S[V]/I)^*\to (S[V])^*$. The image is, at least formally, the space
of solutions of the differential equations given by $I$. We denote by $Sol (I)$ the space of  $C^{\infty}$ solutions of the system of
differential equations given by $I$
\smallskip

Assume now that $S[V]/I$ is finite dimensional. Denote by $\{\phi_1,\ldots ,\phi_k\}\subset U$ the finite set of points which are the support of
$I$. Take the decomposition $$S[V]/I=\oplus_{i=1}^kS[V]/I(\phi_i)$$ where $S[V]/I(\phi_i)$ is local and  supported in $\phi_i$.  Under these
assumptions we get
\begin{theorem}\label{soleqdiff} 1) If $S[V]/I$  is finite dimensional and supported at 0. The image of $i^*$ lies  in $S[U]$ and coincides with $Sol (I)$.

2) If $S[V]/I$   is finite dimensional and supported at  a point $\phi\in U$,
$$Sol (I)=e^{ \langle    \phi\,|\, x \rangle}Sol (\tau_{-\phi} I).$$

3) For a general finite dimensional $S[V]/I$,
$$Sol (I)=\oplus_{i=1}^kSol(I(\phi_i))$$\smallskip
\end{theorem}

\begin{proof} 1) is the duality. As for 2),  clearly, for any ideal $I$ and $\phi$ we have  $Sol(\tau_\phi(I))=\tau_{-\phi}^*(Sol(I))=e^{ \langle    \phi\,|\, x \rangle}Sol(I)$.
It only remains to explain the analytic nature of the statement. This follows, even in the stronger setting of tempered distributions by the
fact the a tempered distribution which is annihilated by high enough derivatives is necesarily a  polynomial (see \cite{DM1} Prop. 3.1).
\end{proof}
Applying this discussion to our setting,  we get that, for generic $\mu$,  the solutions of the differential equations of the ideal
$I_X(\underline \mu)$ are exactly the space with basis the functions $e^{\phi}$  where the elements $\phi$ run on the $d(X)$ points defined
previously, while for $\underline\mu=0$ we have a remarkable space of polynomial solutions.  The nature of this space is explained in the next
section.  \smallskip

In the general case  we have a mixture between these two extreme cases. The set of solutions is a finite set of points and in each such point
$\phi\in U$ we have a subset  $X_{\phi}$ of $X$ formed by those elements $x_i\in X$ such that $\langle\phi\,|\,x_i\rangle=\mu_i$.

Then one can see that this point contributes to the solutions with the functions $e^\phi p$ where $p$ are the polynomials associated to
$X_{\phi}$ at 0.

\section{A realization of $A_X$ }

\subsection{The graded dimension}

Our next  task is to prove that, the graded dimension of the space $D(X)$ is given by:
\begin{theorem}\label{main}
$$ H_X(q)=\sum_{\underline b\in\mathcal B(X)}  q^{N  -n(\underline b) }.$$
\end{theorem}

In order to do this, we want  to realize the ring $A_X$ as the $S[V]$ submodule $Q_X$, of the space of polar parts $P_X$ generated by the class
$u_X$, of the element $d_X^{-1}$. In doing this  we shall exhibit a homogeneous basis made by elements $u_{\underline b}$ as $\underline b$
varies among the bases of $X$,  each of degree $N  -n(\underline b) $.

The spaces $R_X$ and  $P_X$ are naturally graded (as functions).  It is convenient to shift this degree by $N$ so that $u_X$ has degree $0$ and
the generators $u_{\underline b}$ have degree   $ N-s$. If $\underline b$ is a no broken basis, these will be just the elements $u_{\underline
b}$ introduced in section \ref{31} when $\underline \mu=0$.

With these gradation the natural map $\pi:A\to Q_X$ defined by $\pi(f)=fu_X$ preserves degrees.

It is clear that  Theorem (\ref{main})  follows from  the following more precise result on $Q_X$, which also describes  a graded basis for
$A_X$:

\begin{theorem}\label{main1}  \begin{enumerate}\item
The annihilator of $u_X$ is the ideal   $I_X$ generated by the elements $M_{Y}=\prod_{x\in  Y}x,$ as $Y$ runs over the cocircuits. Thus
$Q_X\simeq A_X$ as  graded $A$-modules.

\item  Given a basis  $\underline b:=\{a_{i_1},\dots, a_{i_s}\}$ extracted from $X$,  set $$u_{\underline b}:=(\prod_{a\in X-B(\underline
b)}a)u_X.$$  The elements   $u_{\underline b}$,   as $\underline b$ runs over the bases extracted from $X$, are a basis of $Q_X$.

\end{enumerate}
\end{theorem}
\proof

To see the first part, let us remark that clearly $  I_Xu_X=0$ since, if $Y$ is a cocircuit $M_Yu_X=\prod_{x\in X-Y}x^{-1}$ lies in $R_{X,
s-1}$, hence it   is 0 in the module $P_X$.

Therefore, from Theorem \ref{main2}, it is enough to see that $\dim Q_X\geq d(X)$.
\smallskip

We want to proceed by induction on $s$ and on the cardinality of $X$.

If $X$ consists of a basis of  $V$, clearly both $A_X$ and $Q_X$ are 1-dimensional and the claim is clear.   Otherwise we can assume that
$X:=\{Z,\, y\}$ and $Z$ still spans  $V$. We need thus to compare several of the objects under analysis in the case in which we pass from $Z$ to
$X $.

Let us consider the ring $A/(y)$, polynomial functions  on the subspace of $U$ where $y$ vanishes and   denote by $ \overline Z$ the set  of non
zero vectors in the image of $Z$ (or $X$) in $A/(y)$.\medskip

As in Theorem (\ref{main2}) the set $\mathcal B_X$ of bases  extracted from $X$ can be decomposed into two disjoint sets.

$\mathcal B_Z$   and the bases  containing $y$. This second set is in 1-1 correspondence with the bases of $V/\mathbb Cy$ contained in
$\overline Z$.

We need several lemmas.

First  we obviously have an inclusion $R_Z\subset  R_X$ and also $R_{Z,k}\subset  R_{X,k},\ \forall k$.
\begin{lemma}\label{inc}  $R_{Z,s-1}=R_Z\cap   R_{X,s-1}$ and  we get an inclusion of
$ P_Z$ into $ P_X.$

\end{lemma}
\proof The two statements are equivalent and follow immediately by ordering $Z$ and adding $y$ as last element. Then all the no broken circuit
bases for $Z$ are also no broken circuit bases for $X$ so that $P_Z$ (as module over the differential operators with constant coefficients) is a
free direct summand of $P_X$.\qed

\smallskip
As a consequence of this Lemma let us consider in $P_X$ the map multiplication by $y$. We have clearly $yu_X=u_Z$ thus we obtain an exact
sequence of $A-$modules:
$$ 0\to K\to Q_X\stackrel{y}\to Q_Z\to 0$$
where $K=Q_X\cap Ker (y)$.\smallskip

We need to analyze $K$ and prove that $\dim(K)$ is greater or equal to the number $d_y(X)$, of bases extracted from $X$ and containing $y$.
Since we already know that $\dim(K)+d(Z)\leq d(X)$  and  $d(X)=d(Z)+d_y(X)$  this will prove the claim.\medskip

In order to achieve the inequality $\dim(K)\geq d_y(X)$ we  will find inside $K$ a direct sum of subspaces whose dimensions add up to $d_y(X)$.

Let us first discuss a special case. Assume that $Z$   spans a subspace $V'$ of codimension 1 in $V$,   so $y$ is a vector outside this
subspace.

We clearly have inclusions:
$$R_Z \subset R_X\quad y^{-1}R_{Z,k-1} \subset R_{X,k}\quad \forall k.$$
Also the element  $u_X\in P_X$  is killed by $y$ and $A/(y)$ can be identified to $S[V'].$

\begin{lemma}\label{multip} The multiplication by  $y^{-1}$ induces an isomorphism between $P_Z$ and the kernel of the multiplication by $y$ in $P_{X}$.
\end{lemma}
\begin{proof} Since $y^{-1}R_{Z,k-1} \subset R_{X,k}$ it is clear that the multiplication by $y^{-1}$ induces a map from $P_Z=R_{Z,s-1}/R_{Z,s-2}$ to $P_X=R_{X,s}/R_{X,s-1}$. It is also clear that the image of this map lies in the kernel of the multiplication by $y$.

To see that it gives an isomorphism to this kernel, order the elements of $X$ so that $y$ is the first element. A no broken circuit basis for
$X$  is of the form  $\{y,\underline c\}$ where $\{\underline c\}$ is a no broken circuit basis for $Z$.

Now fix  a set of coordinates $x_1,\dots,x_s$ such that $x_1=y$ and $x_2,\dots,x_s$ is a basis of the span of $Z$. Denote by $\partial_i$ the
corresponding partial derivatives. We have that in each summand $\mathbb  C[\partial_1,\dots,\partial_s]u_{y,\underline c}$  the kernel of
multiplication by $x_1$ coincides with $\mathbb  C[\partial_2,\dots,\partial_s]u_{y,\underline c}$.

The claim follows easily since $\mathbb  C[\partial_2,\dots,\partial_s]u_{y,\underline c}$  is the image, under $y^{-1}$ of $\mathbb
C[\partial_2,\dots,\partial_s]u_{ \underline c}$.  \end{proof}

Notice that $y^{-1}u_Z=u_X$, so that the following lemma is immediate

\begin{lemma}\label{mulq}  Under  the previous hypotheses, the     multiplication by  $y^{-1}$ induces an isomorphism between $Q_Z$ and $Q_X$.\end{lemma}

\smallskip

Let us now pass to the general case. Consider  the set  $\mathcal S_y(X)$  of all complete sublists   of $X$  which span a subspace of
codimension  1 not containing  $y$.

For each  $Y\in \mathcal S_y(X)$ we have, by  Lemma \ref{multip}, that the multiplication by  $y^{-1}$ induces an inclusion   $i_Y:P_Y\to
P_{Y\cup\{y\}}\to P_{X}$ with image in the kernel $Ker(y)\subset P_X$.  Thus we get a map
$$g:=\oplus_{Y \in \mathcal S_y(X)}i_Y:\oplus_{Y \in \mathcal S_y(X)}P_Y\to Ker(y)$$
We then claim that:

\begin{lemma}\label{iso} $g$ is  an isomorphism of $\oplus_{Y \in \mathcal S(X)}P_Y$ onto  $Ker(y)$.

\end{lemma}
\proof As before order the elements of $X$ so that $y$ is the first element. A no broken circuit basis for $X$  is of the form  $\{y,\underline
c\}$ where $\{\underline c\}$ is a no broken circuit basis for $Y:=X\cap\langle \underline c\rangle$.

By construction $Y\in \mathcal S_y(X)$ and we obtain the direct sum decomposition $$P_{X}=\oplus_{Y\in\mathcal S_y(X)}P_{Y\cup\{y\}}.$$ Clearly
$$Ker(y)=\oplus_{Y\in\mathcal S(X)}P_{Y\cup\{y\}}\cap Ker(y)$$ and, by Lemma \ref{multip}, for each $Y\in  \mathcal S_y(X)$,  $i_y$ gives an
isomorphism of $P_Y$ with  $P_{Y\cup\{y\}}\cap Ker(y)$. This proves the lemma.

\qed

In order to finish the proof of Theorem \ref {main1} notice that by Lemma \ref{mulq} and Lemma \ref{iso}, we get an inclusion of $\oplus_{Y \in
\mathcal S_y(X)}Q_Y$ into $K=Q_X\cap Ker (y)$, so that $\dim Q_X\geq \dim Q_Z + \sum_{Y\in\mathcal S_y(X)}\dim P_Y=d(X)$.

This gives the required inequality and implies also that we have a canonical exact sequence of modules:
\begin{equation}\label{seqes}0\to  \oplus_{Y\in\mathcal S_y(X)}Q_Y\to Q_{X}\stackrel{y}\to Q_Z\to 0.\end{equation}

It remains to establish the second part of Theorem \ref {main1}. We prove it by induction, ordering $X$ so that $y$ is the last element.  Let us
write $\mathcal B_X$ as the disjoint union of $\mathcal B_Z$ and the set $\mathcal B^{(y)}_X$ of bases containing $y$. Also set $\mathcal
U_X=\{u_{\underline b}\,|\,\underline b\in\mathcal B_X\}$ and write it as the dijoint union of $\mathcal U'_X=\{u_{\underline b}\,|\,\underline
b\in\mathcal B_Z\}$ and its complement $\mathcal U''_X=\{u_{\underline b}\,|\,\underline b\in\mathcal B^{(y)}_X\}$. Remark that by our
definitions, $y\mathcal U'_X=\mathcal U_Z$ while $\mathcal U''_X$ is the disjoint union $\cup_{Y\in S_y(X)}i_Y(\mathcal U_Y)$. By induction and
formula (\ref{seqes}) this clearly implies that $\mathcal U_X$ is a basis of $Q_X$.\qed
\medskip

It follows from our theorem that in top degree $Q_X$ has as a basis the elements $u_{\underline b}$ as  $\underline b$ runs over the set of no
broken bases.\medskip

In the parametric case, we have the decomposition $d_X^{-1}=\sum_pc_pd_{X_p}^{-1}$  (cf. (\ref{separa})). We set     $Q_X(\underline\mu)$ equal
to the $S[V]$-module generated by $u_X$ and for each $p\in P(X,\underline \mu)$, $Q_X(p)$ equal to the  $S[V]$-module generated by $u_{X_p}$.
\begin{proposition}
1)  $$Q_X(\underline\mu)=  \oplus_{p\in P(X_p,\underline\mu_p)}   Q_{X_p } (\underline\mu_p).$$

2) $Q_X(\underline\mu)$ is isomorphic to $A_X(\underline \mu)$ and the above decomposition coincides with the decomposition of $A_X(\underline
\mu)$ into its local components.
\end{proposition}
\begin{proof}  The first part follows immediately from  Remark \ref{sommandi}.

As for the second since  the annihilator of $u_{X}$  contains  $I_{X,\underline \mu_p}$, we have a map $A_X(\underline \mu)\to
Q_X(\underline\mu)$ which by Fitting decomposition induces a map of the local summands. On each such summand this map is an isomorphism by the
previous theorem since we can translate $p$ to $0$. From this everything follows.
\end{proof}

We can easily prove as corollary a theorem by several authors, see  \cite{AS}, \cite{DR},\cite{J}.
\begin{corollary} Consider in $S[V]$ the subspace $\mathcal P(X)$
spanned by all the products $M_Y:=\prod_{x\in Y}x,\ Y\subset X$ such  that $X-Y$ spans  $V$. Then $\mathcal P(X)$ is in duality with
$D(X)$.\end{corollary}
\begin{proof} Multiply by $d_X^{-1}$. $\mathcal P(X)d_X^{-1}$
is spanned by the polar parts $\prod_{x\in Z}x^{-1},$ $ Z\subset X$  such that $Z$ spans  $V$. We have seen that this space of polar  parts maps
isomorphically to its image into $P_X$ and its image is  clearly $Q_X$. This proves the claim.
\end{proof}
\begin{remark}   The last theorem we have proved is equivalent to saying  that, in the algebra $R_X$   the intersection   $d_X^{-1}A\cap R_{X,s-1}=d_X^{-1}I_X.$

\end{remark}

It is of some interest to analyze the deeper intersections  $d_X^{-1}S[V]\cap R_{X, k}:=U_k$.  We will sketch this point which uses the
structure of the filtration by polar degree as explained in \cite{dp3}.

If $X$ does not necessarily span $V$ we define $I_X$ as the ideal generated by the products $\prod_{x\in Y}x$ where $Y\subset X$ is any   subset
such that the span of $X-Y$ is strictly contained in the span of $X$.  We set $A_X(V)=A/I_X$.

Of course if we fix a decomposition  $V:=\langle X\rangle\oplus T$ we have $A_X(V)=A_X\otimes S[T]$.

Consider the set of all subspaces $W$ spanned by elements of $X$ (including the space $\{0\}$ spanned by the  empty set). We call any such
subspace a {\it subspace of the arrangement} generated by $X$. Set $L=d_X^{-1}A$, $X_W:=X\cap W$,  and $L_k=L\cap R_{X,k}$. We have
\begin{theorem}\label{filtr}  For each $k$ we have that   $L_k/L_{k-1}$ is isomorphic to the direct sum  $\oplus A_{X_W}(V)$ as $W$ varies on the subspaces  of dimension $k$ of the arrangement.

\end{theorem}
\begin{proof}
Under the map $A\to L$ given by $f\to d_X^{-1}f$, we know that  the submodule  $I_X$ maps isomorphically onto $L_{s-1}$. So $L_{s-1}$  is
spanned by the elements $d_Y^{-1}$  where $Y=X_W$ for the hyperplanes $W$ of the arrangement. By (a small generalization of) Lemma \ref{inc} we
know that the filtration by polar order in $R_X$ induces the same filtration in $R_Y$  thus by induction we have that the graded associated to
$d_Y^{-1}S[V]$ is the direct sum of all the pieces $A_{X_W}$ as $W$ varies on the subspaces  of dimension $k$ of the arrangement  generated by
$Y$. This follows from Proposition \ref{ilrema}.
\end{proof}
The previous theorem gives an interesting combinatorial identity once we compute the Hilbert series of $d_X^{-1}S[V]$ directly or as sum of the
contributions given by the previous filtration.
$$\frac{q^{-|X|}}{(1-q)^s}=\sum_{k=0}^{s-1}\sum_{W\in \mathcal W_k}\frac{q^{-|X_W|}}{(1-q)^{s-k}}H_{X_W}(q). $$\medskip

\begin{remark}
  In Fourier transform the previous discussion can be translated into an analysis of the distributional derivatives of the multivariate spline, i.e. into an analysis of the various discontinuities achieved by the successive derivatives on all the  strata of the singular set.
  \end{remark}
  \bigskip

In the applications to the box spline it is interesting, given a {\it set} $X$ of vectors which we list is some way, to consider for each $k\geq
0$ the list $X^k$ in which every element $a\in X$ is repeated $k$ times.  Let us make explicit the relationship between $H_X(q)$ and
$H_{X^k}(q)$.

First the number of bases in the two cases is clearly related by the formula $d(X^k)=d(X)k^s$, to each basis $\underline b:=\{b_1,\dots,b_s\}$
extracted from $X$ we associate  $k^s$ bases $\underline b(h_1,\dots,h_s)$,   indexed by $s$ numbers $h_i\leq k$  expressing the position of the
corresponding $b_i$ in the list of the $k$ repeated terms in which it appears.

Now it is easy to see that:
$$ n(\underline b(h_1,\dots,h_s))=k(n(\underline b)-s)+ h_1+\dots +h_s.    $$
Thus we deduce the explicit formula:
$$  H_{X^k}(q)=\sum_{\underline b\in\mathcal B_X }\sum_{h_1,\dots,h_s}q^{kN- k n(\underline b) +ks - h_1-\dots -h_s }= H_X(q^k)\big({q^s-1\over q-1}\big)^k.$$

\subsection{More differential equations}   Our next task  is to fully characterize the  polynomials $p_{\underline b}( x)$  appearing in formula (\ref{JK2}) by differential equations. In Theorem \ref{exxpr}, we have seen that these polynomials lie in $D(X)$.\medskip

For a given no broken circuit basis $\underline b$, consider the element  $D_{\underline b}:=\prod_{a\notin \underline b}a$.
\begin{proposition}
The polynomials  $p_{\underline b}$ satisfy the system
\begin{equation}\label{dual}D_{\underline b}p_{\underline c}( x_1 ,\dots,  x_s) =    \begin{cases} 1\quad\text{if}\quad \underline b=\underline c\\ 0 \quad\text{if}\quad \underline b\neq \underline c\end{cases}\end{equation}
In particular the polynomials  $p_{\underline b}$ for $\underline b\in\mathcal{NB}$  are linearly independent and   characterized, in
$D(X)_{N-s}$ by the equations (\ref{dual}).
\end{proposition}

\begin{proof} This follows from the identity $D_{\underline b}u_X=u_{\underline b}$, the linear independence of the elements $u_{\underline b}$ and  the fact that the dimension of  $D(X)_{N-s}$  is the cardinality of $\mathcal{NB}$.\end{proof}

As a consequence, using formula (\ref{JK}),  we can characterize by differential equations the multivariate spline $T_X(x )$  on each big cell
$\mathfrak c$  as the function $T$ in  $D(X)_{N-s}$ satisfying the equations:
\begin{equation}\label{dualc}D_{\underline b}T =    \begin{cases} |\det(\underline b)|^{-1}\quad\text{if}\quad  \mathfrak c\subset C(\underline b)\\ 0 \quad\text{otherwise} \end{cases}\end{equation}

We have identified $D(X)$   to the dual of $A_X$  and hence of   $Q_X$. The basis   $u_{\underline b}$, we found in $Q_X$, defines thus a dual
basis  $u^{\underline b}$ in $D(X)$.

\begin{corollary}
When $ {\underline b}\in\mathcal{NB}$   we have, as polynomials   $u^{\underline b}  =p_{\underline b}$
\end{corollary}
\begin{proof}
This follows from the previous Proposition and the definition of the duality given in formula (\ref{pairing}).
\end{proof}

We shall use a rather general notion of the theory of modules. Recall that the socle $s(M)$ of a module $M$  is the sum of all its irreducible
submodules. Clearly if $N\subset M$ is a submodule $s(N)\subset s(M)$ while for a direct sum $s(M_1\oplus M_2)=s(M_1)\oplus s(M_2)$.  If $M$ is
finite dimensional   $s(M)\neq 0$ so that, for a non--zero submodule $N$ we must have $N\cap s(M)\neq 0$.
\begin{proposition} The socle of the $S[V]$-module $Q_X$ coincides with its top degree part, with basis the elements $u_{\underline b}$.
\end{proposition}
\begin{proof}
The socle of the    algebra of constant coefficients differential operators, thought of as a module over the polynomial ring,  is clearly
generated by 1. It follows that the socle of $P_X$ (as  $S[V]$-module), has as  basis the elements $u_{\underline b}$. Since $u_{\underline
b}\in Q_X$ we have $s(Q_X)= s( P_X ),$   the claim follows.
\end{proof}

\begin{theorem} $D(X)$ is spanned by the polynomial   $p_{\underline b}$, as $\underline b$ runs over the set of no broken bases, and all of their derivatives.\end{theorem}

\begin{proof}
Everything follows once we observe:

1) $D(X)$ is in duality with $A_X$ and so with $Q_X$.

2) The polynomials $p_{\underline b}$, $\underline b\in \mathcal {NB}$ are in duality with the elements $u_{\underline b}$.

3) The orthogonal of a proper submodule $N$ of $D(X)$ is a non--zero submodule, thus it must intersect the socle of $Q_X$ in  a non trivial
subspace. In particular $N$ cannot contain all the elements  $p_{\underline b}$.
\end{proof}

We finish this section observing that we can dualize the sequence  (\ref {seqes}) to get an exact sequence
$$0\to D(Z) \stackrel{D_y}\longrightarrow D(X)\longrightarrow\oplus_{Y\in\mathcal S_y(X)}D(Y) \to 0.$$
Here, since the map $Q_{X}\to Q_Z\to 0$ is given by multipliation by $y$, the inclusion $0\to D(Z) \to D(X)$ is the inclusion of    $D_y(D(Z))$.

\part{The discrete case}
\bigskip

\section{ The discrete case}

\subsection{The   case $X$ in a lattice}
Splines are used to  interpolate and approximate functions. For this purpose it is important to understand the class of smoothness of a spline.   In the case of the box--spline it is  easy to prove (see \cite{dhr}):

Let $m$ be the minimum number of elements in a cocircuit in $X$, assume   $m\geq 2$, then
    $B_X$ is of class $C^{m-2}$.

Thus, provided we choose the list $X$ in a suitable way, we can achieve any finite level of smoothness required.\smallskip

A particularly useful case is  when     we have chosen a lattice    $\Lambda\subset V$ spanning $V$ such  that each vector in the list $X$ lies in $\Lambda$.  We assume that the Lebesgue measure on $V$ is normalized in such a way that a fundamental domain for $\Lambda$ has volume 1.\smallskip

In this case, if $C_{sing}(X)$  denotes the set of singular vectors of the cone  $C(X)$, the set of all translates
$\cup_{\lambda\in\Lambda}C_{sing}(X)$ is called the {\it cut region}. It is a union of a finite number of  $s-1-$dimensional bounded polytopes
and all their translates.

From   formula (\ref{boxc}),   and the fact that $D(X)$ is stable under translation we obtain (see also \cite{dhr}):
\begin{proposition} The complement of the cut region is   a union of  all translates of a finite number of cells, each an
   interior of a (compact) polytope.

Over each such cell   the functions $B_X(x-\lambda),\ \lambda\in \Lambda$ are polynomials in the space $D(X)$.

\end{proposition}\smallskip

One easily proves the fundamental fact  (see \cite{dhr}): \begin{theorem}  If $X$ spans $V$, the translates of $B_X$ form a {\it partition of unity}:
\begin{equation}\label{partun}1=\sum_{\lambda\in\Lambda}B_X(x-\lambda).\end{equation}
 \end{theorem}
 \begin{proof}   If $X$ is a basis, $B_X$  is the characteristic function of the parallelepiped
 with basis $X$ divided by its volume and it easily follows that  the identity  (\ref{partun})
 holds outside a set of measure 0.  In general one can use the iterative description of the box spline obtained by integration:
 $$ B_{X,v}(x)=\int_0^1B_X(x-tv)dt,$$ a formula deduced immediately from the definition (\ref{box}).

$$\int_0^1\dots\int_0^1f(\sum_{i=1}^Nt_ia_i)dt_1\dots dt_N=\int_0^1\int_{\mathbb R^s}f(x+tv)B_X(x)dx dt =$$
$$\int_0^1\int_{\mathbb R^s}f(x )B_X(x-tv)dx dt =\int_{\mathbb R^s}f(x )(\int_0^1B_X(x-tv)dt)dx .$$
\end{proof}

\noindent  Given a function $f(\lambda)$ on the lattice $\Lambda$ define
$$B_X*f(x):=\sum_{\lambda\in\Lambda} B_X(x-\lambda)f(\lambda)$$
The space of all functions obtained by this procedure is called the {\it cardinal spline space} and   the space of polynomials $D(X)$ is
characterized by the property of consisting exactly of the polynomials lying in the cardinal spline space  (see \cite{dhr}).  An extensive
portion of the literature on the box spline is devoted to understanding how to use this space for approximation or interpolation of
functions.

\subsection{The Partition functions} \begin{definition}  We identify a function   $f$ on $\Lambda$  with the distribution
$$\sum_{\lambda\in\Lambda}f(\lambda)\delta_\lambda$$
where $\delta_v$ is the Dirac distribution supported at $v$.\end{definition}
Recall that we have defined the partition function on $\Lambda$ by
$$\mathcal P_X(v):=\#\{(n_1,\dots,n_N)\,|\, \sum n_ia_i=v,\ n_i\in\mathbb N\}.$$
This is then identified with the tempered distribution
\begin{equation}\label{parto}\mathcal T_X:=\sum_{v\in\Lambda}\mathcal P_X(v)\delta_v\end{equation}
 We can then compute its Laplace transform obtaining
\begin{equation}\label{laparto}
L\mathcal T_X=\sum_{v\in\Lambda}\mathcal P_X(v)e^{-v}=\prod_{a\in X} \frac{1}{1-e^{-a}}
\end{equation}
We shall think of $\mathcal T_X$ as a discrete analogue of the multivariate spline $T_X$. We also have an analogue of $B_X$. Namely, setting
$$\mathcal Q_X(v):=\#\{(n_1,\dots,n_N)\,|\, \sum n_ia_i=v,\ n_i\in\mathbb \{0,1\}\}.$$
\begin{equation}\label{partob}\mathcal B_X:=\sum_{v\in\Lambda}\mathcal Q_X(v)\delta_v\end{equation}
with Laplace transform
\begin{equation}\label{laparto2}
L\mathcal B_X=\sum_{v\in\Lambda}\mathcal Q_X(v)e^{-v}=\prod_{a\in X} ({1+e^{-a}})=\prod_{a\in X} \frac{1-e^{-2a}}{1-e^{-a}}
\end{equation}
which implies that
\begin{equation}\mathcal B_X(x)=\sum_{S\subset X}(-1)^{|S|}\mathcal T_X(x-2a_S)\end{equation}
with $a_S=\sum_{a\in S}a$.

\subsection{Some basic modules and their transforms} As in Section \ref{duemod}, we want to consider the module $\mathcal L_X$ generated,
in the space of tempered distribution on $V$, by the element $\mathcal T_X$ under the action of  a suitable algebra of operators.

It is convenient to  choose the {\it algebra $\mathcal W(\Lambda)$ of difference operators with polynomial coefficients.}

This algebra is generated by $S[U]$, thought of as   polynomials on $V$ and by the translation operators $\tau_v$, $v\in\Lambda$ defined by
$\tau_vf(a)=f(a+v)$.  From the basic formulas (\ref{L1}), (\ref{L2}) we get that, under  Laplace transform, a polynomial becomes  a differential
operator  with constant coefficients while the translation $\tau_v$ becomes multiplication by $e^v$.\smallskip

Consider the group algebra   $\mathbb C[\Lambda]$ with basis the formal elements $e^a,\ a\in\Lambda$.  Define next the algebra $\mathcal W(V)$
generated by   $S[U]$, thought of as   differential operators with constant coefficients on $U$ and by the functions  $e^v, v\in \Lambda$.
Additively we have $\mathcal W(V)=S[U]\otimes \mathbb C[\Lambda]$.

The algebras $\mathcal W(\Lambda)$ and $\mathcal W(V)$ are isomorphic by the isomorphism $\phi$ defined by
$$\phi(\tau_v)=e^v,\ \ \phi(u)=-D_u$$
for $v\in\Lambda$, $u\in U$. Thus, given a module $M$ over $\mathcal W(\Lambda)$ we shall denote by $\hat M$ the same space considered  as a
module over $\mathcal W(V)$ and think of it as a {\it formal Fourier transform}.\smallskip

Define $S_X:=\mathbb C[\Lambda][\prod_{a\in X} (1-e^{-a})^{-1}$ to be  the localization  of $\mathbb C[\Lambda]$, obtained by inverting
$\delta_X:=\prod_{a\in X} (1-e^{-a})$, and consider $S_X$ as a module over $\mathcal W(V)$.

Given a linearly independent subset $A\subset X$, let $$\Lambda_A:=\{\sum_{a\in A}n_aa,\ n_a\in \mathbb Z\}, \quad C(A):=\{\sum_{a\in A}r_aa,\
n_a\in \mathbb R^+\},$$  be the sublattice and the positive cone that they generate. Set $$\xi_A:=\sum_{v\in \Lambda_A\cap C(A)}\delta_v $$ a
tempered distribution with Laplace transform given by:
\begin{equation}\label{DisLap}
L\xi_A=\prod_{a\in A}{1\over 1-e^{-a}}.
\end{equation}

In analogy with Theorem \ref{Lapp}, using \cite{dp3}, we get,

\begin{theorem}\label{LappTr}   Under Laplace transform,  $\mathcal L_X$    is mapped isomorphically  onto $S_{X}$.
In other words we get a canonical isomorphism of $\hat {\mathcal L}_X$ with  $S_{X}$ as $\mathcal W(U)$-modules.

$\mathcal L_X$ is the space of tempered distributions which are linear combinations of polynomials times $\xi_A$, $A\subset X$ a linearly
independent subset, and their translates under $\Lambda$.
\end{theorem}
\begin{proof} The  proof   is analogous to that of Theorem \ref{Lapp}.  We use the results on partial fractions explained in \cite{dp3} to show that $S_X$ is generated by $\delta_X^{-1}$ as a $\mathcal W(V)$-module, together with the simple formula (\ref{DisLap}). \end{proof}
\begin{remark} There is a similar parametric case where the Laplace transform takes values in the algebra $S_{X,\underline  \mu}:=\mathbb C[\Lambda][\prod_{a\in X}(1-e^{-a-\mu_a})^{-1}]$.
\end{remark}
Notice that the fact that the Laplace transform of $\mathcal L_X$  is an algebra means that $\mathcal L_X$ is closed under convolution.\medskip

\subsection{The toric arrangement\label{ta}}

The algebra $S_X$ has a precise geometric meaning.  In fact $\mathbb C[\Lambda]$ is the coordinate ring of an affine variety whose points are
the    homomorphisms $\Lambda\to \mathbb C^*$   called   {\it complex characters}. We shall denote this variety by $T$.

Notice that each character is obtained as follows. One takes a vector  $\phi\in U_{\mathbb C}=hom(V,\mathbb C)$ (the complexified dual) and
constructs the function   $a\mapsto e^{\langle \phi\,|\,  a\rangle}$. The elements  of $\Lambda^*:=\{\phi\in  U_{\mathbb C}\,|\,\langle
\phi\,|\,a\rangle\in 2\pi i\mathbb Z\}$  give rise to the trivial character. Thus  $T=U_{\mathbb C}/\Lambda^*$ is an algebraic group isomorphic
to $(\mathbb C^*)^s$.\medskip

The expression $  e^{\langle \phi\,|\,  a\rangle}$ has to be understood as a pairing, i.e. a map
$$    e^{\langle \phi\,|\,  a\rangle}:T\times\Lambda\to \mathbb  C^*.$$
This duality expresses the fact that $\Lambda$ is the group of  algebraic characters of $T$.

The class of a vector $ \phi\in  U_{\mathbb C}$ in $T$ will be denoted by $e^\phi$ so that the value of $e^a$ in $e^\phi$ is $e^{\langle
\phi\,|\, a\rangle}$.\smallskip

If we fix an integral basis $e_i$ for $\Lambda$ and set $x_i:=e^{e_i}$,  we see that $\mathbb C[\Lambda]= \mathbb C[x_1^{\pm 1},\dots,x_s^{\pm
1}]$  is the ring of Laurent polynomials, i.e. the coordinate ring of algebraic functions on the {\it standard torus  $(\mathbb
C^*)^s$.}\smallskip

The equation  $\prod_{a\in X}(1-e^{-a})=0$,  defines an hypersurface $Y$ in $T$, the union of the kernels of all the characters   $e^a$, $a\in
X$. The algebra $S_X$ is clearly the coordinate ring of the complement of $Y$ in $T$, itself an affine variety.

Given a basis $\underline b$ extracted from $X$, consider the lattice $\Lambda_{\underline b}\subset\Lambda$ that it generates in $\Lambda$.

We have that $ \Lambda/\Lambda_{\underline b}$ is a finite group of order  $ [\Lambda:\Lambda_{\underline b}]=|\det(\underline b)|.$

Its character group  is the finite subgroup  $T( {\underline b})$ of $T$  which is the intersection of the kernels of the functions   $e^a$ as
$a\in  \underline b .$

We now define the {\it points of the arrangement}
\begin{equation}\label{ipunti}P(X):=\cup_{\underline b\in\mathcal B(X)}T( {\underline b}).\end{equation}

For any point in $P(X)$ we choose once and for all a representative $\phi\in U$ so that the given point equals $e^\phi$. We will denote by
$\tilde P(X)$ the corresponding set of representatives. We now set:
\begin{equation}\label{xfi}X_\phi:=\{a\in X\,|\,  e^{\langle \phi\,|\,a\rangle}=1 \}.\end{equation}

The points of the arrangement form the zero dimensional pieces of the entire {\it toric arrangement}. By this we mean the finite set formed of
all connected components of all the intersections of  the hyperfurfaces of equations $1-e^{-a}=0,\ a\in X$. For  one such intersection  the
connected component through 1 is a subtorus $T'$, these other components are  cosets of $T'$.

\subsection{$S_X$ as a module} As in Section \ref{duemod} and following
\cite{dp3}  for every $k\leq s$,  let us consider the submodule $S_{X,k}\subset S_X$ spanned by the elements $${f\over
\prod_i(1-e^{-a_i})^{h_i}}$$ such  that the vectors $a_i$ which appear in the denominator with positive exponent span a subspace of dimension
$\leq k$.\smallskip

Given a point $e^\phi\in P(X)$ and a basis $\underline b:=\{a_1,\dots, a_s\}$ extracted from $X_\phi$,  we have seen that $e^\phi$ belongs to
the  finite  group $T( {\underline b})$, in duality with  $ \Lambda/\Lambda_{\underline b}$. The finite dimensional group algebra $\mathbb C[
\Lambda/\Lambda_{\underline b}]$  is identified with the functions on  $T( {\underline b})$.

We choose a set $R_{\underline b}$ of representatives in $\Lambda$ for the cosets  $\Lambda/\Lambda_{\underline b}$ by taking for each coset
the unique representative of the form $$\sum_{a\in \underline b}p_aa,\ 0\leq p_a<1, p_a\in\mathbb Q.$$  Character theory tells us that   the
element  $e(\phi):=|\det(\underline b)|^{-1}\sum_{\lambda\in R_{\underline b}}e^{-\langle \lambda\,|\,\phi\rangle}e^\lambda$  has the property
of taking the value 1 on the point $e^\phi$ and 0 on all the other points of  $ T( {\underline b})$.
Moreover:
\begin{proposition}\label{ltd}  $e(\phi)\prod_{a\in\underline b}(1-e^{-a})^{-1}$ is the Laplace transform of the distribution  $$|\det(\underline b)|^{-1}\sum_{v\in C(\underline b)\cap \Lambda}e^{\langle \phi | v\rangle}\delta_v.$$
\end{proposition}
\begin{proof}  By our choice of representatives, if $v\in C(\underline b)\cap \Lambda$ we can write $v=\lambda+\sum_{a\in\underline b} n_aa,\ n_a\in\mathbb N,\ \lambda\in R_{\underline b}$

\end{proof}

In  $\Pi_X:=S_X/S_{X,s-1}$ the  class
$$\omega_{\underline b,\phi}:=\big[\frac {e(\phi)}{\prod_{a\in \underline b} (1-e^{-a})}\big],\quad \mod (S_{X,s-1}),$$  is clearly
independent of the chosen representatives and   is an eigenvector for $\Lambda$ of eigenvalue $e^\phi$. Indeed $$e^\mu(\omega_{\underline
b,\phi})=e^{ \langle  \mu\, \,\phi\rangle}\big[\frac{|\det(\underline b)|^{-1} \sum_{\lambda\in R_{\underline b}}e^{-\langle
\lambda+\mu\,|\,\phi\rangle}e^{\lambda+\mu}}{{\prod_{a\in \underline b} (1-e^{-a})}}\big]=e^{ \langle  \mu\,|\,\phi\rangle}\omega_{\underline
b,\phi}$$ by the independence   from the chosen representatives.   One can show \cite{dp3}:

\begin{enumerate}
\item The $\mathcal W(V)$ module $\Pi_X$ is semisimple of finite length.

\item The isotypic components  of  $\Pi_X$ are indexed by the {\it points of the  arrangement}.

\item  Given $e^\phi \in P(X)$, the corresponding isotypic component   $\Pi_X(\phi)$ is the direct sum of the irreducible modules
$\Pi_{\underline b,\phi}$  generated by the classes $\omega_{\underline b,\phi}$ and indexed by the no broken bases extracted from $X_\phi$.

\item As a module over the ring  $S[U]$  of differential operators on $V$ with constant coefficients, $\Pi_{\underline b,\phi}$ is  free of rank
1  with generator the class $\omega_{\underline b,\phi}$.\end{enumerate}
Although we do not want to reproduce here all the details of         \cite{dp3}, it is not too difficult to reconstruct all these statements, in particular by using the ideas of section                \ref{vtp}.

      In fact, since $\mathcal W(V)=S[U]\otimes \mathbb C[\Lambda]$, we have:
$$\Pi_{\underline b,\phi}=\mathcal W(V)\otimes_{ \mathbb C[\Lambda]} \omega_{\underline b,\phi}=
S[U]  \omega_{\underline b,\phi}$$

Theorem \ref{LappTr} tells us that we can transport our filtration on $S_X$  to one on $\mathcal L_X$. By formula (\ref{DisLap}), in filtration
degree $\leq k$  we have those distributions which are supported in a finite number of translates of the sets $\Lambda\cap C(A)$ where $A$ spans
a lattice of rank $\leq k$.

\begin{remark} In (\cite{dp3}) there is a similar explicit description of  the piece $S_{X,k}/S_{X,k-1}$, for each $k\geq 0$, using the varieties  of the arrangement,   of codimension $k$.   The proofs are in the same lines as the ones we developed in the simpler case  of hyperplanes with Proposition \ref{ilrema}.
\end{remark}
\subsection{Local expression for $\mathcal T_X$\label {locnilp}}

Let us now consider the element $v_X$ class of the generating function $\prod_{x\in X}(1-e^x)^{-1}$ in $ \Pi_X$. Decompose it uniquely as  a sum
of  elements $v_{X_\phi}$  in $\Pi_X(\phi)$. By what we have seen in the previous section,  each one of these elements is expressed as   a sum
$$v_{X_\phi}=\sum_{\underline b\in\mathcal{NB}_{X_\phi}} \mathfrak q_{\underline b,\phi}\omega_{\underline b,\phi}$$ for suitable polynomials
$\mathfrak q_{\underline b,\phi}$.  Thus: \begin{equation}\label{ilv}  v_X=\sum_{\phi\in \tilde P(X)} \sum_{\underline b\in
\mathcal{NB}_{X_\phi}\ |\ \mathfrak c\subset C(\underline b)} \mathfrak q_{\underline b,\phi}\omega_{\underline b,\phi}\end{equation}

We are now ready to state but,  not yet fully justify, the main formula that one can effectively use for computing the  partition function
$P_X$:
\begin{theorem} \label{pexxpr}  Given a point $x$ in the closure of a big cell $\mathfrak c$ we have
\begin{equation}\label{JKP}  P_X(x )=\sum_{\phi\in \tilde P(X)} e^\phi\sum_{\underline b\in \mathcal{NB}_{X_\phi}\ |\ \mathfrak c\subset C(\underline b)}\mathfrak q_{\underline b,\phi}( -x).\end{equation}
\end{theorem}
\proof Using Proposition \ref{ltd} and  formulas (\ref{laparto}) and (\ref{ilv}), we deduce that the right handside of (\ref{JKP})  coincides
with the partition function on the cone $C(X)$ minus possibly a finite number of translates of lower dimensional cones.

In particular we deduce that  the partition function $P_X$ coincides locally with  a quasi polynomial for the lattice $\Lambda$ on the cone
$C(X)$ minus possibly a finite number of translates of lower dimensional cones. Also, since to get this result we have only used the fact that
each vector in $X$ lies in $\Lambda$, we can substitute $\Lambda$ with the  larger lattice $\Lambda/n$ and we also get that the partition
function also coincides locally with  a quasi polynomial for the lattice $\Lambda/n$ on the cone $C(X)$ minus possibly a finite number of
translates of lower dimensional cones.  By the continuity of quasipolynomials we deduce that  our Theorem will follow from the following
Proposition which will be proved in    section  \ref{lespa}:
\begin{proposition}\label{lariduz}  On the closure of  each big cell $\mathfrak c$, the partition function  $P_X$ coincides with a quasi
polynomial for some lattice $\Lambda/n$.

\end{proposition}

\section{Partial fractions}
\subsection{Algebraic identities } We follow the approach of Szenes and Vergne,
used in \cite{SV} (and also in \cite{dp3})  to prove Proposition \ref{lariduz} and hence finish the proof of Theorem  \ref{pexxpr}.

We do not need a very fine analysis, let us consider the set $\tilde X$  formed by all the positive rational multiples of the vectors in $X$.

Observe that all the cones generated by subsets of elements of  $\tilde X$ coincide with the ones  generated by subsets   of  $X$. In other
words the geometry of $C(X)$ does not depend on the denominators that we will introduce.\smallskip

The following lemma is a simple computation
\begin{lemma}\label{unavariab} Let $\underline b$ be a set of linearly independent vectors in a lattice $M$, $e^{\phi}$ a character of the lattice
spanned by $\underline b$, $n_c$ a positive integer for all $c\in \underline b$. If we expand
$$\prod_{c\in \underline b}{(1-e^{\langle\phi|c\rangle-c})^{- n_c}=\sum_{y\in M\cap C(\underline b)}Q(y)e^{-y}},$$
we have for $y=\sum_{c\in \underline b}k_cc$, $k_c\geq 0$ for all $c\in\underline b$
$$Q(y)=\prod_{c\in  \underline b}e^{\langle\phi|y\rangle}\binom{ n_c+k_c-1  }{k_c}.$$
\end{lemma}

Notice that this lemma give immediately Proposition \ref {lariduz} in the special case in which $X$ consists of  the elements of a basis
(contained in the lattice $\Lambda$) each repeated any number of times.

Our next task is to reduce to this case. For this we need to show

\begin{proposition}\label{partialfra} Let $\gamma=\sum_{a\in X}r_aa$ with $r_a\in \mathbb Q$, $0\leq r_a\leq 1$, then the function
$e^{-\gamma}\prod_{a\in X}   (1-e^{\langle\psi|c\rangle-a})^{-1} $, with $e^\psi$ a character, can be written as a linear combination with
constant coefficients of elements of   the form $$\prod_{c\in \underline b}{(1-e^{\langle\phi|c\rangle-c})^{- n_c}}$$ with $\underline b$ a
linearly independent set of elements in $\tilde X$. Furthermore if $e^\psi$ is of finite order, each of the   $e^{\phi}$ has finite order on the
lattice spanned by $\underline b$.
\end{proposition}

Before giving the proof of this Proposition we need to show some simple identities. \begin{lemma}\label{rela4} \begin{equation}{n\over
1-x^n}=\sum_{i=0}^n{1\over 1-\zeta^ix},\quad \zeta:=e^{2\pi i/n}.\end{equation}
\end{lemma} \proof Take an auxiliary variable    $t$ and the logarithmic derivative (relative to   $t$) $$ {nt^{n-1}\over
t^n-x^n}dt=d\log(t^n-x^n)=d\log(\prod_{i=0}^{n-1}(t-\zeta^ix))$$
$$=\sum_{i=0}^{n-1}d\log( t-\zeta^ix))=\sum_{i=0}^{n-1}{1\over ( t-\zeta^ix)}dt.$$ Next set  $t=1$ in the coefficients of $dt$. \qed

\begin{lemma}\label{somma}
\begin{equation}
1-\prod_{i=1}^r z_i=\sum_{\emptyset\subsetneq I\subset \{1,\ldots ,r\}}\prod_{i\in I}(-1)^{|I|+1}(1-z_i)
\end{equation}
\end{lemma}
\proof The proof is by induction on $r$. The case $r=1$ is clear. In general, using the inductive hypothesis, we have
\begin{displaymath}
\sum_{\emptyset\subsetneq I\subset \{1,\ldots ,r\}}\prod_{i\in I}(-1)^{|I|+1}(1-z_i)= \sum_{\emptyset\subsetneq I\subset \{1,\ldots
,r-1\}}\prod_{i\in I}(-1)^{|I|+1}(1-z_i)-
\end{displaymath}
\begin{displaymath}
\sum_{I\subset \{1,\ldots ,r-1\}}\prod_{i\in I}(-1)^{|I|+1}(1-z_i)(1-z_r) = 1-\prod_{i=1}^{r-1} z_i \\ + (1-(1-\prod_{i=1}^{r-1}z_i
))(1-z_r)=\end{displaymath}
\begin{displaymath}
=1-\prod_{i=1}^r z_i
\end{displaymath}
\qed

A variation of this formula is the following:
\begin{lemma}\label{somma1} Set $t= \prod_{i=1}^h z_i\prod_{i=h+1}^rz_i^{-1}$, $a\in\mathbb C^*$
\begin{displaymath}
1-t=\sum_{\emptyset\subsetneq I\subset \{1,\ldots ,h\}}\prod_{i\in I}(-1)^{|I|+1}(1-z_i)-a^{-1}\sum_{\emptyset\subsetneq I\subset \{h+1,\ldots
,r\}}\prod_{i\in I}(-1)^{|I|+1}(1-z_i)+\end{displaymath}\begin{equation}\label {cavol} a^{-1}\sum_{\emptyset\subsetneq I\subset \{h+1,\ldots
,r\}}\prod_{i\in I}(-1)^{|I|+1}(1-z_i)(1-at)
\end{equation}
\end{lemma}
\proof This is immediate from the previous Lemma once we remark that $ 1-t =1-\prod_{i=1}^hz_i-t(1-\prod_{i=h+1}^rz_i$).\qed

From this  we get,
\begin{lemma}\label{quozie} Set $t= \prod_{i=1}^h z_i\prod_{i=h+1}^rz_i^{-1}$ with $0\leq h\leq r$. Then
$${1\over \prod_{i=1}^r(1-z_i)}=\sum_{\emptyset\subsetneq I\subset \{1,\ldots ,h\}}{(-1)^{|I|+1}\over (1-t)\prod_{i\notin I}(1-z_i)}-$$\begin{equation}\label{quoz}-\sum_{\emptyset\subsetneq I\subset \{h+1,\ldots ,r\}}({(-1)^{|I|+1}\over\prod_{i\notin I}(1-z_i)}-{(-1)^{|I|+1}\over (1-t)\prod_{i\notin I}(1-z_i)}).\end{equation}
If $a\in\mathbb C^*$ and $ a\neq 1$
$${1\over \prod_{i=1}^r(1-z_i)(1-at)}={a\over a-1}(\sum_{\emptyset\subsetneq I\subset \{1,\ldots ,h\}}{(-1)^{|I|+1}\over (1-at)\prod_{i\notin I}(1-z_i)})-$$\begin{equation}\label{quoz2}-{1\over a-1}(\sum_{\emptyset\subsetneq I\subset \{h+1,\ldots ,r\}}({(-1)^{|I|+1}\over\prod_{i\notin I}(1-z_i)}-{(-1)^{|I|+1}\over (1-at)\prod_{i\notin I}(1-z_i)})+{1\over \prod_{i=1}^r(1-z_i)}).\end{equation}
\end{lemma}
\proof The first relation follows from (\ref{cavol}) dividing by  $(1-t)(1-z_1)(1-z_2)\cdots (1-z_r)$ and taking $a=1$.

The second writing
$$1={a\over a-1}(1-t)-{1\over a-1}(1-at)$$
and then dividing by  $(1-at)(1-z_1)(1-z_2)\cdots (1-z_r)$.\qed

\begin{lemma}\label{expand2}  Given elements $a_i \in \tilde X$,   $b=\sum_{i=1}^M(n_i/n) a_i$  $0\leq n_j\leq n$, and a character $e^\eta$ of the
lattice spanned by the $a_i$, we can write the element
\begin{equation}\label{produ}{e^b \over\prod_{i=1}^M(1-e^{\langle\eta|a_i\rangle-a_i})}=\prod_{i=1}^M{e^{(n_i/n) a_i}
\over(1-e^{\langle\eta|a_i\rangle -a_i})}\end{equation} as a linear combination with constant coefficients of elements of the form $${1\over
(1-e^{\langle\zeta|c_1\rangle -c_1}) \cdots (1-e^{\langle\zeta|c_r\rangle -c_r}) }$$ where  the $c_i\in \tilde X$ and $e^\zeta$ is a character
of the lattice spanned by the $c_i$. Furthermore if $e^\eta$ is of finite order, also each $e^\zeta$ is of finite order.
\end{lemma}
\proof   Start from the identity
$${e^{-a}\over 1-e^{\langle\eta|a\rangle-a}}={e^{-\langle\eta|a\rangle}-e^{-\langle\eta|a\rangle}(1-e^{\langle\eta|a\rangle-a})
\over 1-e^{\langle\eta|a\rangle-a}}=  {e^{-\langle\eta|a\rangle}\over 1-e^{\langle\eta|a\rangle-a}}-e^{-\langle\eta|a\rangle}$$ from which we
get that, if $k\leq n$ and $e^{\zeta_i}, i=1,\dots, n$ are the $n^{th}$ roots of $e^\eta$:
$${e^{-ka/n}\over 1-e^{\langle\eta|a\rangle-a}}=\prod_{i=1}^k{e^{-a/n}\over 1-e^{\langle\zeta_i|a\rangle-a/n}}\prod_{i=k+1}^n{1\over 1-e^{\langle\zeta_i|a\rangle-a/n}} $$
can be expanded as a  linear combination of products of elements $( 1- e^{\langle\zeta_i|a\rangle-a/n})^{-1}$.  We now apply this procedure to
each term of our product. The fact that if $e^\eta$ is of finite order also the $e^\zeta$ are of finite order follows immediately from their
definition.\qed

\subsection{The main expansions\label{lespa}}
We can now  give

\proof(of Proposition  \ref{partialfra}) From Lemma \ref{expand2} we deduce that we can assume that $\gamma=0$.

Also there is nothing to prove unless we have a linear dependency $na_0=\sum_i^kn_ia_i, $  $n_i\in\mathbb Z$  which we write as  $$b_0:={
a_0\over|\prod_in_i|}=\sum_i\pm {a_i \over |n\prod_{j\neq i}n_j|}=\sum\pm b_i. $$ We substitute now  $a_0$ with $b_0:={ a_0/ |\prod_in_i|}$ and
$a_i$ with $b_i:={a_i / |n\prod_{j\neq i}n_j|} $. Next apply  lemma  (\ref {rela4}), using suitable  characters $e^\zeta$. We get a new
expansion in which the $b_i$ appear in place of the $a_i$.

For each product $\prod_{i=0}^k{( 1-e^{\langle\zeta_i|b_i\rangle-b_i})^{-1}}$ we have
two possibilities. If the equations  $1-e^{\langle\zeta_i|b_i\rangle-b_i}=0$ are incompatible we are in position to apply the formula
(\ref{quoz2})  to $x_i=\zeta_ie^{b_i}, i\geq 1$ and to $t=\zeta_0e^{b_0}$.  If they are compatible we apply formula (\ref{quoz}). In both cases
we substitute the product with a linear combination of products with less terms in the denominator. We then proceed by induction. \qed

We can now finish the Proof of Proposition \ref{lariduz}. Using Proposition  \ref{partialfra} and Lemma \ref{unavariab}, we get that for any
$\gamma=\sum_{a\in X}r_aa$ with $r_a\in \mathbb Q$, $0\leq r_a\leq 1$,  the function $e^{-\gamma}\prod_{a\in X}   (1-e^{-a})^{-1} $ is a
quasipolynomial in  the interior of the big cells. From this we deduce that
$${1\over\prod_{a\in X}   (1-e^{-a})} =e^\gamma {e^{-\gamma}\over\prod_{a\in X}   (1-e^{-a})^{-1}}$$ is a quasipolynomial on the big cells
translated by $\gamma$. Since as $\gamma$ varies these translates cover each big cell, our claim follows.\qed
\subsection{$E-$splines}

In our setting we can also consider the  Euler--Maclaurin sums (cf. \cite{BV4},\cite{BV5}):

$$\mathcal P_{X,\underline\mu}(v):=\sum_{(n_1,\dots,n_N)\in P(v)} e^{-\sum n_i\mu_i},$$ $$ P(v):=\{ (n_1,\dots,n_N)\,|\, \sum n_ia_i=v,\ n_i\in\mathbb N\} .$$

\bigskip

$\mathcal P_{X,\underline\mu}(v)$  must   be understood as the distribution   $\sum_{v\in \Lambda\cup C(X)}  \mathcal
P_{X,\underline\mu}(v)\delta_v$ supported at the points of the lattice $\Lambda\cup C(X)$.

Its Laplace transform is
$$\sum_{v\in  \Lambda\cup C(X)}\mathcal P_{X,\underline\mu}(v)e^{-\langle x\,|\, v\rangle}=\prod_{a\in X} (1-e^{-a-\mu_a})^{-1}.$$
This leads to a similar theory in which now in general there are more points in the arrangement, in the generic case $\delta(X)$  (see
definition \ref{delta}), distinct points, and correspondingly more exponential functions appear in the final formulas.  We leave the details to
the reader and also refer to section \ref{diffe}, where we approach the same topic from a different perspective.

\smallskip

A particularly interesting case is when the numbers $e^{\mu_a}=\nu_a=\chi(a)$ for a character  $\chi$ of $\Lambda$, in particular for a
character  of finite order which is 1 on a sublattice $\Lambda_0$ of finite index. Then we have that, if $\sum _{a\in X} n_aa=v$ we have
$e^{-\sum_{a\in X} n_a\mu_a}=\overline \chi(v)$.  In this case we denote $\mathcal P_{X,\underline\mu}= \mathcal P_{X,\chi}$.  Thus we  have
that
$$\mathcal P_{X,\chi}=\sum_{v\in \Lambda\cup C(X)} \overline\chi(v) \mathcal P_{X }(v)\delta_v$$  and
$$\mathcal P_{X }|_{\Lambda_0}=|\Lambda/ \Lambda_0|^{-1}\sum_{\chi\in \widehat{ \Lambda/ \Lambda_0}}\mathcal P_{X,\chi}$$
where $\mathcal P_{X }|_{\Lambda_0}$ is the restriction of $\mathcal P_{X }$ to ${\Lambda_0}$.  In a similar way, by applying translations we
can restrict to the other cosets.
\section{Difference equations \label{difequ}}

\subsection{Difference operators}

Let us consider the space of complex valued functions $f$ on the set $\Lambda$.  We identify this space with the {\it  algebraic dual} $\mathbb
C[\Lambda]^*$ of $\mathbb C[\Lambda]$  by the formula:
$$\langle f\,|\,e^a\rangle:=f(a).$$
For $ v\in\Lambda$ we define the   difference operator $\nabla_v$, acting on functions $f\in \mathbb C[\Lambda]^*$ as:
$$\nabla_vf(u):=f(u)-f(u-v).$$
Parallel    to the study of $D(X)$, we  start by studying   the  {\it difference equations}  $\nabla_{Y}f=0$, where $\nabla_{ Y}:=\prod_{v\in
Y}\nabla_v$ as $Y\in \mathcal E(X)$ runs over the cocircuits.\smallskip

Let us denote the space of solutions by:
\begin{equation} \nabla(X) :=\{ f:\Lambda\to\mathbb C,\ |\ \nabla_{Y}(f)=0,\ \forall Y\in \mathcal E(X)\}.\end{equation}
As we shall see, this space is not only a formal construct, but it plays an essential role in the theory of partition functions.\medskip

We want to reformulate the fact that a function is a solution of a system of  difference equations as the property for such a function to vanish
on an appropriate ideal $J_X$  of $\mathbb C[\Lambda]$.\smallskip

Notice that, the ideal $I_1$ of functions in $\mathbb C[\Lambda]$ vanishing at $1\in T$ has as linear basis the elements  $1-e^{-a},\ a\in
\Lambda,\ a\neq 0$.  If one takes another point $ e^{\phi}$  the ideal $I_{\phi}$ of functions in $\mathbb C[\Lambda]$ vanishing at $e^\phi$,
has as linear basis the elements  $1-e^{-a+\langle a\,|\,\phi\rangle},\ a\in \Lambda,\ a\neq 0$.\medskip

In fact, for every $a,x\in \Lambda$  we have:
$$\langle  \nabla_af\,|\, e^x\rangle=\langle   f\,|\, (1-e^{-a}) e^x\rangle.  $$ Thus the difference operator $\nabla_a$ is  the dual of the multiplication operator  by $1-e^{-a}$. In this setting, we get a statement analogous to that of Theorem \ref{soleqdiff} on differential equations:
\begin{theorem}
A  function  $f$ on $\Lambda $ satisfies the difference equation

$p(\nabla_{a_1},\dots,\nabla_{a_k})f=0$  with $p$ a polynomial if and only if,  thought of as   element of the dual of $\mathbb C[\Lambda]$, $f$
vanishes on the ideal  of $\mathbb C[\Lambda]$ generated by   the  element $p(1-e^{-a_1},\dots,1-e^{-a_k})$.
\end{theorem}
We have also the {\it twisted difference operators}  $\nabla_a^\phi,\ e^\phi\in T$   defined by $(\nabla_a^\phi f)(x):=f(x)-e^{\langle
a\,|\,\phi\rangle}f(x-a)$, dual to multiplication by $1-e^{-a+\langle a\,|\,\phi\rangle}$.

We need some simple algebraic considerations on subschemes of $T$ supported at a point $e^{\phi}\in T$. Let $S=:\mathbb C[\Lambda]/J$ be the
(finite dimensional) coordinate ring of such a scheme. \smallskip

In this slightly more general setting, we can repeat the discussion made in section \ref{locnilp}. Given an element $a\in \Lambda$,  we have
that $ e^{-\langle a\,|\,\phi\rangle}e^a-1$  is 0 in $\phi$. Thus, in the ring $S$, the class of $e^{-\langle a\,|\,\phi\rangle}e^a-1$ is
nilpotent. The power series of $\log(1+t)=t-t/2+t^2/3-\dots $  can be computed on $t=e^{-\langle a\,|\,\phi\rangle}e^a-1$ and we defines
$\underline a:=\log(1+(e^{-\langle a\,|\,\phi\rangle}e^a-1))+\langle a\,|\,\phi\rangle$.

We clearly have:
\begin{lemma}  The map $i:a\to \underline a$ is a linear map.

$\tilde a:=\underline a-\langle a\,|\,\phi\rangle$ is nilpotent.

$e^{\underline a}=\exp(\tilde a)e^{\langle a\,|\,\phi\rangle}\in S$  equals the class of $e^a$.

$i$ extends to a linear map $i:V\to S$ and then to a  surjective homomorphism $\overline i:S[V]\to S$.
\end{lemma}
\begin{proof}
The class $e^{-\langle a\,|\,\phi\rangle+a}=\exp(\tilde a)=\sum_{k\geq 0}{\tilde a^k / k!}$ reduces to a finite sum lying in the image of  the
homomorphism $\overline i$ which is therefore surjective.
\end{proof}\medskip

We have thus identified $S$ with a quotient $S[V]/I$ by an ideal  of finite codimension in $S[V]$.

\smallskip

The fact   that, for each $a\in\Lambda$ we have $\underline a-\langle a\,|\,\phi\rangle$ is nilpotent, means that $I$ defines a unique point in
which of course $\underline a=\langle a\,|\,\phi\rangle$. That is $I$ defines the point $\phi$  (notice that the choice of $\phi$ is not
unique).\medskip

We thus can identify,  using this algebraic logarithm, the given scheme as a subscheme in the tangent space. This proves the following:

\begin{proposition}\label{diffe1} Let $J\subset $ be an ideal  such that $ \mathbb C[\Lambda]/J $ defines (set theoretically) the point $e^{\phi}$. Under the logarithm isomorphism    $ \mathbb C[\Lambda]/J  $ becomes isomorphic to  a ring $ S[V]/I $ defining the point $\phi$.\smallskip

In particular we get a canonical isomorphism of  the space of solutions of the difference equations given by $J$   with the space of solutions
of the differential equations  given by $I$.

The explicit formula is,  given $f\in S^*= (S[V]/I )^*$:\begin{equation}\label{iddd}f(a):=\langle f\,|\,e^a\rangle=e^{ \langle a\,|\,\phi\rangle
} \langle f\,|\,\exp(\tilde a)\rangle.\end{equation} where $ \langle f\,|\,\exp(\tilde a)\rangle$ is a polynomial in $a$.\end{proposition}

A special role is played by   the points of finite order   $m$, i.e. characters $e^\phi$ on $\Lambda$ whose kernel is a sublattice
$\Lambda_\phi$ of index $m$ (of course this implies $m\phi\in\Lambda^*$).
\smallskip

As we have seen, a function $f(x)$ appearing in the dual of  $ \mathbb C[\Lambda]/J $, when $J$ defines (set theoretically) the finite order
point  $e^\phi$, is of the form    $f(x):=e^{\langle x\,|\,\phi\rangle}g(x)$ with $g(x)$ a polynomial.

Since  $e^\phi$ is of finite order, we know   $e^{\langle x\,|\,\phi\rangle}$ takes  constant values ($m-${th} roots of 1)     on the $m$ cosets
of the sublattice $\Lambda_\phi.$ Thus $f(x)$ is a polynomial only on each such coset.   This is a typical example of what is called a {\it
periodic polynomial} or {\it quasi--polynomial}.

Formally:
\begin{definition}\label{PE}  A function $f$ on a lattice   $\Lambda$ is a periodic polynomial, if there exists a sublattice $\Lambda^0$ of finite index in $\Lambda$, such that $f$, restricted to each coset of $\Lambda^0$ in $\Lambda$, is (the restriction of) a polynomial.
\end{definition}

We deduce the general:

\begin{theorem} Let $J\subset $ be an ideal  such that $ \mathbb C[\Lambda]/J $ is finite dimensional and defines (set theoretically) the points $e^{\phi_1},\dots,e^{\phi_k}$, all of finite order.

Then,     the   space   $ (\mathbb C[\Lambda]/J)^* $ of  solutions of the difference equations associated to $J$ is a direct sum of spaces of
periodic polynomials $e^{\langle a\,|\, \phi_i\rangle}p(a)$ for the points $\phi_i$, each  invariant under translations under $\Lambda$.

\end{theorem}

Although we will not need it, for completeness we now want to prove a converse. Let us thus take a finite dimensional vector space $Q$ spanned
by periodic polynomials $ e^{\langle a\,|\, \phi_i \rangle}p (a)$ (for some points  $e^{\phi_i }$ of finite order)   invariant under
translations under $\Lambda$. We need the following:
\begin{lemma}\label{ind}  Any non zero subspace $M$ of $Q$, invariant under translations under $\Lambda$, contains one of the functions $e^{\langle a\,|\, \phi_i \rangle}$.
\end{lemma}
\proof   Indeed,  a non zero subspace invariant under translations under $\Lambda$  contains a non zero common $\Lambda$-eigenvector. This
eigenvector  must necessarily be a multiple of one of the $e^{\langle a\,|\, \phi_i \rangle}$.\qed

Let $J$ be the ideal of the difference equations satisfied by $Q$.
\begin{theorem}   $J$ is of finite codimension and $Q$ is the space of all solutions of $J$.

\end{theorem}

\begin{proof} We can easily reduce to the case of  a unique point $\phi$.
In fact the commuting set of linear operators $\Lambda$ on $Q$  has a canonical Fitting decomposition  relative to the generalized eigenvalues
which are characters of $\Lambda$. Our hypotheses imply that these points are of finite order.\smallskip

For any $p\in Q$ and $b\in \Lambda$, we set  $\tau_bp=p(x+b)$. Consider   the function $j(p)$ on $\Lambda$ given by
\begin{equation}\label{deff} j(p)(b)= p(b)=  (\tau_bp)(0).\end{equation}
By definition $j(\tau_ap)(b)=j(p)(a+b)$  so $j:Q\to  \mathbb C[\Lambda]^*$ is equivariant. We have $ j( e^{\langle a\,|\, \phi \rangle})(b)=
e^{\langle b\,|\, \phi \rangle} $.   In particular  $e^{\langle a\,|\, \phi \rangle}$  is not in the kernel  of $j$. Since $j$ is equivariant,
Lemma \ref{ind}  implies that $j$   is injective, and we deduce a surjective map  $\mathbb C[\Lambda]\to Q^*$.

We give to $Q^*$ an algebra structure by defining the product of two linear forms as follows. For $f(a) =e^{\langle a\,|\, \phi \rangle}p (a)\in
Q$ we have that, by translation invariance, $f(a+b) =\sum_ie^{\langle a\,|\, \phi \rangle}p_i (a)e^{\langle b\,|\, \phi \rangle}q_i (b)$  with
$e^{\langle a\,|\, \phi \rangle}p_i (a), e^{\langle b\,|\, \phi \rangle}q_i (b)\in Q$. This can be thought of as a {\it coalgebra structure}
$\Delta:Q\to Q\otimes Q$ on $Q$.

By duality this induces an algebra structure for which $j^*$ is a surjective algebra homomorphism. Thus $Q^*$ is naturally isomorphic to
$\mathbb C[\Lambda]/J$ for some ideal $J$.
\smallskip

Now consider    the dual map $\pi:(\mathbb C[\Lambda]/J)^*\to Q$. By the first theorem   $(\mathbb C[\Lambda]/J)^*$ is identified to a space of
periodic polynomials. This identification coincides with the one given by $\pi$, due to  the definition    of $j$. \end{proof}

\bigskip
\subsection{The difference Theorem\label{diffe}}

We go back to the equations $\nabla_{Y}f=0$. As in Section \ref{remarkable},  we can generalize our problem by considering, for any sequence
$\underline \nu=\{\nu_v|\ v\in X\}$ of  non zero complex numbers,   the products $N_{Y}(\underline \nu):=\prod_{v\in  Y}(1-\nu_ve^{-v})$ as $Y$
runs over the cocircuits, and the ideal  $J_X(\underline \nu)$ that they generate. We set $J_X:= J_X(1,\dots,1)$.\smallskip

The space of functions $\nabla(X) $, which  we want to describe, is thus identified to the dual of $\mathbb C[\Lambda]/J_X$.

It turns out that also this space is finite dimensional and we are going to show that its dimension   is the following  {\it weighted analogue}
$\delta(X)$ of $d(X)$.

We use the notations of section \ref{ta}. Given a basis $\underline b$ extracted from $X$ consider the lattice $\Lambda_{\underline
b}\subset\Lambda$ that it generates in $\Lambda$, $ \Lambda/\Lambda_{\underline b}$ is a finite group of order  $ [\Lambda:\Lambda_{\underline
b}]=|\det(\underline b)|,$ its character group  is the finite subgroup  $T( {\underline b})$ of $T$  which is the intersection fo the kernels of
the characters   $e^a$ as $a\in  \underline b .$

We now define,
\begin{equation}\label{delta}\delta(X):=\sum_{\underline b\in\mathcal B(X)}|\det(\underline b)| .\end{equation}

We have:
\begin{proposition}
\begin{equation}\sum_{e^\phi\in P(X)}  d(X_\phi)=\delta(X).\end{equation}
\end{proposition}
\proof $\delta(X)$ counts the number of pairs $e^\phi, \underline b$ such that $e^\phi\in T( {\underline b})$, or equivalently $\underline
b\subset X_\phi$.\qed

\begin{lemma}
The variety defined by the ideal $J_X$ is $P(X)$.
\end{lemma}
\proof  The proof is identical to that of lemma \ref{span} and of the first part of Theorem \ref{main2}, the only difference being the fact,
that when we extract a basis $\underline b$  from $X$ the equations $e^a-1=0$ as $a\in  \underline b  $  define $T( {\underline b})$.  By its
definition (\ref{ipunti}), $P(X)=\cup  T( {\underline b})$.\qed

Now by elementary commutative algebra we know that
$$\mathbb C[\Lambda]/J_X=\oplus_{e^\phi\in P(X)}\mathbb C[\Lambda]/J_X(\phi)$$
where $\mathbb C[\Lambda]/J_X(\phi)$ is its localization at $e^\phi$.

\begin{theorem}\label{logiso} Under the logarithm isomorphism    $\mathbb C[\Lambda]/J_X(\phi)$ becomes isomorphic to  the ring $A_{X_\phi}=S[V]/I_{X_\phi}.$ Thus:
\begin{equation} \mathbb C[\Lambda]/J_X\cong \oplus_{e^\phi\in P(X)}A_{X_\phi}.\end{equation}

In particular   $\dim(\mathbb C[\Lambda]/J_X)=\delta(X)$.

\end{theorem}

\proof  Let us see what happens to the equations $N_{ Y} =\prod_{v\in  Y}(1- e^{-v})$ as $Y$ runs over the cocircuits.

When we localize at $e^\phi$ the elements $1-e^{-v}$  where $v\notin X_\phi$ become invertible and hence can be dropped from the equations. As
for $1-e^{-v}$  where $v\in X_\phi$ we have that $1-e^{-v}=\underline v (1-\sum_{k\geq 1}{(-\underline v)^k / (k+1)!})$.

Obviously $(1-\sum_{k\geq 1}{(-\underline v)^k / (k+1)!})$ is invertible hence we can replace the equation  $N_{Y}$ by $M_{ Y}=\prod_{v\in
Y}\underline v$.  We obtain the equations defining $A_{X_\phi}$  and this proves the claim completely.

\qed

It is convenient to single out, in this direct sum, the term  relative to $e^\phi=1$ which is $D(X)$ and write
\begin{equation} \label{decc}E(X):=\oplus_{e^\phi\in P(X), \ e^\phi\neq 1}e^{\phi}D({X_\phi}),\quad  \mathbb \nabla(X)=D(X)\oplus E(X) .\end{equation}
\subsection{The parametric case}  One could treat in a similar way the parametric case.

We take a sequence $\underline \nu:=\{\nu_a\,|\, a\in X\}$.

The main difference is that now the points of the arrangement  are defined as follows.

Given a basis ${\underline b}$  extracted from $X$ instead of the finite subgroup  $T( {\underline b})$, intersection fo the kernels of the
characters   $e^a$ as $a\in  \underline b  $  we have to consider the set $T_{ \underline \nu}( {\underline b})$  where   $\nu_ae^{-a}-1=0$
this is a coset of the subgroup  $T( {\underline b})$, still consisting of  $|\det(\underline b)|$  elements.

For a generic sequence $\underline \nu$ the ideal $J_X(\underline \nu)$ defines $\delta(X)$ distinct and reduced points, while for special
values we may have less points but in each such point  a  similar space of quasi periodic polynomials.

So now we define (as for (\ref{ipunti}))
$$  P_{ \underline \nu}(X):=\cup_{\underline b\in\mathcal B(X)}T_{ \underline \nu}( {\underline b}).$$

For any $e^\phi\in P_{ \underline \nu}(X)$ set:
\begin{equation}\label{xfin}X_\phi(  \underline \nu):=\{a\in X\,|\,  e^{\langle a\,|\,\phi\rangle}=\nu_a \}.\end{equation}
We have as usual the analogues of  the results in  the non parametric case:

\begin{lemma}
The variety defined by the ideal $J_X(  \underline \nu)$ is $P_{ \underline \nu}(X)$.
\end{lemma}

$$\mathbb C[\Lambda]/J_X(  \underline \nu)=\oplus_{e^\phi\in P_{ \underline \nu}(X)}\mathbb C[\Lambda]/J_X(  \underline \nu)(\phi)$$
where $[\mathbb C[\Lambda]/J_X(  \underline \nu)](\phi)$ is its localization at $e^\phi$.

Finally
\begin{theorem} Under the logarithm isomorphism    $\mathbb C[\Lambda]/J_X(\phi)$ becomes isomorphic to  the ring $S[V]/I_{X_\phi}.$ Thus:
\begin{equation} \mathbb C[\Lambda]/J_X(  \underline \nu)\cong \oplus_{e^\phi\in P_{ \underline \nu}(X)}S[V]/I_{X_\phi}(  \underline \nu).\end{equation}

For all   $\underline \nu$, we have   $\dim(\mathbb C[\Lambda]/J_X(  \underline \nu))=\delta(X)$.

\end{theorem}

Of course also in this case we could perform the polygraph construction, now the base would be a torus $\mathbb C^N$. We still have a flat
family and the polygraph is a Cohen Macaulay variety as in Remark (\ref{rem}).   \smallskip

One can use the parametric case as follows.  Suppose that we pass, from the lattice  $\Lambda$ to a finer lattice $\Lambda/n$ for some positive
integer $n$ then, for $a\in\Lambda$ we write in $\mathbb C[\Lambda/n]$  $$1-\mu e^{-a}=\prod_{\gamma\, |\,\gamma^n=\mu}(1-\gamma e^{-a/n})$$
Thus if we take our list $X=\{a_1,\ldots ,a_N\}$ with  parameters $\underline \nu=\{\nu_1,\ldots ,\nu_N\}$ in $\Lambda$ and consider it in
$\Lambda/n$, we get an ideal $J^{(n}_X(\underline \nu)\subset \mathbb C[\Lambda/n]$.  The same ideal is also associated to the list
$$X^{(n)}=\{\underbrace{{a_1/ n},\ldots ,{a_1/ n}}_n,\ldots,\ldots,\underbrace{{a_N/ n},\ldots ,{a_N/n}}_n\}\ \text{with parameters}$$$$\underline\nu^{(n)}:=\{\gamma_1,\varepsilon\gamma_1,\ldots,\varepsilon^{n-1}\gamma_1,\ldots ,\gamma_N,\ldots  ,\varepsilon^{n-1}\gamma_N\}$$
where, for each $i=1,\ldots N$, $\gamma_i$ is an n-th root of $\nu_i$ and $\varepsilon=\exp(2\pi i/n)$.

One has to remark that $\Lambda/n$ is the character group of a torus $T^{1/n}$  which maps surjectively  $\pi:T^{1/n}\to T$ to $T$, with kernel
$K_n$ the character group of $(\Lambda/n)/\Lambda$.

The points of the arrangement associated to $X$ in $T^{1/n}$ are $\pi^{-1}(P(X))$ and a union of cosets  of $K_n$.

The algebra  $\mathbb C[\Lambda/n]/J_X(\underline\nu) \mathbb C[\Lambda/n]= \mathbb C[\Lambda/n]/J_{X^{(n)}}(\underline\nu^{(n)})$. Thus even if
we start from the trivial parameters $=1$, once we pass to a finer lattice we find as parameters the roots of 1.

\subsection{A realization of $C[\Lambda]/J_X$}

As for the box spline the functions in $D(X)$ play a basic role, so do the functions in $ \nabla(X) $   in the theory of linear diophantine
equations. In fact they  describe combinatorially the {\it partition  function} $\mathcal P_X (v)$ associated to $X$.

We go back to the space  $\Pi_X:=S_X/S_{X,n-1}$  of {\it  polar parts}. This we have decomposed, as $\Pi_X=\oplus_{\phi\in P(X)}
P_{X_\phi}(\underline \phi) $, through the points $\phi\in P(X)$  into local modules $P_{X_\phi}(\phi)$  of polar parts for the affine
arrangement centered at $\phi$.

The element $v_X$ class of the generating function $\prod_{a\in X}(1-e^{-a})^{-1}$ in $ \Pi_X$ decomposes into a sum of local elements
$v_{X_\phi}$.

We consider next the submodule $\mathcal Q_X$ generated in $\Pi_X$ by $v_X$ under $\mathbb C[\Lambda]$ and we deduce, reducing to Theorem
(\ref{main1}) and using Theorem \ref{logiso} that:

\begin{theorem}\label{main3}
The annihilator of $v_X$ is the ideal   $J_X$ generated by the elements $N_{Y}=\prod_{a\in  Y}(1-e^{-a}),$ as $Y$ runs over the cocircuits. Thus
$\mathcal Q_X=\mathbb C[\Lambda]/J_X$.

We obtain a canonical commutative diagram, made  of isomorphisms,  compatible with all the identifications made:
\[\begin{CD}\mathbb C[\Lambda]/J_X @>\cong>>\oplus_{\phi\in P(X)} A_{X_\phi}\\
@V\cong VV @V\cong VV\\\mathcal Q_X @>\cong>>\oplus_{\phi\in P(X)} Q_{X_\phi}
\end{CD}\]

\end{theorem}

\subsection{From volumes to partition functions\label{vtp}}We apply the previous theory  to the partition function   $\mathcal P_X(a)$.

Consider the class of  a product $\prod_{a\in X}(1-e^{-a})^{-1}$ in $\Pi_X$. For each point $\phi$ of the arrangement we project to the
corresponding isotypic component $\Pi_X(\phi)$ which is generated by the elements $\omega_{\underline b,\phi}$.

Since $\omega_{\underline b,\phi}$  is an eigenvector for $\Lambda$ of eigenvalue $e^\phi$, using the commutation relations with the
derivatives, we clearly have that, for all  $a\in\Lambda$, the element  $e^a-e^{\langle \phi\,|\,  a\rangle}$ acts locally nilpotently on   the
modules $\Pi_X(\phi)$ .    It follows that the action of $e^a$ on $\Pi_X(\phi)$ is  given by an operator of the form $e^{\langle a\,|\,
\phi\rangle}  u_a$ with $u_a=e^{\tilde a}$ with $\tilde a$  locally nilpotent.

The map   $a\mapsto \tilde a+\langle a\,|\, \phi\rangle$ is linear.  It induces a linear map of    $V=\Lambda\otimes_{\mathbb Z}\mathbb C$  into
linear operators on $\Pi_X(\phi)$ whose image consists of mutually commuting operators.   In this way we obtain  a $S[V]$ module structure on
$\Pi_X(\phi)$. The commutation relations with the operators of $S[U]$    then give us an action of the usual Weyl algebra   $W(U)$.

\begin{proposition}\label{voltopar} As a  $W(U)$ module  $\Pi_X(\phi)$  is  isomorphic to the module of polar parts for the list $X_\phi$  with parameters
$\underline \phi:=\{\phi_a=-\langle \phi\,|\,  a\rangle\}$ for each $a\in X_{\phi}$ by an isomorphism $j_{\phi}$ characterized by
$$j_{\phi}(\omega_{\underline b,\phi})=u_{\underline b}$$
for each no broken basis extracted from $X_{\phi}$.
\end{proposition}
\begin{proof}
The fact that the class  $\omega_{\underline b,\phi}$ is an eigenvector is an eigenvector for $\Lambda$ of eigenvalue $e^\phi$ and the
definition of the  $S[V]$ module structure imply that $\omega_{\underline b,\phi}$ is an eigenverctor for $V$ of eigenvalue $\phi$. Thus we get
a surjective map of $W(U)$ modules
$$\gamma_{\phi}:P_{X_\phi}(\underline \phi)\to (\Pi_X)_{\phi}$$
with $\gamma_{\phi}(u_{\underline b})=\omega_{\underline b,\phi}$ for each no broken basis in $X_{\phi}$. Since both modules are free
$S[U]$-modules and $\gamma_\phi$ takes a basis for $P_{X_\phi}(\underline \phi)$ to a basis for $ (\Pi_X)_{\phi}$ we deduce that it is an
isomorphism. $j_\phi$ is the inverse of $\gamma_\phi$.
\end{proof}

  By the previous discussion $e^a$ acts
on $\Pi_X(\phi)$ by the operator $e^{\langle \phi\,|\,  a\rangle}e^{\tilde a}$ with $\tilde a$ locally nilpotent.  We distinguish two cases. If
$e^{\langle \phi\,|\,  a\rangle}\neq 1$, then $(1-e^{-a})^{-1}$ gives an invertible operator in $\Pi_X(\phi)$.

If on the other hand,   $e^{\langle \phi\,|\,  a\rangle}=1$, the operator ${\tilde a/ (1-e^{-\tilde a})}$ is invertible (its power series is
written in term of Bernoulli numbers). Thus the image under $j_\phi$ of the component of $\prod_{a\in X}( 1-e^{-a})^{-1}$ equals
$$\!\prod_{a\notin X_\phi}{1\over 1-e^{-\tilde a -\langle a\,|\,\phi\rangle}}\!\prod_{a\in X_\phi}
{\tilde a\over 1-e^{-\tilde a}}(a- \langle \phi\,|\,  a\rangle) ^{-1}\!=  Q_\phi\prod_{a\in X_\phi}\! (a-\langle \phi\,|\,  a\rangle)^{-1}$$
where $$  Q_\phi= \prod_{a\notin X_\phi}{1\over 1-e^{-\tilde a -\langle \phi\,|\,  a\rangle}}\prod_{a\in X_\phi} {\tilde a\over 1-e^{-\tilde a}}
$$ is an invertible operator, on  $P_X(\phi)$, expressed as a power series in the operators $\tilde a$ and hence locally given by a
polynomial.  We deduce
\begin{proposition}\label{laide}
Under the isomorphism $j$ of coordinates $j_\phi$: $$j:\oplus_{\phi\in P(X)}\Pi_X(\phi)\to \oplus_{\phi\in P(X)}P_{X_\phi}(\underline\phi)$$ we
have the transformation
\begin{equation}\label{laidenti}j:\prod_{a\in X}( 1-e^{-a})^{-1}\mapsto  \sum_{\phi\in P(X)}Q_\phi\prod_{a\in X_\phi}\! (a-\langle \phi\,|\,  a\rangle)^{-1}.\end{equation}
\end{proposition}

We now apply Laplace transform and deduce the final general formula expressing the partition function as a sum of transforms of the local
multivariate splines.  In order to justify our results let us make some remarks. We have made an identification of  $\Pi_X(\phi)$  with
$P_{X_\phi}(\underline\phi)$ as $W(U)$ modules.

In particular  we have that  $\Pi_X(\phi)$ is a free module, over the algebra of differential operators with polynomial coefficients,  in the
classes of the elements  $  {e(\phi)}{\prod_{a\in \underline b} (1-e^{-a})^{-1}}$  while $\Pi_X(\phi)$ is a free module, over the same algebra,
in the classes of  $ {\prod_{a\in \underline b} (a-\langle \phi\,|\,  a\rangle)^{-1}}$.

Take then $M_1,M_2$ to be the two free submodules, over the algebra of differential operators with polynomial coefficients, of $S_X$ re.  $R_X$
generated by these elements.   Define by $j:M_2\to M_1$ to be the module isomorphism mapping $ {\prod_{a\in \underline b} (a-\langle \phi\,|\,
a\rangle)^{-1}}$  to $  {e(\phi)}{\prod_{a\in \underline b} (1-e^{-a})^{-1}}$.

From proposition  \ref{ltd}   we see that $M_1$  is the span of the Laplace transforms of distributions of type
$$ \sum_{\phi\in P(X),\underline b\in X_\phi} p_{\underline b}(x) |\det(\underline b)|^{-1}\sum_{v\in C(\underline b)\cap \Lambda}\phi(v)\delta_v  $$  while $M_2$, by (\ref{basic}) and   (\ref{laplaT}),  is the span of the Laplace transforms of distributions of type  $$\sum_{\phi\in P(X),\underline b\in X_\phi} p_{\underline b}(x) |\det(\underline b)|^{-1}e^{\langle\phi\,|\, x\rangle}\chi_{C(\underline b)}.$$

If we restrict both types of distributions on the set of regular points of $C(X)$ we see that, $M_2$ can be identified with a space of smooth
functions, which are locally linear combinations of polynomials times exponentials, while   $M_2$  can be identified  to functions of the
regular points in $\Lambda$ which are locally quasi polynomials.

\begin{theorem}\label{restr} Under these identifications the map $j$ consists simply into restricting the functions to the points in $\Lambda$.
\end{theorem}
\begin{proof}
Since at the level of functions the map is a linear isomoprhism compatible with the multiplication by the polynomials $S[V]$, it is enough to
verify the statement on the generator $ \sum_{v\in C(\underline b)\cap \Lambda}\phi(v)\delta_v  $  and  $e^{\langle\phi\,|\,
x\rangle}\chi_{C(\underline b)} $ for which it is clear.
\end{proof}

We deduce: \begin{theorem} \label{toddo} On the intersection of $\Lambda$ with the open set of regular points we have:
\begin{equation} \label{conclu} \mathcal P_X=\sum_{\phi\in P(X)}   \widehat {Q}_\phi  T_{X_\phi,\underline \phi}\end{equation}
\end{theorem}
\begin{proof} The explicit formula is a consequence of formula (\ref{laidenti}) plus the previous discussion which implies that the two sides coincide, since both functions are restrictions of quasi polynomials.\end{proof}

From this formula we can deduce one valid everywhere using the method of Jeffrey--Kirwan residues.

\begin{theorem} \label{pexxpr2}  Given a point $x$ in the closure of a big cell $\mathfrak c$ we have
\begin{equation}\label{JKPd}  P_X(x )=\sum_{\phi\in \tilde P(X)} \widehat {Q}_\phi\big(e^\phi\sum_{\underline b\in \mathcal{NB}_{X_\phi}\ |\ \mathfrak c\subset C(\underline b)} \mathfrak p_{\underline b,X_\phi}( -x)\big).\end{equation}
\end{theorem}
\begin{proof}
Both terms of the equality are continuous on the closure of $\mathfrak c$, they coincide in the interior by the previous proposition, hence they
are equal.
\end{proof}

\section{Some Applications}

In this short section we would like to give a streamlined presentation of some of the applications. All the results are taken from the papers of
Dahmen and Micchelli, or from \cite{dhr}, with minor variations of the proofs.
\smallskip

For further details and more information the reader should look at the original literature.

The main steps are the following.

\subsection{Discrete convolution}

Let us first  analyse the  discrete convolution:
$$  B_X*_dp=\sum_{\lambda\in \Lambda}  B_X(x+\lambda)p(-\lambda)$$ and prove that:
\begin{theorem} When $p\in D(X)$ we have that also $  B_X*_dp\in  D(X)$.

This defines a linear isomorphism  $F$ of $D(X)$ to itself, given explicitly by the invertible differential operator $F_X:=\prod_{a\in X}
{1-e^{- D_a}\over D_a}$.\end{theorem}

\proof A way to understand this convolution is by applying the Poisson summation formula to the function of $y$,   $B_X(x+y)p(-y)$.    Its
Laplace transform is obtained from the Laplace transform $e^{x}\prod_{a\in X} (1-e^{-a})/ a$ of $B_X(x+y)$ by applying the polynomial  $\hat
p(x)$ as differential operator.

In our definition of Laplace transform we have
$$Lf(\xi)= (2\pi)^{n/2}\hat f(i\xi)$$
where $\hat f$ denotes the usual Fourier transform. So  Poisson summation formula gives:
$$\sum_{\mu} L\phi(\mu)=\sum_{\lambda}  \phi(\lambda)$$
where $\mu$ runs in the {\it dual lattice} $\Lambda^*$,  of elements for which $\langle \mu\,|\,\lambda\rangle\in 2\pi i\mathbb Z,$ $ \forall
\lambda\in\Lambda.$

Thus if we are in the situation that $L\phi(\mu)=0, \forall \mu\neq 0,  \ \mu\in \Lambda^*$ we have

$$  L\phi(0)=\sum_{\lambda}  \phi(\lambda)$$

The main observation of  Dahmen and Micchelli  is that:
\begin{lemma} The Laplace    transform of $B_X(x+y)p(-y)$ vanishes  at all points $\mu\neq 0,   \ \mu\in \Lambda^*.$
\end{lemma}
\proof  We may assume that $p(x)$ is homogeneous of some degree  $k$.   The evaluation of $p(x)$ against $\prod_{a\in X}{1-e^{-a}\over a}$ can
be understood as follows. Each factor ${1-e^{-a}\over a}$ can be expanded in power series.  We select from at most $k$ factors the homogeneous
parts of some degrees $h_i>0$ so that $\sum_ih_i=k$ and evaluate $p$ against the resulting monomial, then we multiply by the remaining factors
and sum over such choices.   Now, if the factors we have chosen correspond to a cocircuit  the evaluation of $p$ on the monomial is 0. If
instead this is not the case we  have still a product $\prod_{a\in A}{1-e^{-a}\over a}$  where the elements $A$ span. Thus, if $\mu\neq 0$ there
is at least one $a\in A$ which does not vanish at $\mu$. But now  clearly ${1-e^{-a} }$  and hence $(1-e^{-a})/ a$  vanishes at $\mu$ and the
Lemma follows. \qed

We go back to the proof of our Theorem. We have shown that in our case,  Poisson summation degenerates  to the computation at 0. Taking a
polynomial in derivatives and computing against a function and then at 0 is just duality thus:
$$p(x) (e^{ x}\prod_{a\in X} {1-e^{- a}\over  a})(0)=\langle  p\,|\,  e^{ x}\prod_{a\in X} {1-e^{- a}\over  a}\rangle=
  \langle  F_{ X}p\,|\,  e^{ x} \rangle =F_{ X}p( x). $$

Since $D(X)$ is stable under derivatives and $F_X$ is clearly invertible both our claims follow.\qed

\subsection{Paving the box}
Take the box  $B(X):=\{\sum_{a\in X} t_aa \,| \, 0\leq t_a\leq 1\}$, which is the support of the box spline.

We start by giving  a nice decomposition of $B(X)$ into suitable parallelepipeds.  In order to present it we need the following:
\begin{lemma}\label{faccio}
If a point $v:=\sum_{a\in X}  v_aa, \ 0\leq v_a\leq 1$   is in the boundary of $B(X)$   the set  $A:=\{x\,|\, 0< v_a< 1\}$ does not span  $V$.
\end{lemma}
\begin{proof}  If $A$ spans let us extract a basis  $b_1,\dots,b_s$ from $A$, take $0<\epsilon$ small. Then the set of points   $v+\sum_{i=1}^s  t_ib_i,\ |t_i|<\epsilon$ is an open ball contained in  $B(X)$.

 \end{proof}   Given $\lambda\in\Lambda$ and a set of linearly
independent vectors $\underline b:=\{b_1,\dots,b_h\}$ from $X$ define:
$$\Pi_{\lambda}(\underline b):=\{\lambda+\sum_{i=1}^ht_i b_i\},\  0\leq t_i\leq 1.$$

\begin{proposition}  $B(X)$  can be paved with parallelepipeds, of the form $\Pi_{\lambda}(\underline b)$ where $\underline b$ runs on the set of all bases
extracted from $X$ and  $\lambda\in \Lambda$ depends on $\underline b$.
\end{proposition}
\proof Suppose that  $X=\{Z,y\}$ and we have paved $B(Z)$   (by induction) so that its boundary is paved by some faces of these parallelepipeds.
Consider
$$ B(Z)_y:=\{p\in B(Z)\,|\, p+ty\notin B(Z),\ \forall t>0\}.$$

We have  a map $\pi:B(Z)_y\times [0,1]\to B(X),\ (p,t)\mapsto p+ty$. We easily see that:
\begin{enumerate}
\item $\pi$ is a homeomorphism to its image. \item $B(Z)\cap  \pi(B(Z)_y\times [0,1])=B(Z)_y$ \item  $B(X)=B(Z)\cup  \pi(B(Z)_y\times [0,1])$
\end{enumerate}

We are now going to pave $B(Z)_y$   by $s-1$ dimensional parallelepipeds $\Pi_{\lambda}(\underline c)$ with $\underline c=\{c_1,\dots,c_{s-1}\}$
running on all  sets of $s-1$ linearly independent vectors in $Z$ which are together with $y$ form a basis.

Consider the set $\mathcal H_y$ of all hyperplanes $H$ generated by subsets of $Z$ and not containing $y$. Given such $H$ take the unique linear
form $\phi_H$ vanishing on $H$ and such that $\langle \phi_H|y\rangle=1$. Set
$$\lambda_{H,y}=\sum_{x\in Z|\phi_H(x)>0}x.$$
We claim that
$$B(X)_y=\cup_{H\in \mathcal H_y}B(Z\cap H)+\lambda_{H,y}$$
and this is a paving. Remark that $\phi_H$ takes its maximum value on $B(Z)$ in the point $\lambda_{H,y}$. Since if $b>0$ and $v\in B(Z\cap
H)+\lambda_{H,y}$, $\phi_H(v+by)=\phi_H(\lambda_{H,y})+b>\phi_H(\lambda_{H,y})$ we get the inclusion $B(Z\cap H)+\lambda_{H,y}\subset B(X)_y$.

To see the converse observe that by Lemma \ref{faccio}, $B(X)_y$ is a union of polytopes of the form $B(Z\cap H)+\mu$ with $H\in \mathcal H_y$
and $\mu=\sum_{a\in A}\varepsilon_a a$ with $A= Z-(Z\cap H)$ and $\varepsilon_a\in\{0,1\}$. Fix $a\in A$. Write $y=h+\phi_H(a)^{-1} a$. Take $v$
in the relative interior of $ B(Z\cap H)$. We have
$$v+ty=(v+th)+(\mu-\varepsilon_aa)+(\varepsilon_a+t\phi_H(a)^{-1})a.$$
If $t$ is sufficiently small then $v+th\in B(Z\cap H)$. Furthermore  if $\varepsilon_a=0$ and $\phi_H(a)>0$ or $\varepsilon_a=1$ and
$\phi_H(a)<0$, $0<\varepsilon_a+t\phi_H(a)^{-1}<1$. Thus this point lies in $B(Z)$ giving a contradiction.

 By induction, each $B(Z\cap H)$ is paved by $s-1$ dimensional parallelepipeds $\Pi_{\lambda}(\underline c)$ with $\underline c=\{c_1,\dots,c_{s-1}\}$
running on all bases of $H$ extracted from $Z\cap H$. Since $\{\underline c,y\}:=\{c_1,\dots,c_{s-1}, y\} $  is a basis of $X$ and all bases
containing $y$ are so obtained, we get the desired paving
$$   \pi(B(Z)_y\times [0,1])=\cup_{\lambda,\underline c} \Pi_{\lambda}(\{\underline c,y\}).$$
From this our claim is immediate.\qed

As a simple application we obtain
\begin{proposition}Let $x_0$ be a regular point. Then
$(B(X)-x_0)\cap \Lambda$ consists of $\delta(X)$  points.
\end{proposition}
It is easy to see that each parallelepiped $\Pi_{\lambda}(\underline b)$ translated by a regular point, intersects $\Lambda$ in $|\det(
\underline b)|$  points. Summing over all parallelepipeds we obtain our claim.\qed

\subsection{Linear independence} We now assume that we are in  the unimodular case. In this case $\delta(X)=d(X)$.
Choose a regular point $x_0$ consider  the $d(X)$ points $P(x_0):=(B(X)-x_0)\cap \Lambda=\{p_1,\dots,p_{d(X)}\}$.
\begin{proposition}
Evaluation of polynomials in the points in $P(x_0)$ establishes  a linear isomorphism between $D(X)$ and the $\mathbb C^{d(X)}$(or $\mathbb
R^{d(X)}$ if we restrict to real polynomials).
\end{proposition}
\proof  Since $\dim(D(X))=d(X)$  it suffices to  prove that, a polynomial $p(x)\in D(X)$ vanishing on these points is identically zero. From the
formula  $p(x)=\sum_{a\in \lambda} B(x+a)p(-a)$ we see that,  for given $x$  the only terms appearing in the sum are the ones where $x+a\in
B(X)$    or $a\in B(X)-x$. Thus for $x$ in a small open neighborhood of the given point    $x_0$ the only terms appearing are the ones in which
$a\in \{p_1,\dots,p_{d(X)}\}$. If by way of contradiction, a polynomial $p\in D(X)$  vanishes on these points  we have that  $p(x)=0$ proving
the Proposition.\qed

Now we can prove the {\it Theorem on linear independence of the translates of the box spline}:
\begin{theorem}
For $X$ unimodular and any, non identically 0, function $f(\lambda)$ on $\Lambda$ we have:
$$\sum_{a\in \Lambda} B(x+a)f(-a)\neq 0.$$\end{theorem}
\proof  Assume $f(a_0)\neq 0$ for some $a_0\in\Lambda$.  We can find  then, a regular point $x_0$, such that $a_0\in \mathcal  P(x_0)$ is one of
the points $p_i$ previously defined.  Thus there is a non zero polynomial $p(x)\in D(X)$ coinciding, on $\mathcal  P(x_0)$  with $f$.   Then
$\sum_{a\in \Lambda} B(x+a)f(-a)=p(x)\neq 0$  on the set $B(X)-x_0$.\qed

Unimodularity is a necessary condition for this theorem in fact one has:

\begin{proposition}  Discrete convolution maps $\nabla(X)$ into $D(X)$ with kernel $E(X)$ (see Formula (\ref{decc})).

\end{proposition}
\begin{proof}  The first statement follows from the following identity, for $Y\subset X$ we have:
$$  D_Y(B_X*_da)=B_{X/Y}*_d \nabla_Ya$$

The  second by the fact that the kernel of the discrete convolution is invariant under translation thus it suffices to verify it for the
functions $e^\phi$ as $\phi$ varies over the points of the arrangement.   Let $Y\subset X$ be a basis of a sub--lattice $\Lambda_0$ on which
$\phi$ is 1. We have $B_X*_d\phi=B_{X/Y}* B_Y*_d\phi$ and $B_Y*_d\phi=0$ since  $\sum_{a\in\Lambda /\Lambda_0}\phi(a)=0$ (the sum is over coset
representatives).

\end{proof}

We apply the previous theory  to the partition function   $\mathcal P_X(a)$.

\begin{proposition}
$$T_X(x)=\sum_{a\in\Lambda}\mathcal  P_X(a)B_X(x-a)$$
\end{proposition}
\proof Compute the Laplace transform $L(\sum_{a\in\Lambda} \mathcal P_X(a)B_X(x-a))$:
$$ =\sum_{a\in\Lambda}\mathcal  P_X(a)e^{-a}L(B_X )=\prod_{a\in X}{1\over  1-e^{-a}} \prod_{a\in X} {1-e^{-a}\over a}=\prod_{a\in X} {1\over a}=LT_X.$$  \qed

In  the unimodular case, where we have  the  linear independence of translates of $B_X$, we recover the results of section  \ref{vtp}.

On  $  D(X)$,   $F_X$ is invertible and its inverse is: $$Q:=\prod_{x\in X} { D_x \over 1-e^{-D_x} } .$$

We have  by definition $QFp=p=FQp$ on $D(X)$. Take a big cell $\mathfrak c$   over which $T_X$ coincides with some polynomial
$p_{\mathfrak c}\in D(X)$. Set
$q_{\mathfrak c}:=Qp$, we have  on $\mathfrak c$

$$T_X(x)=FQp_{\mathfrak c}=\sum_{a\in \Lambda}Qp_{\mathfrak c}(a)B(x-a)$$
since $T_X(x)=\sum_{a\in\Lambda}\mathcal  P_X(a)B_X(x-a)$ we have by  linear independence:

\begin{equation}\label{secfor}Qp_{\mathfrak c}(a)=\mathcal P_X(a).\end{equation}

Notice that $Q$ is like a {\it Todd operator}, its factors can be expanded using the Bernoulli numbers $B_n$ by the defining formula:
$$ { D_x \over 1-e^{-D_x} }=\sum_{k=0}^\infty  {B_n\over n!}(-D_x)^n$$ Formula \ref{secfor} allows us to pass from the formula for the volume to one for the partition function (a special case of formula \ref{pexxpr2}).

\part{Residues}
\bigskip

\section{Residues}
\subsection{Cohomology}

For the computations of volumes and partition functions we want to apply a cohomological method, like the usual method in one complex variable
for computing definite integrals.

With the notations of section \ref{duemod}, we start  from the affine algebraic  variety $\mathcal A_X$, with coordinate  ring  $R_X$. Using the
De Rham's  Theorem     due to    Grothendieck, in order to compute the cohomology of $\mathcal A_X$ with complex coefficients we can use the
algebraic de Rham complex $(\Omega^*_X,d)$ (cf. \cite{gr}). Here $d$ is the usual de Rham differential while for each $0\leq k\leq s$, the algebraic
differential forms of degree $k$ are:
$$\Omega^k_X:=\{\sum_{1\leq i_1<i_2<\dots<i_k\leq s}f_{i_1, \dots,i_k}(x)dx_{i_1}\wedge   \dots \wedge dx_{i_k}\,|\, f_{i_1,  \dots,i_k}(x)\in R_X\} $$ where $\{x_1,\ldots ,x_s\}$ are coordinates in $U$.

Let us look more closely to what happens in degree $s$, the top degree. Define a homomorphism, uniquely determined by the choice of the volume
form $dx_1\wedge\cdots \wedge dx_s$,
$$i_X:\Omega^s_X\to R_X$$
by $i_X(fdx_{1}\wedge   \dots \wedge dx_{s})=f$. We have that $i_X(d(\Omega^{n-1}_X))$ is the space $\partial(R_X)$ spanned by all partial
derivative of elements in $R_X$.

\begin{definition}
Let   $H_X$ be the  space spanned by   the    elements  $M_{\underline b}:=  \prod_{a\in \underline b}a^{-1}$, as $\underline b$ varies among
the bases  extracted  from the list  $X$.\end{definition}

From the computation of the cohomology of $\mathcal A_X$ for which we refer to \cite{te}, we get the following result which also follows
immediately from the discussion of the next section and the expansion in partial fractions.
\begin{theorem}\label{teofon} \begin{enumerate}
\item  We have the decomposition in direct sum:
$$R_X=H_X\oplus \partial(R_X)$$
In particular, using $i_X$, the space  $H_X$ can be identified to the space $H^n( \mathcal A_X)$. \item   The  elements   $M_{\underline b} $,
as   $\underline b$  varies   in the set $\mathcal {NB}(X)$ of no broken bases are a basis  of    $H_X$. \end{enumerate}

\end{theorem}

If we take any basis $\underline b=\{b_1,\ldots b_s\}$ of $V$ we have
$$db_1\wedge\cdots \wedge db_s=det(\underline b)dx_1\wedge\cdots \wedge dx_s$$ (this is gives the definition of $det(\underline b)$).  It is then  clear that
$$ \det(\underline b)\prod_{i=1}^s{1\over b_i}dx_1\wedge\dots\wedge dx_s =d\log(b_1)\wedge\cdots\wedge d\log(b_s).$$

Using $i_X$ we deduce:

\begin{proposition}\label{basette}
The    cohomology classes of the  forms   $$\omega_{\underline b}:=d\log(a_1(x))\wedge\dots\wedge d\log(a_n(x))=d\log(a_1 )\wedge\dots\wedge
d\log(a_n )$$  as   $ a_1 ,\dots,a_n  $  varies   in  $\mathcal {NB}(X)$, form a  basis    of  the cohomology  $H^n( \mathcal A_X)$.
\end{proposition}

\subsection{Local residue\label{primloc}}  From the expansion in partial fractions follows easily   that  $R_X=H_X+ \partial(R_X)$, thus in fact Theorem \ref{teofon}   is a consequence of Proposition \ref{basette}. This can be proved directly  using the following
  local computation.\smallskip

Let $\underline b=\{b_1,\ldots, b_s\}$ be a basis extracted from $X$. To $\underline b$ we associate an injection
$$j_{\underline b}:R_X\to \mathbb C[[u_1,\ldots , u_s]][(u_1\cdots u_s)^{-1}]$$
 ($\mathbb C[[u_1,\ldots , u_s]][(u_1\cdots u_s)^{-1}]$ is the ring of formal Laurent series in the variables $u_1,\ldots ,u_s$), defined by
 $$j_{\underline b}(f(b_1,\dots ,b_s))=f(u_1,u_1u_2,\ldots ,u_1u_2\cdots u_s).$$
 This is well defined since given $a\in X$, if $k=\gamma_{\underline b}(a)$ is the maximum index  such that $a\in\langle b_k,\ldots ,b_s\rangle$, we have
 $$a=\sum_{j=k}^s\alpha_jb_j=u_1\cdots u_k(\alpha_k+\sum_{j=k+1}^s\alpha_j\prod_{i=k+1}^ju_i)$$
 with $\alpha_k\neq 0$ so that $a^{-1}$  can be expanded as a Laurent series. Clearly this map extends at the level of differential forms.

 Given a top differential form $\omega\in \Omega^s_X$ we   define $res_{\underline b}(\omega)$ as the coefficient of $(u_1\cdots u_s)^{-1}du_1\wedge\cdots \wedge u_s$
  in  $j_{\underline b}(\omega)$.

The main technical result from which  all the important results follow is given by

\begin{lemma}\label{duali}Let $\underline b $ and $\underline c$ be two no broken bases extracted from $X$ then
$$res_{\underline b}(\omega_{\underline c})=\delta_{\underline b,\underline c}.$$
\end{lemma}
\begin{proof} First notice that $j_{\underline b}(\omega_{\underline b})=d\log u_1\wedge\cdots \wedge d\log u_s$ so that
$res_{\underline b}(\omega_{\underline b})=1$.

Now notice that if $\underline c\neq \underline b$ there are two distinct elements $c,c'\in \underline c$ such that $\gamma_{\underline b}
(c)=\gamma_{\underline b}(c')$. Indeed since in a non broken basis the fist element is always also the first element in the list $X$,
$\gamma_{\underline b}( c_1)=1$. If $\langle c_2,\ldots c_s\rangle\neq \langle b_2,\ldots b_s\rangle$ then there exists an index $i>1$ with
$\gamma_{\underline b}(c_i)=1$. If $\langle c_2,\ldots c_s\rangle= \langle b_2,\ldots b_s\rangle$ then   both $\underline c-\{c_1\}$ and
$\underline b-\{b_1\}$ are no broken bases of $\langle c_2,\ldots c_s\rangle$ extracted from $X\cap \langle c_2,\ldots c_s\rangle$ and
everything follows by induction.

We than write each $c$ in the form $\prod_{i=1}^ku_i  f$  with $f(0)\neq 0$ so that $d\log c=d\log (\prod_{i=1}^ku_i)+d\log f$.  Expanding the product $j_{\underline b}(\omega_{\underline c})$   we get a linear combination of forms, all terms containing only factors of type $d\log (\prod_{i=1}^ku_i)$  vanish since two elements are repeated,   the others
 are a
product of a closed form by a form  $d\log(f)$ with $f(0)\neq 0$  which is exact, so all  these terms are exact and the residue is 0.

\end{proof}

As a first consequence we deduce:

\begin{theorem}\label{indiplin} The cohomology classes of the forms  $\omega_b$, as $\underline b$ varies among the no broken bases are linearly independent.

For any top differential form $\psi$, denote by $[\psi]$ its cohomology class we have:
\begin{equation}\label{espa} [\psi] = \sum_{\underline b\in \mathcal{NB}(X)} res_{\underline b}(\psi)\omega_{\underline b}\end{equation}

\end{theorem}
\proof

Since clearly the map $j_{\underline b}$ takes exact forms to exact forms we deduce that $res_{\underline b}$ factors through $H^s(\mathcal
A_X)$. In view of this everything follows from Lemma \ref{duali}.

\begin{remark} It is not difficult to  verify that  we get the same cohomology if, instead of  taking   forms  with    coefficients in  $R_X$   we take   coefficients of   type  $f/d^k$ with  $f$    any function holomorphic   around 0.
\end{remark}

We are now going to define  the {\it  total  residue  $Tres$} and   the residues  {\it $res_{\underline b}$.}

With the use of $i_X$ these operators  can be  defined  either on  algebraic  differential forms   of   degree $s$  or   on    functions.
\bigskip

\begin{definition}
Given $f\in R_X$, $Tres( f)$ is  the cohomology class  of   the form  $fdx_1\wedge\dots\wedge dx_n$.
\end{definition}
We can now reformulate formula (\ref{espa}) as: $$ Tres( f)=\sum_{\underline b\in\mathcal {NB}(X)} res_{\underline b}(f) [\omega_{\underline b}].$$
\subsection{Residues and Laplace transform}

Recall that in formula (\ref{JK}) we have introduced the polynomials $p_{\underline b,X} $, in order  to give  an expression of the multivariate
spline $T_X$ .
\begin{theorem}\label{corrare} For every $\underline b\in\mathcal {NB}(X)$, \begin{equation}\label{resiform} p_{\underline b,X}
(-y)= det(\underline b) res_{\underline b}({e^{\langle y|x\rangle} \over d(x)}).\end{equation} with $d(x)=\prod_{a\in X}a$. \end{theorem}
\begin{proof} We begin remarking that formula (\ref{resiform}) makes sense since, if we expand $e^{\langle y|x\rangle}/ d(x)$ with respect to
the variables $y=\{y_1,\ldots ,y_s\}$, we get a power series whose coefficients lie in $R_X$.

In order to    prove    this  Theorem we need some properties of   $Tres$. The first property of   $Tres$, which  follows  from the definition
is  that, given a  function $f$ and  and an index   $1\leq i\leq s,$ we have {\it $Tres(\pd{f}{x_i})=0$}, hence for two  functions $f,g$:
$$  Tres(\pd{f}{x_i}g)=-Tres(f \pd{g}{x_i}).$$ In other words for a polynomial $P$:
\begin{equation}\label{agg}  Tres(P(\partial_x)(f) g)=Tres(f  P(-\partial_x)(g)).\end{equation}

We shall  use the preceding relation    (\ref{agg})    for the function      $f=e^{\langle y|x\rangle}$  for which    we have:
\begin{equation} P(\partial_x)e^{\langle y|x\rangle}= P(y)e^{\langle y|x\rangle},\end{equation}

The second simple property is  that, given a  basis    $\underline b$ extracted from  $X$  and     a  function $f$ regular at    0, we have:
\begin{equation}\label{eva}  Tres( {f\over\prod_{a\in\underline b} a(x) })=f(0)Tres({1\over\prod_{a\in\underline b} a(x)}).\end{equation}
We get
$$ Tres({e^{\langle y|x\rangle} \over d(x)})=Tres(\sum_{\underline b\in\mathcal {NB}(X)}e^{\langle y|x\rangle}p_{\underline b,X} (\partial_x){1\over\prod_{a\in\underline b} a(x)})$$$$=\sum_{\underline b\in\mathcal {NB}(X)}Tres({p_{\underline b,X}(-\partial_x)(e^{\langle y|x\rangle})\over\prod_{a\in\underline b} a(x)})=\sum_{\underline b\in\mathcal {NB}(X)}p_{\underline b,X}(-y)Tres( {1\over\prod_{a\in\underline b} a(x)})$$\begin{equation}\label{aggi}=\sum_{\underline b\in\mathcal {NB}(X)}{1\over  det(\underline b) }p_{\underline b,X}(-y)[\omega_{\underline b}].\end{equation}
From this the theorem   follows.
\end{proof}

\subsection{Partition functions}
 The same method can be applied to partition functions.  In (\cite{dp3}) we have computed the full cohomology
of the toric arrangement, nevertheless this is not strictly necessary for the residue computations which are essentially local.
In the formula (\ref{JKP}) one has contributions localized at a point of the arrangement $e^\phi$ and for a no broken basis
$\underline b\subset X_\phi$.
They can be computed again using the operator $res_{\underline b,\phi}$ in this case.

This means that, we restrict a function or a form to a neighbourhood of $\phi$, use logarithmic coordinates so that the divisors $1-e^{-a}=1$
which appear in this neighbourhood coincide with the linear hyperplanes of the arrangement associated to $X_\phi$ and finally  compute the local
residue at $\underline b$ for this local hyperplane arrangement.

By formula (\ref{ilv}) and Proposition \ref{laide}, we get that the class of the Laplace transform  $\prod_{a\in X}(1-e^{-a})^{-1}$ in the
module $\Pi_X$ of polar parts decomposes into a sum of local factors, corresponding to the summands $\Pi_{X,\phi}$ and then $\Pi_{X,\phi}=\oplus
S[U][ \prod_{a\in\underline b}(a-\langle \phi\,|\, a\rangle  )^{-1}] $ under the logarithm map.   With this isomorphism $j$  we have
$$j(\prod_{a\in X}(1-e^{-a})^{-1})=\sum_{\phi\in P(X)} \sum_{\underline b\in\mathcal{NB}_{X_\phi}}\mathfrak q_{\underline b,X_\phi}(  y)[ \prod_{a\in\underline b}(a-\langle \phi\,|\,  a\rangle  )^{-1}]$$ which by \ref{restr} gives a formula for the partition function.
We claim that we have  \begin{equation}\label{resper} \mathfrak q_{\underline b,X_\phi}( -y) = det(\underline b) res_{\underline
b,\phi}({e^{\langle y|z\rangle} \over \prod_{a\in X} (1-e^{-a(z)-\langle \phi\,|\,  a\rangle})}). \end{equation} In order to prove this, we
start by  making a change of coordinates $x=z+ \langle y|\phi\rangle$ so that we center the point at 0 and we pick out the factor
$e^\phi=e^{\langle y|\phi\rangle}$  and then  showing,  using the language of proposition \ref{laide},  that:
\begin{lemma} $res_{\underline b,\phi}(f)= res_{\underline b,\phi}(j(f))$ for any $f\in M_2$.
\end{lemma}
\begin{proof} Due to the properties of $Tres$ it is enough to prove it on the generators where we see that, around the point $\phi$  we have
$$  \frac {e(\phi)}{\prod_{a\in \underline b} (1-e^{-a}) }=|\det(\underline b)|^{-1}\sum \phi(\lambda)e^{-\lambda }\prod_{a\in \underline b} (1-e^{-a})^{-1}$$ $$=|\det(\underline b)|^{-1}\sum \phi(\lambda)e^{-\lambda} \prod_{a\in \underline b}a^{-1}h(x),\ \qquad h(\phi)=1.$$ Since by definition $\phi(\lambda)=e^{ \langle\lambda\,|\,\phi\rangle}$,  at the point  $\phi$ we have that $$|\det(\underline b)|^{-1}\sum \phi(\lambda)e^{-\langle\lambda\,|\,\phi\rangle} =1.$$
This proves the claim by formula (\ref{eva}).\end{proof}

Now  we go back to formula \ref{resper}.   We have that
\begin{equation}\label{resper2}   res_{\underline b,\phi}(j\big({e^{\langle y|x\rangle} \over \prod_{a\in X} (1-e^{-a(x)})}\big))=res_{\underline b,\phi}(e^{\langle y|x\rangle}\mathfrak q_{\underline b,X_\phi}(  y)[ \prod_{a\in\underline b}(a-\langle \phi\,|\,  a\rangle  )^{-1}]) \end{equation}
Change coordinates, so to center $\phi$ at 0, getting  $x=z+\phi$ and so
$$res_{\underline b,\phi}(e^{\langle y|x\rangle}\mathfrak q_{\underline b,X_\phi}(  y)[ \prod_{a\in\underline b}(a-\langle \phi\,|\,  a\rangle  )^{-1}])=  det(\underline b) ^{-1}   \mathfrak q_{\underline b,X_\phi}(-  y) e^\phi$$\bigskip

 {\bf Summarizing}   These formulas, together with the local computation of residues form an effective algorithm to compute the functions we
 have been studying.

It remains to discuss a last    algorithmic point.

In order to    compute the Jeffry--Kirwan residue, at a given a  point     $p\in C(A)$,  it is necessary to  determine  a  big cell $c$  for
which $p\in\overline c$.

\smallskip

In general, the determination of the big cells is  a very complex  problem. For our computations it suffices much less.

\smallskip
 
Let us take   thus simply a  point     $q$ internal to    $C(A)$ and  not  laying   on  any hyperplane     generated by  $n-1$  vectors of
$X$. This  is not   difficult to do,  and let us consider the segment  $qp$.
\smallskip

This segment  intersects   these   hyperplanes      in a finite number  of   points, thus we can   determine  an $\epsilon$ sufficiently small
for which all the points $tp+(1-t)q, 0<t<\epsilon $ are  regular.
\smallskip

If we take   one of   these points $q_0$ it lays  in  a  cell for which $p$ is  in the closure.

At  this  point,   for every no broken basis, we must verify in simple way  if  $q_0$ lays or not    in the cone generated by the basis.

\section{ Minimal models\label{iresidues}}
\vskip8pt

\subsection{Geometry of residues}

For the moment our definition of  residues $res_{\underline b} \psi$ is  purely algebraic. In fact   its true geometric meaning is based upon a
general definition of the notion of   residue in  several dimensions.\medskip

This section is quite independent  of the rest of the paper and can be used as an introduction to the theory developed in    \cite{dp} and
\cite{dp1}.
\smallskip

The first point to be understood is that, the non linear coordinates $u_i$ used in
  section \ref{primloc}, represent local coordinates around a {\it point at infinity}  of a suitable geometric model of  a completion of the variety $\mathcal A_X$.  In fact we are thinking of models  proper over the space $U\supset \mathcal A_X$ in which the complement of  $\mathcal A_X$ is a divisor with normal crossings. In this respect the local computation  done in  section \ref{primloc}, corresponds to a model in which all the subspaces of the arrangement have been blown up, but there is a subtler model which gives rise to a more intricate combinatorics but possibly to  more efficient computational algorithms, due to its {\it minimality}.

\section{Irreducibles and  nested sets \label{nested}}

\subsection{Irreducibles and decompositions}

The notions  that  we are about to give are of combinatorial nature   (cf. \cite{fei}) but  we develop  them in  the  following context:

As usual let us consider a list   $X:=\{a_1,\dots,a_N\}$ of non zero vectors in $V$ which in this section we assume to be a complex vector
space.

Given  a sublist  $A\subset X$  the list {$\overline A:=X\cap \langle A\rangle$}  will be called the completion of  $A$. Thus  $A$ is  complete
if  $A=\overline A$.

The space of vectors $\phi\in U$ such that $\langle a|\phi\rangle=0$ for every $a\in A$ will be denoted by $A^\perp$. Notice that clearly
$\overline A$ equals to the list of vectors $a\in X$ which vanish on $A^\perp$.

From this we see that  we get a bijection between the complete sublists   of   $X$  and subspaces of the arrangement defined by  $X$.

We give the main
\begin{definition}   Given a     complete set   $A\subset X$,   a  {\bf decomposition }   is a decomposition  $A=A_1\cup  A_2$  in    non empty sets, such  that:  $$\langle A\rangle=\langle A_1\rangle\oplus \langle A_2\rangle. $$
Clearly the two  sets   $A_1,A_2$ are necessarily complete.

We shall say   that  a     complete set   $A$ is   irreducible    if it does not have    a non trivial decomposition.\end{definition}

If $A=A_1\cup A_2$ is   a decomposition   of   a     complete set   and $B\subset A$ is     complete we have $B=B_1\cup B_2,$ where $
B_1=A_1\cap B,\ B_2=A_2\cap B$. Also $\langle B\rangle=\langle B_1\rangle\oplus \langle B_2\rangle, $ and we have
\begin{lemma}  $B=B_1\cup B_2$ is   a decomposition, unless one of the two  sets   is empty.
\end{lemma}

We deduce immediately:

\begin{proposition}\label{contai}
{If  $A=A_1\cup  A_2$  is   a decomposition   and $B\subset A$ is  irreducible, then $B\subset A_1$ or $B\subset A_2$} \end{proposition} From
this get:
\begin{theorem}
Every   set   $A$   can  be decomposed as  $A=A_1\cup A_2\cup\dots \cup A_k$ with the $A_i$ irreducible  and:
$$\langle A\rangle=\langle A_1\rangle\oplus \langle A_2\rangle\oplus\dots\oplus \langle A_k\rangle.  $$

This decomposition is    unique up to   order.\end{theorem} \proof The existence of an irreducible decomposition follows by a simple induction.

Let   $A=B_1\cup B_2\cup\dots \cup B_h$ be    a second decomposition. Proposition \ref{contai} implies that every $B_i$ is   contained in an
$A_j$ and viceversa.

Thus the $A_j$'s and the $B_i$'s    are the same  up to   order.\qed

$A=A_1\cup A_2\cup\dots \cup A_k$ is called the {\it decomposition into irreducibles} of    $A$.

\begin{example}\label{configu}  An interesting example is that of the      configuration space of $s$-ples of point in a line (or the root system $A_{s-1}$). In this case $X=\{z_i-z_j|1\leq i<j\leq s\}$.

In this case, irreducible sets are in bijection with subsets of $\{1,\ldots ,s\}$ with least 2  elements. Indeed give one such subset $S$ it
corresponds to the irreducible $I_S=\{z_j-z_j|\{i,j\}\subset S\}$.

Given a complete set   $C$,  the irreducible decomposition of $C$    corresponds a sequence of disjoint subsets $ S_1,\dots,S_k$    of
$\{1,\ldots ,s\}$ with least 2  elements.
\end{example}

\subsection{Nested sets}
We define now     the basic notion of   {\it   nested set}.\smallskip

We say   that  two  sets   $A,B$ are {\it comparable} if  one is  contained in the other.

\begin{definition}  We shall say that   a family $\mathcal S$ of   irreducibles  $A_i$ is  {\sl nested} if given   elements  $A_{i_1},\dots, A_{i_h}\in \mathcal S$ mutually incomparable we have that
$C:=A_1\cup A_2\cup\dots\cup A_i$ is complete and  $C:=A_1\cup A_2\cup\dots\cup A_i$  is its decomposition into irreducibles.
\end{definition}
\begin{remark}\label{oss}  If $A_1,\dots,A_k$ is  nested we have that  $\cup_iA_i$ is  complete. In fact this    union can be obtained   taking the
maximal elements,  that  are necessarily non  comparable, and  then applying  the  definition of   nested.\end{remark}\smallskip

We are in  particular    interested  in   {\it maximal nested sets}, which we  denote  by    MNS.

Nested sets can be inductively constructed combining the following inductive procedures.

We can use two simple inductive ways to construct  nested sets. The proof is left to the reader:

 \begin{enumerate} \item Suppose we are given a nested set  $\mathcal S$,
an minimal element   $A\in\mathcal S$ and a nested  set $\mathcal P$ whose elements are contained in $A$. Then we have that $\mathcal
S\cup\mathcal P$  is nested. \item Suppose we are given a nested set  $\mathcal S$, a    complete set  $A$ containing each element of $\mathcal
S$. Then if
 $A=A_1\cup \dots \cup A_k$  is the decomposition of $A$  into irreducibles,
 $\mathcal S\cup \{A_1, \dots, A_k\}$ is nested.
\end{enumerate}

In the case of Example \ref{configu} one can easily verify that nested sets correspond to families of subsets in $\{1,\ldots ,s\}$ each
containing at least two elements and such that  any  two subsets of  the family  are either comparable or disjoint.
\smallskip

\begin{theorem}\label{mns} Assume that $\langle X\rangle=V$.
Let  $\mathcal S:=\{A_1,\dots,A_k\}$ be a   MNS in $X$.

Given $A\in \mathcal S$ let  $B_1,\dots, B_r$  be the  elements  of   $\mathcal S$ contained  properly in $A$, and maximal with  this property.

\begin{enumerate}
\item    $C:=B_1\cup \dots\cup B_r$ is  complete and decomposed by the $B_i$.

\item  $\dim\langle A\rangle=\dim\langle C\rangle+1$.

\item   $k=\dim(V)$.\end{enumerate}
\end{theorem}

\proof  (1) is  the  definition of   nested set  since the   $B_i$, being maximal, are necessarily non  comparable.

\smallskip

(2) Let us consider $\langle C\rangle=\oplus_{i=1}^r \langle B_i\rangle\subset \langle A\rangle$.

\noindent  $\langle C\rangle\neq  \langle A\rangle$ otherwise, since $C$ is complete, by  the definition of   nested set,  we must  have  $A=C$.
This is absurd since  $A$ is  irreducible and     the $B_i$'s are properly contained  in $A$.

Therefore there exists an  element  $a\in A\,|\,a\notin \langle C\rangle$. Let us denote by $A':=X\cap\langle C,a\rangle$. We have $C\subsetneq A'\subset A$.
We claim that  $A=A'$.

Otherwise, as  one can  easily see, adding all the   irreducibles that  decompose $A'$ to the family  $\mathcal S$, we obtain    a  nested
family  that  contains properly $\mathcal S$. This contradicts  the maximality  of $\mathcal S$.
 Clearly $A=A'$ implies that $ \langle A\rangle= \langle C\rangle\oplus \mathbb Ca$ and thus $\dim\langle A\rangle=\dim\langle
C\rangle+1$.
\smallskip

(3)   We proceed by induction  on $s=\dim(V)$.

If  $s=1$    there is  nothing to prove, there is  a  unique   set    complete  and     irreducible namely $X$.

Let $s>1$. Decompose   $X=\cup_{i=1}^hX_h$ into irreducibles.

We have that   a  MNS in $X$ is  the union of   MNS in each $X_i$. Then  $s=\dim(\langle X\rangle)=\sum_{i=1}^h\dim(\langle X_h\rangle)$.

Thus we can assume that   $X$ is irreducible. In this case we have that  $X\in \mathcal S$ for every MNS $ \mathcal S$.

Let $B_1,\dots, B_s$  be  the  elements  of   $\mathcal S$  properly contained   in $X$ and maximal with  this property.

The set   $\mathcal S$ consists  of   $X$  and of  the subsets $\mathcal S_i:=\{A\in \mathcal S\,|\, A\subset B_i\}.$

Clearly  $\mathcal S_i$ is   a MNS  relative to the set   $B_i$ (otherwise we could add  an  element   to    $\mathcal S_i$  and  to
$\mathcal S$ contradicting the maximality of    $\mathcal S$).

By   induction   $\mathcal S_i$ has $\dim \langle B_i\rangle$  elements  and thus by  (2) the claim follows.  \qed

In the example   of  configuration spaces $X$ spans the hyperplane in $\mathbb C^s$ where the sum of coordinates equals to zero. We have
described irreducibles via the corresponding subsets in $\{1,\ldots ,s\}$. Under this correspondence, one can see easily that a maximal nested
set $ \mathcal S$ is formed by $s-1$ elements  and for any $A\in \mathcal S$ with $a>2$ elements either $A$ contains a unique maximal element
$B\in\mathcal S$ with necessarily $a-1$ elements or exactly two maximal elements $B_1,B_2\in\mathcal S$ with $A=B_1\cup B_2$.

We can   present  such  a MNS in a convenient  way as  a     planar binary  rooted tree with $s$ leaves labelled by $\{1,\ldots ,s\}$. Every
internal  vertex  of the the corresponds   to the set  of  numbers that  appear on its    leaves.

For example the graph

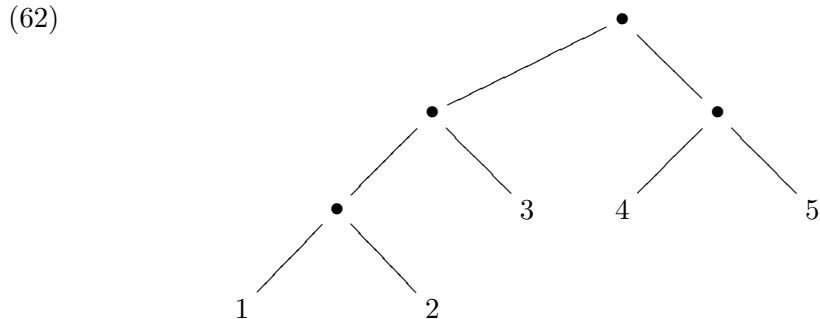
\begin{equation}\label{example}
{\qquad \qquad \xymatrix{&&&& \bullet  \ar@{-}[dr]\ar@{-}[dll] &&\\&&\bullet  \ar@{-}[dr]\ar@{-}[dl]& & &  \bullet \ar@{-}[dr]\ar@{-}[dl]&&\\&
\bullet  \ar@{-}[dr]\ar@{-}[dl] &&3&4& &5 \\ 1& &2  & &  &&}}
\end{equation}
\bigskip

represents the MNS \ \  $ \{1,2\}, \{1,2,3\},\{4,5\},\{1,2,3,4,5\}$.

Given  a MNS   $\mathcal S$   let us define a  map:
$$p_{\mathcal S}:X \to \mathcal S$$ as  follows.
Since $\cup_{A\in \mathcal S}A=X$ every element $a\in X$ lies in at least an $A\in \mathcal S$. Also if an   element  $a\in X$  appears  in two
elements $A,B\in \mathcal S$ the two elements must be necessarily comparable, thus there exists a minimum among the two. It follows that for any
element $a\in X$ there is a minimum element $p_{\mathcal S}(a)\in \mathcal S$.

Now  a new   definition:

\begin{definition}
We shall say that  a  basis   {$\underline b:=\{a_1,\dots,a_s\}\subset X$} of   $\mathbb C^s$, is  {\it  adapted} to the MNS    $\mathcal S$ if
the map $a_i\mapsto p_{\mathcal S}(a_i)$ is a  bijection.
\end{definition}   Such a basis     always exists. It suffices to take, as  in the proof   of   \ref{mns}, for every $A\in
 \mathcal S$ an  element  $a\in A-\cup_iB_i$,  where  the $B_i$ are the    elements  of      $\mathcal S$ properly contained  in $A$.
\medskip

Given any basis  $\underline b:=\{a_1,\dots,a_s\}\subset X$, we  shall build  a MNS  $\mathcal S_{\underline b}$ to    which it is  adapted, in
the   following way. Consider for any $1\leq i\leq s$  the    complete set  $A_i:=X\cap\langle\{a_i,\dots,a_s\}\
\rangle=\overline{\{a_i,\dots,a_s\}}$. Clearly $A_1=X\supset A_2\supset \dots\supset A_s$.

  For    each $i$   consider all the   irreducibles in  the decomposition of   $A_i$. Clearly for different $i$ we can obtain also  several
times the same irreducible, in any case  we have:

\begin{theorem}
The family  {\it $\mathcal S_{\underline b}$} of   all the (distinct)   irreducibles  that  appear  in the decompositions  of the  sets   $A_i$
form  a MNS  to    which the basis  $\underline b$ is  adapted.
\end{theorem}
\begin{proof}
By induction.  Decompose   $X=A_1=B_1\cup B_2\cup\dots\cup B_k$ into irreducibles, by construction:

$$s=\dim\langle A_1 \rangle=\sum_{i=1}^k \dim\langle B_i\rangle.$$

We have that  $A_2=(A_2\cap B_1)\cup (A_2\cap B_2)\cup\dots\cup ( A_2\cap B_k)$ is   a decomposition   of   $A_2$, not  necessarily into
irreducibles.\smallskip

Since $\dim\langle A_2 \rangle= s-1$ we have:
$$s-1=\dim\langle A_2 \rangle=\sum_{i=1}^k \dim\langle A_2\cap B_i\rangle.$$

Therefore   $\dim\langle A_2\cap B_i\rangle<\dim\langle   B_i\rangle$ for exactly one index $i_0$. In other words  we must  have that $A_2\cap
B_i=B_i$ for all the   $i\neq i_0$. For such an index  necessarily $a_1\in B_{i_0}$.

By induction,  the family of   all the     (distinct) irreducibles   that  appear  in the decompositions  of the  sets   $A_i,\ i\geq 2$ form  a
MNS  for $\langle A_2\rangle$, with adapted   basis  $\{a_2,\dots,a_n\}$. To     this   set   we must thus only  add  $B_{i_0}$ in order to
obtain $\mathcal S_{\underline b}$. Thus  $\mathcal S_{\underline b}$ is a nested set with  $s$  elements,  hence  maximal and  the basis
$\underline b$ is adapted.
\end{proof}
\begin{remark} One can easily verify that, conversely, every MNS $\mathcal S$ is of the form $\mathcal S_{\underline b}$ for some adapted basis.
\end{remark}

\subsection{Non linear coordinates} Now we pass to the fundamental   geometric   construction.

Given  a MNS    $\mathcal S$  and     a  basis     $\underline b:=\{a_1,\dots,a_n\}  $   adapted to        $\mathcal S$ let us consider the
$a_i$ as   a {\it system of      linear coordinates     on  $U$.}\smallskip

If $p_{\mathcal S}(a_i)=A$ we denote also $a_i:=a_A$. We build now new     coordinates $z_A,\ A\in \mathcal S$  using the  monomial
expressions:\smallskip

\begin{equation}\label{mon}a_A:=\prod_{B\in \mathcal S,\ A\subseteq B} z_B.\end{equation}

Given  $  A \in \mathcal S$,  let  $ \mathcal S_A:=\{B\subseteq A,\ B\in \mathcal S\}$. Clearly $ \mathcal S_A$ is a MNS for $A$ (in place of
$X$)and the   elements  $a_B$ with  $B\in   \mathcal S_A$ form a  basis of $A$,    adapted to $\mathcal S_A$.

Take $a\in X$ with   $p_{\mathcal S}(a)=A$ and write $a=\sum_{B\in   \mathcal S_A} c_Ba_B,\ c_B\in\mathbb C$. Let us substitute now the
expressions (\ref{mon}):
$$a_A:=\sum_{B\in   \mathcal S_A} c_B \prod_{C\in \mathcal S,\ B\subseteq C} z_C=$$
\begin{equation}\label{monn}\prod_{B\in \mathcal S,\ A\subset B}z_B(c_A+\sum _{B\in   \mathcal S_A, B\neq A} c_B \prod_{C\in \mathcal S,\ B\subseteq C\subsetneq A} z_C).\end{equation}
Since $A$ is  the minimum   set   of     $\mathcal S$  containing $a$ we must  have {\it $c_A\neq 0$}.

At this point we can proceed as in Section \ref{primloc}  and define an embedding of $R_X$ into the ring of Laurent series in the variables
$z_A$, $A\in \mathcal S$. Thus for a top form $\psi$ we can define a local residue  $res_{\mathcal S,\underline b}\psi$ (or $res_{\mathcal
S}\psi$ if $\underline b$ is clear from the context).

\subsection{Proper nested sets}

We are going to tie the concepts  of the previous section to    that of   no broken bases. Our goal is to define a bijection between $\mathcal
{NB}(X)$ and a suitable family of MNS which we shall call {\it proper nested sets}. We need some preliminary steps.

Assume thus that  the basis   $\underline b:=\{a_1,\dots,a_n\}\subset X$ is no broken, we get: \begin{lemma} $a_i$ is  the minimum   element  of
$p_{\mathcal S_{\underline b}}(a_i)$ for every $i$.
\end{lemma}
\begin{proof}  Let $A:=p_{\mathcal S_{\underline b}}(a_i)$, by the    definition of   $S_{\underline b}$,    we must  have that  $A$ decomposes one of the  sets   $A_k=\langle a_k,\dots,a_n\rangle \cap X$. Necessarily it must  be  $k\leq i$.

On the other hand  $a_i$ belongs  to   one of the irreducibles of   $A_i$, that  therefore is  contained in each irreducible $B$ of   $A_k,\
k<i$  that  contains $a_i$.

It follows that $A$ must be one of the irreducibles decomposing  $A_i$. By definition of   no broken basis,  $a_i$ is  the minimum   element  of
$A_i=\langle a_i,\dots,a_n\rangle \cap X$ hence also the minimum element of $A$.\end{proof}

This property suggests us to    define:
\begin{definition}
A MNS $\mathcal S$ is said to be {proper}, if the  elements

$a_S:=\min a,\ a\in S\,|\, S\in \mathcal S$, form a  basis.
\end{definition}
\begin{lemma}\label{oss2}  If  $\mathcal S$  is   {proper},  $p_{\mathcal S}(a_S)=S$. Thus the    elements  $a_S:=\min a,\ a\in S\,|\, S\in \mathcal S$ form a basis
adapted to   $\mathcal S$.\end{lemma}
\begin{proof}  Let   $U\in \mathcal S$ with  $a_S\in U$.
If  $U\subset S$, then  $a_S=\min a\in U$  and thus $a_S=a_U$.

Since the    elements   $a_S$ are a  basis,     this  implies that  $S=U$.\end{proof}

If $\mathcal S$  is  proper, we order its  subsets $S_1,\dots,S_n$ using the increasing order of the  elements  $a_S$. We  then set
$a_i:=a_{S_i}$, and we have:

\begin{theorem}
(1)  The  basis  $\underline b:=\{a_1,\dots,a_n\}$ is  no broken.

(2) In  this  way we establish  a  1-1 correspondence  between no broken   bases  and proper MNS.\end{theorem}  \begin{proof}

(1) From the Remark \ref{oss}, setting $A_i:=\cup_{j\geq i}S_j$, we have that  $A_i$ is  complete and decomposed by the $S_j$ which are maximal.

Clearly, by definition,  $a_i$ is  the minimum   of    $A_i$. It suffices to prove   that $A_i=\langle a_i,\dots,a_s\rangle\cap X$. Hence that
$\langle A_i\rangle=\langle a_i,\dots,a_s\rangle $  since $A_i$ is complete.

We prove it  by induction. If $i=1$ the maximality of $\mathcal S$ implies that  $A_1=X$ so $a_1$ is  the minimum   element  in  $X$. Now
$a_1\notin A_2$  so $A_2\neq X$ and, since $A_2$ is complete $\dim(\langle A_2\rangle)<s$.

Clearly  $\langle A_2\rangle\supset \langle a_2,\dots,a_s\rangle$ so $\langle A_2\rangle=\langle a_2,\dots,a_s\rangle$.

At this point it follows that, by induction,   $\langle a_2,\dots,a_s\rangle$  is a no broken basis  in $A_2$ and thus, since  $a_1$ is minimum
in $X$,    $\langle a_1, a_2,\dots,a_s\rangle$  is a no broken basis.

\smallskip

(2) From the proof, it follows that  the two constructions, of  the MNS associated to   a  no broken basis  and of  the  no broken basis
associated to a proper MNS, are    inverse  of each other and thus the 1-1 correspondence  is  established. \end{proof}

\section{Residues and cycles}

We can  now complete our analysis computing the residues. Recall that, if   $\underline b:=\{b_1,\dots,b_n\}\subset X$ is  a  basis,  we denote
with:
$$\omega_{\underline b}:=d\log(b_1)\wedge \dots\wedge d\log(b_n)$$
\vspace{-10pt}
\begin{theorem}\label{residues} Given   two  no broken  bases  $ {\underline b},  {\underline c}$ we have:$$ res_{\mathcal S_{\underline b}}\omega_{\underline c}=\begin{cases}1\quad\text{if}\quad {\underline b}= {\underline c}\\0\quad\text{if}\quad {\underline b}\neq {\underline c}\end{cases}$$
\end{theorem}
\begin{proof}  We prove first that  $res_{\mathcal S_{\underline b}} \omega_{\underline b}=1$. \smallskip

By definition $b_i=\prod_{A_i\subset B} z_B$ hence $d\log(b_i)=\sum_{A_i\subset B}d\log( z_B).$

When we expand the product we get   a sum  of   products  of   type  $d\log( z_{B_1})\wedge d\log( z_{B_2})\wedge\dots\wedge d\log( z_{B_n})$
with  $A_i\subset  B_i$.

Now, the  unique  non decreasing   and     injective map, of   $\mathcal S_{\underline b}$ in itself  is  the identity.

Therefore in this sum,  all the monomials vanish except for the monomial  $d\log( z_{1})\wedge d\log( z_{2})\wedge\dots\wedge d\log( z_{n})$,
that  has residue 1.

Let us pass now  to the second    case  ${\underline b}\neq {\underline c}$. This follows immediately   from the following Lemma.\end{proof}

\begin{lemma} 1. If  ${\underline b}\neq {\underline c}$, the basis  $ {\underline c}$ is not   adapted to   $\mathcal S_{\underline b}$.
\smallskip

2. If  a  basis    $ {\underline c}$ is not   adapted to  $\mathcal S=\{S_1,\ldots,S_n\} $,     a MNS,  we have $res_{\mathcal S }\omega_
{\underline c}=0$.

\end{lemma}
\begin{proof}     (1) Let  $\underline c=\{c_{ 1},\ldots,c_{ n}\}$ be adapted to    $\mathcal S_{\underline b}$, we want prove   that  we have ${\underline b}= {\underline c}$.
We know that   $c_1=a_1=b_1$, is  the   minimum element    of   $X$.

Let  $A$ be the irreducible component  of   $X$   containing $a_1$. This,   by  definition, is  an  element  of   $\mathcal S_{\underline b}$.
We have   thus $A=S_1$. We claim that $p_{ \mathcal S_{\underline b}}(a_1)=A$.

This follows from the fact that  $a_1=b_1$ and   $ \mathcal S_{\underline b}$ is  proper.   \smallskip

Set $X':=S_2\cup \dots\cup S_n$, $X'$ is    complete. The set $S_2, \dots, S_n$ coincides with  the  proper MNS   $\mathcal S_{\underline b'}$
associated to the no broken basis $\underline b':=\{b_2,\dots,b_n\}$ of   $\langle X'\rangle $.\smallskip

Moreover clearly,  $\underline c'=\{c_{ 2},\ldots,c_{ n}\}$ is   adapted to    $\mathcal S_{\underline b'}$. Therefore ${\underline b'}=
{\underline c'}$ by induction. Hence ${\underline b}= {\underline c}$.

Using the first part, the proof of (2) follows the same lines as the proof of Lemma \ref{duali} and we leave it to the reader.\end{proof}

\begin{remark} \label{concl}  We end this section pointig out that by what we have proved it follows that for any no broken basis $\underline b$ we
have that
$$res_{\underline b}=res_{\mathcal S_{\underline b}}$$
with $res_{\underline b}$ defined in Section \ref{primloc}.
\end{remark}
The advantage of this new definition is that one sometimes uses monomial transformations of smaller degree and this could provide more efficient
algorithms.

\subsection{A minimal model}   Although we do not use it explicitly it may useful to understand the origin of the non linear coordinates we have been using and the entire theory of irreducibles and nested sets that we have built in
\cite{dp}.

Start from  a family of hyperplanes  in $U=V^*$, given by a list $X\subset V$  of linear equations, with complement
$\mathcal A_X\subset U$.

In \cite{dp} we construct  a  minimal smooth variety  $Z_X$ containing $\mathcal A_X$ as an open set with complement a normal crossings divisor,
plus a proper map $\pi:Z_X\to U $ extending the identity of  $\mathcal A_X$.

The smooth irreducible components of the boundary of  $\mathcal A_X$ are indexed by the {\it irreducible subsets}. To describe the intersection pattern between these divisors, in   \cite{dp} we
developed the general theory of nested sets.

 Maximal nested sets  correspond to special points at infinity, intersections of these boundary
divisors. In the papers \cite{S} and  \cite{SV},  implicitly the authors use the points at infinity coming from complete flags which correspond,
in the philosophy of  \cite{dp}, to a {\it maximal}   model with normal crossings. It is thus not a surprise that by passing from a maximal to a
minimal model the combinatorics gets simplified and the constructions become more canonical.\smallskip

Let us recall without proofs the main construction of  \cite{dp}.

For each  irreducible $S\subset X $ we have an orthogonal  subspace $S^{\perp}\subset U $ where
$S^{\perp}=\{a\in U \,|\, x(a)=0,\ \forall x\in S\}.$

From the collections of   projective spaces $\mathbb P(U /S^{\perp})$ of lines in $U /S^{\perp}$ we deduce a map $i:\mathcal A_\Delta\to U^*\times_{S\in\mathcal
I}\mathbb P(U /S^{\perp})$. Set $Z_X$  equal to the closure of the image  $i(\mathcal A_X)$ in this product. In \cite{dp} we have seen that
$Z_X$ is a smooth variety containing a copy of $\mathcal A_X$ and the complement of $\mathcal A_X$ in $Z_X$ is a union of smooth irreducible
divisors $D_S$, having transversal intersection, indexed by the elements   $S\in \mathcal I$.

Still in \cite{dp} we showed that a family $D_{S_i}$ of divisors indexed by irreducibles $S_i$ has non empty intersection (which is then smooth
irreducible) if and only if the family is {\it nested}.  In particular a maximal nested set $\mathcal N$ identifies a special {\it point at
infinity} $p_{\mathcal N}$, intersection of the $s=\dim(U)$ divisors corresponding to the irreducibles in $\mathcal N$. The non linear coordinates are
indeed coordinates in a local chart around $p_{\mathcal N}$ in which $p_{\mathcal N}=0$  and the boundary divisors are given by the vanishing of
the $s$ coordinates.

\section{Final considerations}

As we have seen, the actual problem of computing explicitly the functions which we have discussed can be approached through several different ways. We have proposed three approaches,   by expansion into partial fractions, by solving a system of linear equations interpreting the defining differential equations and finally with the method of residues. It is not
really clear to us which is the most efficient.

A lot of computer work, aiming at computing Clebsch--Gordan coefficients using these methods
has been done by several authors using  the method of residues   but whether this is really the fastest algorithm is yet unclear.

Due to limited space we have   discussed few examples.  For our motivations  one of the most interesting examples is the set $X$ of
positive roots of a given root system.
There are still many things to be uncovered in this case.


\begin{thebibliography}{99}




\bibitem{AS}{A.A.  Akopyan; A.A. Saakyan, }{\it A system of
differential equations that is related to the polynomial class of  translates of a box spline.} (Russian) Mat. Zametki 44 (1988), no. 6,
705--724, 861; translation in Math. Notes 44 (1988), no. 5-6,  865--878 (1989) \newline






\bibitem{WV}{ W. Baldoni-Silva, M. Vergne, }{\it Residues formulae for
volumes and Ehrhart polynomials of convex polytopes,}{   preprint math.CO/0103097}\newline




\bibitem{Be}{ E.T. Bell, }{\it Interpolated denumerants and Lambert
series}{ , Amer. J. Math. 65 (1943), 382--386.}\newline


\bibitem{BV1}{ M. Brion, M. Vergne, }{\it Residue formulae, vector
partition functions and lattice points in rational polytopes.} { J.A.M.S. (v. 10) 4 (1997),  797--833.}\newline


\bibitem{BV}{ M. Brion, M. Vergne, }{\it Arrangement  of hyperplanes. I.
Rational functions and Jeffrey-Kirwan residue.}{ Ann. Sci. \'Ecole Norm. Sup. (4) 32 (1999), no. 5, 715--741.}\newline


\bibitem{BV6}{ M. Brion, M. Vergne, }{\it Arrangement of hyperplanes.
II. The Szenes formula and Eisenstein series. } Duke Math. J. 103  (2000), no. 2, 279--302. \newline


\bibitem{BV2}{ M. Brion, M. Vergne, }{\it An equivariant Riemann--Roch theorem for complete, simplicial toric varieties.} J. Reine  Angew. Math. 482 (1997), 67--92. \newline


\bibitem{BV3}{ M. Brion, M. Vergne, }{\it Lattice points in simple
polytopes.} J. Amer. Math. Soc. 10 (1997), no. 2, 371--392. \newline


\bibitem{BV4}{ M. Brion, M. Vergne, }{\it Une formule d'Euler--Maclaurin pour les polytopes convexes rationnels. } C. R. Acad. Sci.  Paris S\'er. I Math. 322 (1996), no. 4, 317--320. \newline


\bibitem{BV5}{ M. Brion, M. Vergne, }{\it Une formule d'Euler--Maclaurin pour les fonctions de partition. } C. R. Acad. Sci. Paris  S\'er. I Math. 322 (1996), no. 3, 217--220. \newline

\bibitem{Br}{T.  Brylawski, }{\it  The broken-circuit complex.} Trans. Amer. Math. Soc. 234 (1977), no. 2, 417--433.\newline




\bibitem{Co}{S.C.
Coutinho, } { \it A primer of algebraic $D$-modules. }{ London Mathematical Society Student Texts, 33. Cambridge University Press, Cambridge,
1995. }\newline





\bibitem{DM1}{
W. Dahmen, C.Micchelli, }{\it Translates of multivariate splines. }{ Linear Algebra appl. (1983), no. 52, 217--234.}\newline




\bibitem{DM2}{ W. Dahmen, C.Micchelli, }{\it
On the solution of certain systems of partial difference equations  and linear independence of the translates of box splines }{ Trans. Amer.
Math. Soc. 292 (1985), no. 2, 305--320.}\newline



\bibitem{DM3}{
W. Dahmen, C.Micchelli, }{\it The number of solutions to linear Diophantine equations and  multivariate splines. }{ Trans. Amer. Math. Soc. 308
(1988), no. 2, 509--532.}\newline


\bibitem{DM4}{
W. Dahmen, C.Micchelli, }{\it On multivariate $E-$splines. }{ Advances in Math. v. 76, N. 1, (1989) pp. 33-93.}\newline

\bibitem{dhr0}{C. De Boor, K. H\"ollig,  }{\it B-splines from
parallelepipeds.}{ J. Analyse Math,. 42 (1982) pp. 99--115. }\newline



\bibitem{dhr}{C. De Boor, K. H\"ollig, S. Riemenschneider, }{\it Box
splines.}{ Applied Mathematical Sciences 98 (1993). }\newline

\bibitem{dp}{C. De Concini, C. Procesi, }{\it Wonderful Models of
subspace arrangements.}{ Selecta Math. (N.S.) 1 (1995), no. 3, 459--494. }\newline

\bibitem{dp1}{C. De Concini, C. Procesi, }{\it Nested sets and Jeffrey
Kirwan residues.}{in: Geometric methods in Algebra and Number Theory,  Bogomolov, F., Tschinkel. Y., Eds., vol. 235, 139--150. }\newline

\bibitem{dp2}{ C. De Concini, C. Procesi, }{\it On the geometry of
graph arrangements.} {preprint math.CO/0412130. }\newline


\bibitem{dp3}{C. De Concini, C. Procesi, }{\it Toric arrangements.}
{preprint math.AG/0505351. }\newline




\bibitem{DR}{N. Dyn;   A. Ron, }{\it A. Local approximation by certain
spaces of exponential polynomials, approximation order of exponential  box splines, and related interpolation problems.} Trans. Amer. Math.
Soc. 319 (1990), no. 1, 381--403.\newline

\bibitem{Er} {E. Ehrhart, }{\it Polyn\~ omes arithm\'etiques et m\'ethode des poly\`edres en combinatoire, }Birkh\"auser, Basel, 1977\newline


\bibitem{fei}{ E. M. Feichtner, B. Sturmfels, }{\it Matroid polytopes,
nested sets and Bergman fans,}{ preprint CO/0411260.}\newline

\bibitem{FKT}{S. Fomine; A. Kolmogorov; V.M. Tihomirov, }{\it   El\'ements de la th\'eorie des fonctions et de l'analyse fonctionnelle.} (French) Avec un compl\'ement sur les alg\'ebres de Banach, par V. M. Tikhomirov. Traduit du russe par Michel Dragnev. \'Editions Mir, Moscow, 1974. 536 pp.
\newline

\bibitem{gr}{ A. Grothendieck, }{\it  On the de Rham cohomology of
algebraic varieties.}{ Inst. Hautes   \'Etudes Sci. Publ. Math. No. 29 (1966), 95--103.}\newline




\bibitem{Ha}{M. Haiman, }{\it  Hilbert schemes, polygraphs and the
Macdonald positivity conjecture.} J. Amer. Math. Soc. 14 (2001), no.  4, 941--1006 \newline

\bibitem{JKi} {L.C. Jeffrey; F.C.,  Kirwan, } {\it Localization for nonabelian group actions.} Topology 34 (1995), no. 2, 291--327. \newline

\bibitem{J}{ Jia, Rong Qing, }{\it Translation-invariant subspaces and
dual bases for box splines. }(Chinese) Chinese Ann. Math. Ser. A 11  (1990), no. 6, 733--743. \newline






\bibitem{KP}{A.G. Khovanski\u\i; A.V. Pukhlikov, }{\it The Riemann-Roch theorem for integrals and sums of quasipolynomials on virtual polytopes. }(Russian) Algebra i Analiz 4 (1992), no. 4, 188--216; translation in St. Petersburg Math. J. 4 (1993), no. 4, 789--812
\newline



\bibitem{te}{P. Orlik, H. Terao,  }{\it Arrangements of   hyperplanes.}{  Grundlehren der Mathematischen Wissenschaften  300.
Springer-Verlag, Berlin, 1992.}\newline




\bibitem{ron}{A. Ron,  }{\it Exponential box splines.}{   Constr. Approx. 4 (1988), pp. 357-378}\newline
\bibitem{S}{A. Szenes, }{\it
Iterated residues and multiple Bernoulli polynomials. }{ Int. Math. Res. Not. (1998), No.18, 937-956.}\newline

\bibitem{S1}{A. Szenes, }{\it
Residue theorem for rational trigonometric sums and Verlinde's  formula. }{Duke Math. J. 118 (2003), no. 2, 189--227.}\newline

\bibitem{SV1}{ A. Szenes, M. Vergne, }{\it Residue formulae for
vector partitions and Euler-Maclaurin sums}{ preprint, CO/0202253.}
\newline

\bibitem{SV}{ A. Szenes, M. Vergne, }{\it Toric reduction and a
conjecture of Batyrev and Materov }{ preprint, AT/0306311.}
\newline

\bibitem{Ve}{M.  Vergne, }{\it   Residue formulae for Verlinde sums, and for number of integral points in convex rational polytopes.} European women in mathematics (Malta, 2001), 225--285, World Sci. Publishing, River Edge, NJ, 2003.
\newline



\bibitem{Yo}{K. Yosida, }{\it Functional analysis. Sixth edition. }
Grundlehren der Mathematischen Wissenschaften, 123. Springer-Verlag,  Berlin-New York, 1980. xii+501 pp..\newline




\bibitem{W1}{H. Whitney, } {\it  A logical expansion  in mathematics,} Bull. A.M.S., 38  (1932), pp. 572--579.\newline



\bibitem{W2}{H. Whitney, }{\it   On the abstract properties of linear  dependence,} Amer. J.  Math., 57  (1935), pp.  509--533.\newline


\bibitem{Wi}{H. S. Wilf, } {\it Which polynomials are chromatic?,} Colloquio Internazionale sulle Teorie Combinatorie  (Roma, 1973), Tomo I, pp. 247--256, Accad. Naz. Lincei, Rome, 1976.


\end{thebibliography}
\end{document}